\documentclass{amsart}
\usepackage{amsfonts,amssymb,stmaryrd,amscd,amsmath,latexsym,amsbsy}
\numberwithin{equation}{section}
\theoremstyle{plain}
\newtheorem{thm}{Theorem}[section]

\newtheorem{prop}[thm]{Proposition}

\newtheorem{defi}[thm]{Definition}
\newtheorem{lem}[thm]{Lemma}
\newtheorem{cor}[thm]{Corollary}
\newtheorem{eg}[thm]{{Example}}

\theoremstyle{remark}
\newtheorem{rema}[thm]{Remark}
\title{Macdonald difference operators and Harish-Chandra series}
\author{Gail Letzter}
\address{Mathematics Department, Virginia Tech, Blacksburg, VA
24061-0123, USA}
\email{letzter@math.vt.edu}
\author{Jasper V. Stokman}
\address{KdV Institute for Mathematics, Universiteit van Amsterdam,
Plantage Muidergracht 24, 1018 TV Amsterdam, The Netherlands.}
\email{jstokman@science.uva.nl}
\begin{document}
\begin{abstract}
We analyze the centralizer of the Macdonald difference 
operator in an appropriate
algebra of Weyl group invariant difference operators. We show that
it coincides with Cherednik's commuting algebra of difference
operators via an analog of the Harish-Chandra isomorphism.
Analogs of Harish-Chandra series are defined and
realized as solutions to the system of basic hypergeometric 
difference equations associated to the
centralizer algebra. These Harish-Chandra
series are then related to both Macdonald polynomials and Chalykh's
Baker-Akhiezer functions.
\end{abstract}
\maketitle

%%%%%%%%%%%%%%%%%%%%%%%%%%%%%%%%%%%%%%%%%%%%%%%%%%%%%%%%%
%%                                                     %%
%%          Introduction                               %%
%%                                                     %%
%%%%%%%%%%%%%%%%%%%%%%%%%%%%%%%%%%%%%%%%%%%%%%%%%%%%%%%%%

\section{Introduction}

Important examples of eigenfunctions of 
the Macdonald \cite{M2} $q$-difference operator are
Macdonald \cite{M2} polynomials and Chalykh's \cite{C} 
Baker-Akhiezer functions.
In this paper we commence with a 
detailed study of the general spectral analysis of the Macdonald difference
operator. We construct eigenfunctions that are essentially characterized by 
the requirement that
they behave as plane waves deep in a 
distinguished Weyl chamber. We call the eigenfunctions 
difference Harish-Chandra
series and relate them to the Macdonald polynomials 
and to the Baker-Akhiezer functions.
We work in the general set-up of Macdonald's 
recent book \cite{M}, which includes
Koornwinder's \cite{Ko} extension of 
the Macdonald theory.

The refined analysis of the spectral problem of the 
Macdonald difference operator involves
the so-called system of basic hypergeometric difference equations.
Its definition is based on Cherednik's key observation that the 
affine Hecke algebra $H$ admits a realization as difference-reflection
operators, in which the Macdonald difference operator arises as the 
difference reduction of
a particular central element of $H$ (see, e.g., \cite{Ch}).
As such, the Macdonald difference operator is part of a 
commutative algebra $\mathbb{D}$ of difference operators,
obtained as the difference reduction of the full center of 
the affine Hecke algebra. It follows from this that $\mathbb{D}$ is isomorphic
to an algebra $A_0^\prime$ of Weyl group invariant regular functions on a 
complex torus $T^\prime=\hbox{Hom}_{\mathbb{Z}}(L^\prime,\mathbb{C}^\times)$,
where $L^\prime$ is either the co-weight lattice or the co-root lattice of 
the underlying irreducible
root system $R$. A concrete realization of this
isomorphism is given by a difference analog 
$\gamma: \mathbb{D}\rightarrow A_0^\prime$ of the Harish-Chandra isomorphism
(cf., e.g., \cite[\S 3]{Ch1}),
which assigns to the difference operator $D\in\mathbb{D}$ its 
asymptotic leading term deep in a distinguished Weyl chamber, 
up to a suitable twist.

The system of basic hypergeometric difference equations is then given by
\begin{equation}\label{system}
Df=\chi(D)f,\qquad \forall D\in\mathbb{D},
\end{equation}
where $\chi:\mathbb{D}\rightarrow \mathbb{C}$ is an algebra homomorphism.
The algebra homomorphisms $\chi: \mathbb{D}\rightarrow \mathbb{C}$
are naturally parameterized by the orbit space $T^\prime/W_0$, by 
composing the Harish-Chandra isomorphism $\gamma$
with the evaluation map $A_0^\prime\ni p\mapsto p(t)$ 
for $W_0t\in T^\prime/W_0$.

The system \eqref{system} of difference equations 
can be considered over various classes of (formal)
trigonometric functions on the complexification $V_{\mathbb{C}}$
of the ambient Euclidean space $V$ of the root system $R$ (the corresponding
period lattice of the trigonometric functions is
of the form $\frac{2\pi\sqrt{-1}}{\log(q)}L$ where $L$ 
is dual to $L^\prime$ in a suitable
sense, as we will make precise in the main text).
The Macdonald polynomials are the solutions of \eqref{system} which are
Weyl group invariant trigonometric polynomials (the corresponding $\chi$ 
form the so-called polynomial spectrum of the system \eqref{system}).
The difference Harish-Chandra series we construct are 
formal power series solutions of \eqref{system}. 

Before giving a precise outline of the contents of the paper, we first mention
three alternative contexts of the system \eqref{system} of basic 
hypergeometric difference equations
which are important as guiding principle 
for the present work. Firstly, Cherednik \cite{Ch0}, \cite{Ch1}
has related the system \eqref{system} of basic hypergeometric 
difference equations to quantum affine Khniznik-Zamolodchikov type
equations using the so-called Cherednik-Matsuo correspondence 
(see also \cite{Ma}, \cite{Kato} and \cite[Chpt. 2]{Chbook}).
This ties \eqref{system} to $q$-holonomic
systems, which also arise in integrable quantum 
field theories and in representation theory
of quantum affine algebras, see, e.g., \cite{EFK} for a survey.

Secondly, for special multiplicity labels the Macdonald 
difference operator arises as the radial
component of the quantum Casimir acting on a quantum compact 
symmetric space (cf., e.g., \cite{N}, \cite{L}).
In this situation the ring
$\mathbb{D}$ corresponds to the image of the center 
of the quantum universal enveloping algebra under the corresponding
radial component map for most symmetric spaces (see, e.g., \cite{L2}),
and the associated Macdonald polynomials arise as the corresponding 
quantum analogs of the elementary spherical functions.
We expect that the present study of the system \eqref{system} 
will provide a basic step towards the understanding of harmonic analysis on
quantum analogs of noncompact Riemannian symmetric spaces, cf. 
Harish-Chandra's well known classical approach
(see, e.g., \cite{GV} for an overview). Some initial steps towards
the harmonic analysis on quantum noncompact Riemannian symmetric
spaces can be found in e.g. \cite{KS1}, \cite{KS2}, \cite{St} and
\cite{L2}.

Thirdly, the Macdonald difference operator is essentially the Hamiltonian
of the quantum relativistic integrable system of Calogero-Moser type 
(cf., e.g.,  \cite{RS}).
For generic multiplicity labels $\mathbb{D}$ then 
corresponds to the associated algebra of quantum conserved integrals.
The algebra of quantum conserved integrals
is known to be strictly larger than $\mathbb{D}$ for 
special values of the multiplicity labels, in which 
case one speaks of algebraic integrability.
In this situation, Chalykh's \cite{C} Baker-Akhiezer 
functions arise as particular solutions of \eqref{system}.

In the last context the classical analog of the 
system \eqref{system} has been studied in full extent by
Heckman and Opdam \cite{HO}, \cite{O}, \cite{H1}, see also 
the overview in the first part of the book \cite{HS}.
In this paper we develop the theory in close parallel 
to the theory of Heckman and Opdam.

We now proceed to give a detailed description of the content of the paper.
In Section \ref{section1} we recall Cherednik's basic representation of the
affine Hecke algebra, we give the corresponding construction of the 
Macdonald difference operator
and we define the associated commutative 
ring $\mathbb{D}$ of difference operators, following closely Macdonald's
book \cite{M}.

In Section \ref{section2} we analyze the structure 
of the commutative algebra $\mathbb{D}$ of difference
operators in detail.
We consider an algebra $\mathbb{D}_{\mathcal{R}}(L^\prime)^{W_0}$ 
of difference operators containing $\mathbb{D}$.
It consists of $W_0$-invariant difference operators 
with step-sizes from the lattice $L^\prime$
and with coefficients in a suitable $W_0$-invariant 
algebra $\mathcal{R}$ of rational
trigonometric functions, where $W_0$ is the 
Weyl group of the underlying root system $R$.
The functions from $\mathcal{R}$ satisfy the 
essential additional property that they converge deep
in a distinguished Weyl chamber of $V$.
Using a rank reduction argument in an 
analogous manner to the differential theory (see, e.g., \cite{HS})
we prove that the Harish-Chandra homomorphism $\gamma$ also defines 
an algebra isomorphism from
the centralizer of the Macdonald 
difference operator in $\mathbb{D}_{\mathcal{R}}(L^\prime)^{W_0}$
onto $A_0^\prime$. As a consequence, we obtain the 
result that $\mathbb{D}$ equals the centralizer
of the Macdonald difference operator in 
$\mathbb{D}_{\mathcal{R}}(L^\prime)^{W_0}$.
We furthermore obtain a simple criterion 
(Corollary \ref{ok2}) when the centralizer of
the Macdonald difference operator in $\mathbb{D}_{\mathcal{R}}(L^\prime)$
is strictly larger than $\mathbb{D}$.

In Section \ref{HCsection}
we construct the difference analogs of the Harish-Chandra series.
For root systems $R$ of type A, the Harish-Chandra series solutions have been
considered before by Etingof and Kirillov Jr. \cite{EK} and 
Kazarnovski-Krol \cite{K}.
For generic spectral parameters we give a basis of 
the corresponding solution space
consisting of Harish-Chandra series. We relate the 
Harish-Chandra series to the Macdonald polynomial
when the algebra homomorphism $\chi$ is in the polynomial spectrum.
We furthermore show that the difference Harish-Chandra 
series reduce to Chalykh's \cite{C} Baker-Akhiezer function
when the system \eqref{system} is algebraically integrable 
in the sense of \cite{C} and \cite{ES}.
Finally, in Section \ref{sectionCN}
we provide a list of notations used throughout the paper.

In this paper we entirely focus on the algebraic theory of 
the difference Harish-Chandra series. The analytic theory, 
in particular the convergence of the difference 
Harish-Chandra series deep in a distinguished Weyl chamber, 
is part of research in progress of the second author with Michel van Meer. 
Various other natural topics, such as the basic 
hypergeometric function (the analog of the spherical function),
connection matrices, duality, bispectrality and 
quantum group interpretations, are also subject to future research.

{\it Acknowledgments:} The first author was partially supported by
NSA grant no. H98230-05-1-0077. The second author was supported by
the Netherlands Organization for Scientific Research (NWO) in the
VIDI-project ``Symmetry and modularity in exactly solvable models''.

%%%%%%%%%%%%%%%%%%%%%%%%%%%%%%%%%%%%%%%%%%%%%%%%%%%%%%%%%%
%%                                                      %%
%%  Macdonald difference operators                      %%
%%                                                      %%
%%%%%%%%%%%%%%%%%%%%%%%%%%%%%%%%%%%%%%%%%%%%%%%%%%%%%%%%%%

\section{Macdonald difference operators}\label{section1}

%%%%%%%%%%%%%%%%%%%%%%%%%%%%%%%%%%%%%%%%%%%%%%%%%%%%%
%%                                                 %%
%%      Basic root system data                     %%
%%                                                 %%
%%%%%%%%%%%%%%%%%%%%%%%%%%%%%%%%%%%%%%%%%%%%%%%%%%%%%

\subsection{Root system data}\label{Rsd}

Macdonald polynomials, which are naturally attached to affine root systems,
have been successfully analyzed
through the study of Cherednik's \cite{Ch} double affine
Hecke algebra. Their nonreduced extensions, the so-called 
Macdonald-Koornwinder
polynomials (see \cite{Ko}), have been incorporated 
in the theory using a suitable extension
of the double affine Hecke algebra, see Noumi \cite{N} and Sahi \cite{S}.
In order to capture the closely related theories all at once, we
follow Macdonald's recent book \cite{M} and adapt its conventions 
and notations as much as possible throughout the paper.

The basic structure underlying the Cherednik-Macdonald 
theory then consists of a pair $(R,R^\prime)$ of finite,
reduced irreducible crystallographic root systems in an Euclidean 
space $(V,\langle\cdot,\cdot\rangle)$, a
pair $(L,L^\prime)$ of lattices in $V$, 
and a pair $(S,S^\prime)$ of irreducible affine root systems.
The definition of these pairs depends on three different cases, 
to which we will refer to as case {\bf a}, {\bf b}
and {\bf c} throughout the paper. Before listing 
the pairs for each of the three cases,
we first introduce some general notation.

Let $(V,\langle\cdot,\cdot\rangle)$ be a finite dimensional Euclidean space.
Let $\widehat{V}$ be the space of affine linear, real functions on
$V$. Denote $c\in\widehat{V}$ for the constant function one.  We identify
$\widehat{V}\simeq V\oplus\mathbb{R}c$ as real vector space via the scalar 
product $\langle\cdot,\cdot\rangle$ by associating 
to $v+sc\in V\oplus \mathbb{R}c$
($v\in V$, $s\in\mathbb{R}$) the affine linear 
functional $v^\prime\mapsto \langle
v,v^\prime\rangle+s$ on $V$. Let $D: \widehat{V}\rightarrow V$ be
the gradient map, defined by $D(v+sc)=v$ for $v\in V$ and $s\in \mathbb{R}$.
We extend $\langle \cdot,\cdot\rangle$ to a positive semi-definite
bilinear form on $\widehat{V}$ by
\[\langle f,g\rangle:=\langle Df,Dg\rangle,\qquad f,g\in
\widehat{V}.
\]
We define for $f\in \widehat{V}$ with $Df\not=0$ the associated
co-vector by $f^\vee=2f/\|f\|^2\in \widehat{V}$.

Let $R\subset V$ be a finite, reduced irreducible 
crystallographic root system,
with associated Weyl group $W_0\subset \textup{O}(V)$ generated by
the orthogonal reflections $s_\alpha$ in the hyperplanes 
$\alpha^\perp\subset V$ ($\alpha\in R$).
Let $R^\vee=\{\alpha^\vee\}_{\alpha\in R}$ be the associated co-root
system. Denote $Q=Q(R)$ and $P=P(R)$ for the root
lattice and the weight lattice of $R$ respectively,
which are $W_0$-invariant lattices in $V$ satisfying $Q\subset P$.
We write $Q^\vee=Q(R^\vee)$ and $P^\vee=P(R^\vee)$ for the
co-root lattice and co-weight lattice of $R$ in $V$.
We define the affine Weyl group $W_{Q^\vee}$ and the extended
affine Weyl group $W_{P^\vee}$ of $R$ as the
corresponding semi-direct product groups
\[W_{Q^\vee}=W_0\ltimes Q^\vee,\qquad W_{P^\vee}=W_0\ltimes
P^\vee.
\]
The canonical action of $W_0$ on $V$ extends to a faithful action of the
extended affine Weyl group on $V$ with the lattice $P^\vee$ acting by
translations,
\[t(\lambda)(v)=v+\lambda,\qquad v\in V
\]
for $\lambda\in P^\vee$.
The space $\widehat{V}$ inherits a left 
$W_{P^\vee}$-module structure by transposition of the
$W_{P^\vee}$-action on $V$.

The set
\[S(R)=\{\alpha+rc\,\, | \,\, \alpha\in R,\,\, r\in\mathbb{Z} \}
\]
defines an irreducible, reduced affine root system in $\widehat{V}$ 
with underlying finite, gradient
root system $D(S(R))=R$. The associated affine Weyl group,
generated by the orthogonal reflections in the affine
hyperplanes $f^{-1}(0)$ ($f\in S(R)$), is isomorphic to
$W_{Q^\vee}$. The affine root system $S(R)$ is invariant 
under the action of the
extended affine Weyl group $W_{P^\vee}$. Observe furthermore that
$S(R)^\vee=\{f^\vee \,\, | \,\, f\in S(R)\}$ is an irreducible
reduced affine root system in $\widehat{V}$ with underlying gradient root
system $R^\vee$.

The affine root
systems $S(R)$ and $S(R)^\vee$ described above are all reduced. The
nonreduced irreducible affine root systems are root subsystems of the
affine root system of type $C^\vee C$, which we now proceed to
describe. Denote $\{\epsilon_i\}_{i=1}^n$ for the
standard orthonormal basis of $\mathbb{R}^n$. Write 
$R_C\subset \mathbb{R}^n$ for the root
system of type $C_n$ given by the set of roots $\pm\epsilon_i\pm\epsilon_j$
($1\leq i<j\leq n$) and $\pm 2\epsilon_m$ ($m=1,\ldots,n$), where
all sign combinations are allowed (we suppress the dependence on 
$n\in\mathbb{Z}_{>0}$
in the notations). The nonreduced affine root
system of type $C^\vee C_n$ can now be realized as the set
\[S_{nr}:=\{\pm\epsilon_i\pm\epsilon_j+rc,\, \pm \epsilon_m +\frac{rc}{2},\,
\pm 2\epsilon_m+rc \,\, | \,\, 1\leq i<j\leq n,\,\,
1\leq m\leq n,\,\, r\in\mathbb{Z} \}
\]
in $\widehat{\mathbb{R}^n}$ (all sign combinations are allowed).
Note that the associated gradient root system $D(S_{nr})$ is the
nonreduced root system of type $\textup{BC}_n$, containing $R_C$ as the
root subsystem of roots of squared length greater than or equal to 
$2$ in $D(S_{nr})$.

We refer to the three cases \cite[(1.4.1)]{M}, \cite[(1.4.2)]{M} and
\cite[(1.4.3)]{M} of the Cherednik-Macdonald theory as
cases {\bf a}, {\bf b} and {\bf c} respectively. We list below for each case
the pairs ($R,R^\prime$), $(L,L^\prime)$ and $(S,S^\prime)$.

\begin{enumerate}
\item[Case {\bf a.}] $(R,R^\prime)=(R,R^\vee)$, $(L,L^\prime)=(P(R),P(R^\vee))$
and $(S,S^\prime)=(S(R), S(R^\vee))$ with $R\subset V$ a finite, reduced, 
irreducible
crystallographic root system, normalized such 
that long roots have squared length two.
\item[Case {\bf b.}] $(R,R^\prime)=(R,R)$, $(L,L^\prime)=(P(R^\vee),P(R^\vee))$
and $(S,S^\prime)=(S(R)^\vee,S(R)^\vee)$ with $R\subset V$ a finite,
reduced, irreducible crystallographic root system,
normalized such that long roots have squared
length two.
\item[Case {\bf c.}] $V=\mathbb{R}^n$,
$(R,R^\prime)=(R_C,R_C)$,
$(L,L^\prime)=(Q(R_C^\vee), Q(R_C^\vee))$ and
$(S,S^\prime)=(S_{nr}, S_{nr})$.
\end{enumerate}

Note that for case {\bf c}, long roots in $R$ have squared norm four.
Note furthermore that for case {\bf a} and {\bf b} we have
$L=P(R^{\prime\vee})$ and $L^\prime=P(R^\vee)$, while for case
{\bf c} $L=Q(R^{\prime\vee})$ and $L^\prime=Q(R^{\vee})$.
We can thus
associate to each case a pair of (extended) affine Weyl groups
\[(W,W^\prime)=(W_{L^\prime}, W_L),
\]
which are the extended affine Weyl groups of $R$ and $R^\prime$
for case {\bf a} and {\bf b}, and the affine Weyl group of $R$ and
$R^\prime$ for case {\bf c}.

For each of the three cases, the (extended) affine Weyl group $W$
(respectively $W^\prime$) preserves $S$ (respectively $S^\prime$).
In case {\bf a} and {\bf b} the decomposition of $S$ into $W$-orbits
coincides with the decomposition $S=S_s\cup S_l$ of $S$ into the set
$S_s$ of short roots and the set $S_l$ of long roots. By the normalization
of the underlying finite root system, $S_l$ are the roots in $S$ whose squared
length equals two with respect to the semi-positive definite form 
$\langle\cdot,\cdot\rangle$
of $\widehat{V}$. In case {\bf c} the affine root system $S$ 
has five $W$-orbits, namely
\[\mathcal{O}_1=\{\pm \epsilon_m+rc\, | \, 1\leq m\leq n,\, r\in \mathbb{Z}\},
\quad \mathcal{O}_2=2\mathcal{O}_1,\quad
\mathcal{O}_3=\mathcal{O}_1+\frac{c}{2},\quad \mathcal{O}_4=2\mathcal{O}_3
\]
and $\mathcal{O}_5=\{\pm\epsilon_i\pm\epsilon_j+rc\, | \, 1\leq i<j\leq n,
r\in\mathbb{Z} \}$.

Let $S_1\subset S$ be the reduced affine root subsystem of
indivisible affine roots in $S$. Then $S_1$ is
$W$-invariant, $S_1=S$ for case {\bf a} and {\bf b} while
$S_1=S(R_C)^\vee$ for case {\bf c}, which has three $W$-orbits
$\mathcal{O}_1$, $\mathcal{O}_3$ and $\mathcal{O}_5$.

We fix a basis $\{\alpha_1,\ldots,\alpha_n\}$ of the root system
$R$. We associate to it a basis $\{a_0,a_1,\ldots,a_n\}$ of
the corresponding affine root system $S$ by
\begin{equation*}
\{a_0,a_1,\ldots,a_n\}=
\begin{cases}
\{-\varphi+c,\alpha_1,\ldots,\alpha_n\},\qquad &\hbox{case {\bf a}},\\
\{-\varphi^\vee+c,\alpha_1^\vee,\ldots,\alpha_n^\vee\},\qquad
&\hbox{case {\bf b}},\\
\{-\varphi^\vee+\frac{c}{2},\alpha_1^\vee,\ldots,\alpha_n^\vee\},\qquad
&\hbox{case {\bf c}},
\end{cases}
\end{equation*}
where $\varphi\in R$ is the highest root with respect to the given basis
$\{\alpha_1,\ldots,\alpha_n\}$ of $R$ (which is always a long root).
Note that $\{a_0,\ldots,a_n\}$ is also a basis of $S_1$,
and that each affine root $a\in S_1$
is $W$-conjugate to a simple root $a_j$ ($j\in\{0,\ldots,n\}$).
We set $\Delta=\{a_1,\ldots,a_n\}$, which is a basis of the
root subsystem $R^{\prime\vee}=D(S_1)$ of indivisible roots in
the gradient root system $D(S)\subset V$ of $S$.
In particular, $\Delta$ is a $\mathbb{Z}$-basis of $Q(R^{\prime\vee})$.

We denote $s_i=s_{a_i}\in W_{Q^\vee}\subset W$ ($i=0,\ldots,n$) for the
reflection associated to the simple root $a_i\in S_1$ ($i=0,\ldots,n$).
The affine Weyl group $W_{Q^\vee}$ is a Coxeter group with respect to 
the simple reflections
$s_i$ ($i=0,\ldots,n$).
Let $S_1^+$ (respectively $S_1^-$) be the positive (respectively
negative) affine roots in $S_1$ with respect to the basis $\{a_i\}_{i=0}^n$
of $S_1$. The length $l(w)$ of $w\in W$ is defined by
\[l(w)=\#\bigl(S_1^+\cap w^{-1}S_1^-\bigr).
\]
The finite abelian subgroup $\Omega=\{w\in W \, | \, l(w)=0\}$ of $W$ is
isomorphic to $L^\prime/Q(R^\vee)$, and we have the semi-direct product
decomposition $W\simeq W_{Q^\vee}\rtimes\Omega$.
The action of $\Omega$ on $\widehat{V}$
restricts to a faithful action on the finite set $\{a_i\}_{i=0}^n$
of simple roots, hence we may and will view $\Omega$ as a permutation
subgroup of the index set $\{0,\ldots,n\}$.
We denote $\mathbb{C}[W]$ and $\mathbb{C}[\Omega]$ for the complex
group algebra of $W$ and $\Omega$, respectively.

%%%%%%%%%%%%%%%%%%%%%%%%%%%%%%%%%%%%%%%%%%%%%%%%%%%%%%%%%%%%%%%%
%%                                                            %%
%%    The algebra of difference-reflection operators          %%
%%                                                            %%
%%%%%%%%%%%%%%%%%%%%%%%%%%%%%%%%%%%%%%%%%%%%%%%%%%%%%%%%%%%%%%%%

\subsection{The algebra of difference-reflection
operators}\label{DRO}

Denote $A=\mathbb{C}[L]$ and $A^\prime=\mathbb{C}[L^\prime]$ for
the group algebras of the lattices $L$ and $L^\prime$ over $\mathbb{C}$,
respectively. We denote the canonical basis of $A$
(respectively $A^\prime$) by $\{z^\lambda \, | \, \lambda\in L\}$
(respectively $\{\xi^{\lambda^\prime} \, | \, \lambda^\prime\in
L^\prime\}$), so that
\[z^\lambda z^\mu=z^{\lambda+\mu},\qquad z^0=1
\]
and similarly for the $\xi^{\lambda^\prime}$. We fix throughout
the paper $0<q<1$. We extend the definition of the monomials $z^\lambda$
in $A$ (respectively $\xi^{\lambda^\prime}$ in $A^\prime$) by defining
\[z^{\lambda+rc}:=q^rz^\lambda,\qquad
\xi^{\lambda^\prime+rc}:=q^r\xi^{\lambda^\prime}
\]
for $\lambda\in L$, $\lambda^\prime\in L^\prime$ and
$r\in\mathbb{R}$. In particular, for $a\in S$ and $b\in S^\prime$ we have
$Da\in Q(R^{\prime\vee})\subseteq L$ and $Db\in Q(R^\vee)\subseteq
L^\prime$, hence $z^a\in A$ and $\xi^b\in A^\prime$ are well defined.

Let $V_{\mathbb{C}}=V\oplus \sqrt{-1}\,V$ be the complexification of
$V$, and extend $\langle\cdot,\cdot\rangle$ to a complex bilinear
form on $V_{\mathbb{C}}$. We view $A$ and $A^\prime$ as
subalgebras of the algebra of complex analytic functions on $V_{\mathbb{C}}$
by interpreting the canonical basis elements
$z^\lambda$ and $\xi^{\lambda^\prime}$ as complex plane
waves on $V_{\mathbb{C}}$,
\begin{equation}\label{planewave}
z^\lambda(v)=q^{\langle \lambda,v\rangle},\qquad
\xi^{\lambda^\prime}(v)=q^{\langle\lambda^\prime,v\rangle}
\end{equation}
for $\lambda\in L$, $\lambda^\prime\in L^\prime$ and $v\in
V_{\mathbb{C}}$. The (extended) affine Weyl groups $W$ and $W^\prime$
act on $A$ and $A^\prime$ as algebra automorphisms by
transposing their respective actions on $V_{\mathbb{C}}$. Concretely,
for $w=vt(\lambda^\prime)\in W$ and $w^\prime=vt(\lambda)\in W^\prime$ 
with $v\in W_0$,
$\lambda^\prime\in L^\prime$ and $\lambda\in L$
we have
\[w(z^\mu)=q^{-\langle\lambda^\prime,\mu\rangle}z^{v\mu},\qquad
w^\prime(\xi^{\mu^\prime})=q^{-\langle\lambda,\mu^\prime\rangle}\xi^{v\mu^\prime}
\]
for $\mu\in L$ and $\mu^\prime\in L^\prime$. Furthermore,
$w(z^a)=z^{wa}$ and $w^\prime(\xi^b)=\xi^{w^\prime b}$ for $a\in S$
and $b\in S^\prime$.
We write $\mathcal{Q}$ and $\mathcal{Q}^\prime$
for the quotient fields of $A$ and $A^\prime$ respectively. The
$W$-action on $A$ (respectively $W^\prime$-action on
$A^\prime$) extends uniquely to a $W$-action on $\mathcal{Q}$
(respectively $W^\prime$-action on $\mathcal{Q}^\prime$) by field
automorphisms. We write $(wf)(z)\in\mathcal{Q}$ for the 
rational function obtained by acting by $w\in W$
on $f(z)\in\mathcal{Q}$.

%%%%%%%%%%%%%%%%%%%%%%%%%%%%%%%%%%%%%%%%%%%%%%%%%%%%%%%%%%%%%%%%%
\begin{defi}\label{Rdef}
Let $\mathcal{R}$ be the
complex subalgebra of $\mathcal{Q}$
generated \textup{(}as an algebra\textup{)} by the elements
\[\frac{1}{1-rz^\alpha},\qquad \alpha\in R^{\prime\vee},\,\, r\in\mathbb{C}.
\]
\end{defi}
%%%%%%%%%%%%%%%%%%%%%%%%%%%%%%%%%%%%%%%%%%%%%%%%%%%%%%%%%%%%%%%%%
Note that $(1-rz^a)^{-1}\in\mathcal{R}$ for all $r\in\mathbb{C}$ and $a\in S$.
In particular, $\mathcal{R}$ is a $W$-module subalgebra of $\mathcal{Q}$.
Furthermore,
\begin{equation}\label{half}
\frac{1}{1-rz^\alpha}=1-\frac{1}{1-r^{-1}z^{-\alpha}},
\qquad \alpha\in R^{\prime\vee},\,\,
r\in\mathbb{C}^\times,
\end{equation}
hence $\mathcal{R}$ is already generated by
the elements $\frac{1}{1-rz^{-\alpha}}$ with $r\in\mathbb{C}$ and
with roots $\alpha\in R^{\prime\vee,+}:=
R^{\prime\vee}\cap\mathbb{Z}_{\geq 0}\Delta$.

We associate to the $W$-module algebra 
$\mathcal{R}$ the smash-product algebras
\[\mathbb{D}_{\mathcal{R}}(L^\prime):=\mathcal{R}\# t(L^\prime)\subset
\mathcal{R}\# W=:\mathbb{D}_{\mathcal{R}}(W).
\]
In other words, $\mathbb{D}_{\mathcal{R}}(W)$ is the complex, unital,
associative algebra such that
\[\mathbb{D}_{\mathcal{R}}(W)\simeq
\mathcal{R}\otimes_{\mathbb{C}}
\mathbb{C}[W]
\]
as complex vector spaces, such that the canonical linear
embeddings $\mathcal{R}, \mathbb{C}[W]\hookrightarrow 
\mathbb{D}_{\mathcal{R}}(W)$
are unital algebra embeddings and such that the cross relations
\[(f_1\otimes w_1)(f_2\otimes w_2)=f_1w_1(f_2)\otimes w_1w_2
\]
holds for $f_i\in\mathcal{R}$ and $w_i\in W$ (and similarly for
$\mathbb{D}_{\mathcal{R}}(L^\prime)$).

%%%%%%%%%%%%%%%%%%%%%%%%%%%%%%%%%%%%%%%%%%%%%%%%%%%%%%%%%%%%%%
\begin{defi}
We call $\mathbb{D}_{\mathcal{R}}(W)$ \textup{(}respectively
$\mathbb{D}_{\mathcal{R}}(L^\prime)$\textup{)} the algebra of
difference-reflection \textup{(}respectively difference\textup{)}
operators with coefficients in $\mathcal{R}$.
\end{defi}
%%%%%%%%%%%%%%%%%%%%%%%%%%%%%%%%%%%%%%%%%%%%%%%%%%%%%%%%%%%%%%%%%%%%%%%%
The terminology is justified by the canonical, faithful action of
$\mathbb{D}_{\mathcal{R}}(W)$
(respectively $\mathbb{D}_{\mathcal{R}}(L^\prime)$) on $\mathcal{Q}$
as difference-reflection (respectively difference) operators by
\[\bigl(f\otimes w\bigr)(g):=fw(g),\qquad f\in \mathcal{R},\,\, w\in
W,\,\, g\in \mathcal{Q}.
\]

Since $\mathcal{R}$ is a $\mathbb{C}[t(X)]$-submodule algebra of $\mathcal{Q}$
for all lattices $X\subset V$, we can similarly define the algebra
$\mathbb{D}_{\mathcal{R}}(X)=\mathcal{R}\# t(X)$ of difference
operators with coefficients in $\mathcal{R}$ and step-sizes from
$X$. For lattices $L^\prime\subseteq X\subset V$, 
$\mathbb{D}_{\mathcal{R}}(L^\prime)$
canonically embeds as a subalgebra into $\mathbb{D}_{\mathcal{R}}(X)$.

A second extension of $\mathbb{D}_{\mathcal{R}}(L^\prime)$
involves the coefficients of the difference
operators. Recall that
$Q(R^{\prime\vee})=\mathbb{Z}\Delta\subseteq L$. We define
$\mathbb{C}[[z^{-\Delta}]]$ to be the algebra of formal power series
\begin{equation}\label{f}
f(z)=\sum_{x\in\mathbb{Z}_{\geq 0}\Delta}C_xz^{-x},\qquad
C_x\in\mathbb{C}.
\end{equation}
A lattice $X\subset V$ acts
by algebra automorphisms on $\mathbb{C}[[z^{-\Delta}]]$ as
\[t(\nu)\bigl(f(z)\bigr)=\sum_{x\in\mathbb{Z}_{\geq
0}\Delta}C_xq^{\langle\nu,x\rangle}z^{-x},\qquad \nu\in X,
\]
with $f(z)$ given by \eqref{f}. The algebra of difference
operators with step-sizes from $X$ and coefficients from
$\mathbb{C}[[z^{-\Delta}]]$ is the associated smash-product
algebra $\mathbb{D}_{\mathbb{C}[[z^{-\Delta}]]}(X):=
\mathbb{C}[[z^{-\Delta}]]\# t(X)$.
Again, for lattices $L^\prime\subseteq X\subset
V$, $\mathbb{D}_{\mathbb{C}[[z^{-\Delta}]]}(L^\prime)$ canonically
embeds into $\mathbb{D}_{\mathbb{C}[[z^{-\Delta}]]}(X)$.

The $W$-module algebra $\mathcal{R}$ can be canonically
embedded as a $t(X)$-module subalgebra in $\mathbb{C}[[z^{-\Delta}]]$ 
for any lattice $X\subset V$
using, for $r\in\mathbb{C}$ and for $\alpha\in R^{\prime\vee,+}$, 
the formal series expansion
\[\frac{1}{1-rz^{-\alpha}}=\sum_{m=0}^{\infty}r^mz^{-m\alpha}.
\]
Correspondingly, we have a canonical embedding
\[\mathbb{D}_{\mathcal{R}}(X)\hookrightarrow
\mathbb{D}_{\mathbb{C}[[z^{-\Delta}]]}(X)
\]
of algebras. We identify $\mathbb{D}_{\mathcal{R}}(X)$ with its image
under this embedding in the remainder of the paper.

We end this subsection by constructing an action 
of $\mathbb{D}_{\mathbb{C}[[z^{-\Delta}]]}(L^\prime)$
as difference operators on a space of formal linear 
combinations of complex plane waves
on $V_{\mathbb{C}}$, which will be of importance in the
construction of the Harish-Chandra series in Section
\ref{HCsection}.

Let $M$ be the complex commutative subalgebra of analytic functions on
$V_{\mathbb{C}}$ spanned by the complex plane waves $z^u$ ($u\in
V_{\mathbb{C}}$),
\[z^u(v):=q^{\langle u,v\rangle},\qquad v\in V_{\mathbb{C}}.
\]
The algebra $A$ naturally identifies with the subalgebra of $M$
spanned by the plane waves $z^\lambda$ ($\lambda\in L$), cf. 
\eqref{planewave}.

Denote $\overline{M}$ for the complex vector space of
formal series
\[\sum_{u\in\mathcal{C}}K_uz^u \qquad (K_u\in\mathbb{C})
\]
for subsets $\mathcal{C}\subset V_{\mathbb{C}}$ which are contained in
some finite union of sets of the form $\lambda-\mathbb{Z}_{\geq 0}\Delta$ 
($\lambda\in V_{\mathbb{C}}$).
In other words, an element $F(z)\in \overline{M}$ is a 
finite linear combination of formal
series of the form
\begin{equation}\label{leadingexponent}
F_\lambda(z)=\sum_{x\in\mathbb{Z}_{\geq
0}\Delta}K_\lambda(x)z^{\lambda-x}\qquad (K_\lambda(x)\in\mathbb{C})
\end{equation}
where $\lambda\in V_{\mathbb{C}}$.
Note that $\overline{M}$ is canonically a $\mathbb{C}[[z^{-\Delta}]]$-module,
and $M$ embeds in $\overline{M}$ as a $\mathbb{C}[z^{-\Delta}]$-module.
%%%%%%%%%%%%%%%%%%%%%%%%%%%%%%%%%%%%%%%%%%%%%%%%%%%%%%%%%%%%%%%%%%%%%%%%%%
\begin{lem}\label{M}
The action of $\mathbb{C}[[z^{-\Delta}]]$ on $\overline{M}$
extends to an action of the algebra 
$\mathbb{D}_{\mathbb{C}[[z^{-\Delta}]]}(L^\prime)$ of
difference operators with coefficients in $\mathbb{C}[[z^{-\Delta}]]$ on
$\overline{M}$, with $\lambda^\prime\in L^\prime$ acting as
\[t(\lambda^\prime)\bigl(\sum_{u\in\mathcal{C}}K_uz^u\bigr)=\sum_{u\in\mathcal{C}}
K_uq^{-\langle\lambda^\prime,u\rangle}z^u.
\]
\end{lem}
%%%%%%%%%%%%%%%%%%%%%%%%%%%%%%%%%%%%%%%%%%%%%%%%%%%%%%%%%%%%%%%%%%%%%%%%%%%
\begin{proof}
This follows from a direct computation.
\end{proof}
%%%%%%%%%%%%%%%%%%%%%%%%%%%%%%%%%%%%%%%%%%%%%%%%%%%%%%%%%%%%%%%%%%%%%%%%%%%%

%%%%%%%%%%%%%%%%%%%%%%%%%%%%%%%%%%%%%%%%%%%%%%%%%%%%%%%%%%%%%%%%%%%%%%%
%%                                                                   %%
%%      Cherednik's difference-reflection operators                  %%
%%                                                                   %%
%%%%%%%%%%%%%%%%%%%%%%%%%%%%%%%%%%%%%%%%%%%%%%%%%%%%%%%%%%%%%%%%%%%%%%%

\subsection{Cherednik's difference-reflection operators}\label{Chers}

Cherednik's dif\-fe\-ren\-ce-re\-flec\-tion operators realize the
(extended) affine Hecke algebra inside $\mathbb{D}_{\mathcal{R}}(W)$. We
recall the construction in this subsection.
The affine Hecke algebras and its realization in $\mathbb{D}_{\mathcal{R}}(W)$
depend on multiplicity labels, which we now first introduce.

The complex vector space $\mathbb{C}(S)^W$ consisting of $W$-invariant
functions $\underline{k}: S\rightarrow \mathbb{C}$ is
called the space of multiplicity labels associated to $S$.
We denote $k_a=\underline{k}(a)$
for $a\in S$ and $k_i=\underline{k}(a_i)$ for $i=0,\ldots,n$.
For case {\bf a} and {\bf b} the vector space $\mathbb{C}(S)^W$
is at most two dimensional, while it is five dimensional for case {\bf c}.

Depending on the case under consideration, we associate to
a multiplicity label $\underline{k}\in\mathbb{C}(S)^W$ a
dual multiplicity label $\underline{k}^\prime\in\mathbb{C}(S^\prime)^{W^\prime}$
as follows.
\begin{enumerate}
\item[Case {\bf a.}] $\underline{k}^\prime(b)=\underline{k}(Db^\vee)$ 
for all $b\in S^\prime$.
\item[Case {\bf b.}] $\underline{k}^\prime(b)=\underline{k}(b)$ 
for all $b\in S^\prime$.
\item[Case {\bf c.}] We write $\kappa_j$ for the value of $\underline{k}$
at the $W$-orbit $\mathcal{O}_j$ of $S$ ($j=1,\ldots,5$).
We define the dual multiplicity label $\underline{k}^\prime$ by
assigning the following value $\kappa_j^\prime$ to the $W^\prime$-orbit
$\mathcal{O}_j$ of $S^\prime$ ($j=1,\ldots,5$),
\begin{equation*}
\begin{split}
\kappa_1^\prime&=\frac{1}{2}(\kappa_1+\kappa_2+\kappa_3+\kappa_4),\\
\kappa_2^\prime&=\frac{1}{2}(\kappa_1+\kappa_2-\kappa_3-\kappa_4),\\
\kappa_3^\prime&=\frac{1}{2}(\kappa_1-\kappa_2+\kappa_3-\kappa_4),\\
\kappa_4^\prime&=\frac{1}{2}(\kappa_1-\kappa_2-\kappa_3+\kappa_4),\\
\kappa_5^\prime&=\kappa_5.
\end{split}
\end{equation*}
\end{enumerate}
Denote $\mathbb{C}(S_1)^W$ for the space of
multiplicity labels associated to $S_1$. Associated to a multiplicity label
$\underline{k}\in\mathbb{C}(S)^W$
we define an invertible multiplicity label 
$\underline{\tau}\in\mathbb{C}(S_1)^W$
and a dual invertible multiplicity label 
$\underline{\tau}^\prime\in\mathbb{C}(S_1)^W$
as follows. For cases {\bf a} and {\bf b}, we set
\[\tau_a=q^{\frac{1}{2}k_a}=\tau_a^\prime,\qquad a\in S_1=S.
\]
For case {\bf c}, we again write $\kappa_j$ for the
value of $\underline{k}$ on the $W$-orbit $\mathcal{O}_j$
($j=1,\ldots,5$). The invertible multiplicity labels $\underline{\tau}$
and $\underline{\tau}^\prime$ associated to $S_1$ are then defined by
\begin{equation*}
\begin{split}
\tau_a&=q^{\frac{1}{2}(\kappa_1+\kappa_2)},\quad
\tau_a^\prime=q^{\frac{1}{2}(\kappa_1-\kappa_2)},\qquad
a\in \mathcal{O}_1,\\
\tau_a&=q^{\frac{1}{2}(\kappa_3+\kappa_4)},\quad
\tau_a^\prime=q^{\frac{1}{2}(\kappa_3-\kappa_4)},\qquad a\in
\mathcal{O}_3,\\
\tau_a&=\tau_a^\prime=q^{\frac{1}{2}\kappa_5},
\qquad\qquad\qquad\qquad\quad\,\,\,\,
a\in\mathcal{O}_5.
\end{split}
\end{equation*}
We write $\tau_i=\tau_{a_i}$ and 
$\tau_i^\prime=\tau_{a_i}^\prime$ for $i=0,\ldots,n$.

The (extended) affine Hecke algebra depends on an invertible
multiplicity label $\underline{\tau}\in\mathbb{C}(S_1)^W$
(hence depends at most on three continuous, complex parameters).

%%%%%%%%%%%%%%%%%%%%%%%%%%%%%%%%%%%%%%%%%%%%%%%%%%%%%%%%%%%%%%%%%%%%%%%%%
\begin{defi}\label{HA}
Let $\underline{\tau}\in\mathbb{C}(S_1)^W$ be an 
invertible multiplicity label.
The \textup{(}extended\textup{)}
affine Hecke algebra $H(\underline{\tau})$
is the unital complex associative algebra generated by elements $T_i$
\textup{(}$i=0,\ldots,n$\textup{)} and $\omega$
\textup{(}$\omega\in\Omega$\textup{)} satisfying
\begin{enumerate}
\item[{\bf i.}] The linear map $\mathbb{C}[\Omega]\rightarrow
H$, defined by $\omega\mapsto\omega$,
is an algebra morphism.
\item[{\bf ii.}] The $T_i$ \textup{(}$i=0,\ldots,n$\textup{)}
satisfy the braid relations
\[T_iT_jT_i\cdots=T_jT_iT_j\cdots
\]
with $m_{ij}$ factors on both sides, for indices $i\not=j$
such that $s_is_j\in W_{Q^\vee}$ has finite order $m_{ij}$.
\item[{\bf iii.}] $\omega T_i\omega^{-1}=T_{\omega (i)}$ for $i=0,\ldots,n$
and $\omega\in\Omega$.
\item[{\bf iv.}] $(T_i-\tau_i)(T_i+\tau_i^{-1})=0$ for $i=0,\ldots,n$.
\end{enumerate}
\end{defi}
%%%%%%%%%%%%%%%%%%%%%%%%%%%%%%%%%%%%%%%%%%%%%%%%%%%%%%%%%%%%%%%%%%%%%%%%
For $w\in W$ and a reduced expression $w=\omega s_{i_1}s_{i_2}\cdots 
s_{i_{l(w)}}$
with $\omega\in\Omega$ and $i_j\in \{0,\ldots,n\}$ we write
\[T_w=\omega T_{i_1}T_{i_2}\cdots T_{i_{l(w)}}\in
H(\underline{\tau}).
\]
The $T_w$ ($w\in W$) are well defined
(independent of the reduced expression) and form a
linear basis of $H(\underline{\tau})$.

Fix a multiplicity label $\underline{k}\in\mathbb{C}(S)^W$
and denote $\underline{\tau},\underline{\tau}^\prime
\in\mathbb{C}(S_1)^W$ for the invertible multiplicity labels associated
to $\underline{k}$ as described above. The basic representation of 
$H(\underline{\tau})$
which we now proceed to define, 
depends also on $\underline{\tau}^\prime$ (in other words, it depends on
two additional complex parameters for case {\bf c}). To define the
basic representation,
we set for $a\in S_1$,
\begin{equation}\label{c}
c_a(z):=\frac{(1-\tau_a\tau_a^\prime
z^a)(1+\tau_a\tau_a^\prime{}^{-1}z^a)}{\tau_a(1-z^{2a})}\in\mathcal{R}.
\end{equation}
Note that $c_a(z)\in\mathcal{R}$ ($a\in S_1$) 
satisfies the elementary identity
\begin{equation}\label{ceq}
c_a(z)+c_{-a}(z)=\tau_a+\tau_a^{-1},
\end{equation}
and that
\[c_a(z)=\frac{(1-\tau_a^2z^a)}{\tau_a(1-z^a)}
\]
if $\tau_a=\tau_a^\prime$. We define the difference-reflection operators
\[T_i(\underline{k}):=\tau_i+c_{a_i}(z)(s_i-1)
\in \mathbb{D}_{\mathcal{R}}(W),\qquad i=0,\ldots,n.
\]
The following theorem is due to Cherednik \cite{Ch}, and due to Noumi
\cite{N} for case {\bf c}. The present uniform formulation
is from \cite[(4.3.10)]{M}.
%%%%%%%%%%%%%%%%%%%%%%%%%%%%%%%%%%%%%%%%%%%%%%%%%%%%%%%%%%
\begin{thm}\label{fundamental}
{\bf i)} For $i=0,\ldots,n$ the $W$-module subalgebra
$A\subset \mathcal{Q}$ is invariant
under the action of the difference-reflection operators
$T_i(\underline{k})\in\mathbb{D}_{\mathcal{R}}(W)$ on $\mathcal{Q}$.\\
{\bf ii)} The
assignment $T_i\mapsto T_i(\underline{k})|_{A}$
\textup{(}$i=0,\ldots,n$\textup{)}, together with the usual action
of $\omega\in\Omega$ on the $W$-module algebra $A$, uniquely
extends to a faithful representation of the
\textup{(}extended\textup{)} affine Hecke algebra
$H(\underline{\tau})$ on $A$.
\end{thm}
%%%%%%%%%%%%%%%%%%%%%%%%%%%%%%%%%%%%%%%%%%%%%%%%%%%%%%%%%%%%
Since $D\in\mathbb{D}_{\mathcal{R}}(W)$ is uniquely determined by
its action on $A\subset\mathcal{Q}$, we obtain
%%%%%%%%%%%%%%%%%%%%%%%%%%%%%%%%%%%%%%%%%%%%%%%%%%%%%%%%%%%%
\begin{cor}
The assignment $\omega\mapsto\omega$ and
$T_i\mapsto T_i(\underline{k})$
for $\omega\in\Omega$ and $i=0,\ldots,n$ uniquely extends to an injective
algebra homomorphism $\pi_{\underline{k}}:
H(\underline{\tau})\hookrightarrow \mathbb{D}_{\mathcal{R}}(W)$.
\end{cor}
%%%%%%%%%%%%%%%%%%%%%%%%%%%%%%%%%%%%%%%%%%%%%%%%%%%%%%%%%%%%%
We identify $H(\underline{\tau})$ with its image in 
$\mathbb{D}_{\mathcal{R}}(W)$
under the above algebra embedding $\pi_{\underline{k}}$.
In particular, we write $T_i$ for
$T_i(\underline{k})\in\mathbb{D}_{\mathcal{R}}(W)$
($i=0,\ldots,n$) if no confusion can arise.

We now proceed to define a reduction map $\beta$ which
maps difference-reflection operators to difference operators. 
This map is vital
for constructing the Macdonald difference operators 
from the representation $\pi_{\underline{k}}$,
as well as for constructing a large family of 
difference operators commuting with the Macdonald
difference operator.

Since $\mathbb{D}_{\mathcal{R}}(L^\prime)$ is a $W_0$-module algebra
by $W_0$-conjugation in $\mathbb{D}_{\mathcal{R}}(W)$,
\[D\mapsto wDw^{-1},\qquad D\in \mathbb{D}_{\mathcal{R}}(L^\prime),\,\,
w\in W_0,
\]
we have a canonical isomorphism
\begin{equation}\label{mult}
\mathbb{D}_{\mathcal{R}}(W)\simeq \mathbb{D}_{\mathcal{R}}(L^\prime)\# W_0
\end{equation}
of algebras. We now define the linear map
$\beta: \mathbb{D}_{\mathcal{R}}(W)\rightarrow 
\mathbb{D}_{\mathcal{R}}(L^\prime)$
by
\[\beta(D)=\sum_{w\in W_0}D_w,\qquad D\in \mathbb{D}_{\mathcal{R}}(W),
\]
where $D=\sum_{w\in W_0}D_ww$ ($D_w\in \mathbb{D}_{\mathcal{R}}(L^\prime)$)
is the unique decomposition of $D\in \mathbb{D}_{\mathcal{R}}(W)$ along the
isomorphism \eqref{mult}.
Let $H_0=H_0(\underline{\tau})\subset H(\underline{\tau})$ 
be the finite Hecke algebra, generated by
$T_j$ ($j=1,\ldots,n$). Consider the subalgebras
\begin{equation*}
\begin{split}
\mathbb{D}_\mathcal{R}(W)^{H_0}&=\{D\in
\mathbb{D}_\mathcal{R}(W)\,\,\, |\,\,\, \lbrack
D,h\rbrack=0\quad \forall\,h\in H_0 \},\\
\mathbb{D}_\mathcal{R}(L^\prime)^{W_0}&=
\{D\in\mathbb{D}_\mathcal{R}(L^\prime)\,\,\, 
|\,\,\, wDw^{-1}=D \quad \forall\,w\in W_0 \}.
\end{split}
\end{equation*}
Note that the center $Z(H(\underline{\tau}))$ of $H(\underline{\tau})$
is contained in $\mathbb{D}_{\mathcal{R}}(W)^{H_0}$. The following lemma, 
which is the
difference analog of \cite[Lemma 1.2.2]{HS}, 
gives for case {\bf a} a slight extension of \cite[Thm. 3.3]{Ch1}.

%%%%%%%%%%%%%%%%%%%%%%%%%%%%%%%%%%%%%%%%%%%%%%%%%%%%%%
\begin{lem}\label{resn}
The map $\beta$ restricts to an algebra homomorphism
\[\beta:\mathbb{D}_\mathcal{R}(W)^{H_0}\rightarrow
\mathbb{D}_\mathcal{R}(L^\prime)^{W_0}.
\]
In particular, the operators 
$\beta(D)$ \textup{(}$D\in\mathbb{D}_{\mathcal{R}}(W)^{H_0}$\textup{)}
preserve the subalgebra of $W_0$-invariant elements in $\mathcal{Q}$.
\end{lem}
%%%%%%%%%%%%%%%%%%%%%%%%%%%%%%%%%%%%%%%%%%%%%%%%%%%%%%%%%
\begin{proof}
Let $D\in\mathbb{D}_{\mathcal{R}}(W)^{H_0}$ and write
$D=\sum_{w\in W_0}D_ww$
with $D_w\in\mathbb{D}_{\mathcal{R}}(L^\prime)$. Fix $j\in\{1,\ldots,n\}$.
Since
\begin{equation*}
\begin{split}
\beta\bigl(Dc_{a_j}(z)(s_j-1)\bigr)&=0,\\
\beta\bigl(c_{a_j}(z)(s_j-1)D\bigr)&=
c_{a_j}(z)\bigl(s_j\beta(D)s_j-\beta(D)\bigr),
\end{split}
\end{equation*}
we obtain from the fact that $[D,T_j]=0$ in
$\mathbb{D}_{\mathcal{R}}(W)$,
\[0=\beta\bigl(\lbrack D,
c_{a_j}(z)(s_j-1)\rbrack\bigr)=
-c_{a_j}(z)\bigl(s_j\beta(D)s_j-\beta(D)\bigr),
\]
hence $s_j\beta(D)s_j=\beta(D)$.
We conclude that $\beta(D)\in\mathbb{D}_{\mathcal{R}}(L^\prime)^{W_0}$.

For $D,D^\prime\in \mathbb{D}_\mathcal{R}(W)^{H_0}$, written as
$D=\sum_{v\in W_0}D_vv$ and $D^\prime=\sum_{w\in W_0}D_w^\prime w$ with
$D_v,D_w^\prime\in\mathbb{D}_{\mathcal{R}}(L^\prime)$, we now have
\begin{equation*}
\begin{split}
\beta(DD^\prime)&=
\sum_{u,v\in W_0}D_vvD_{v^{-1}u}^\prime v^{-1}\\
&=\sum_{v\in W_0}D_vv\beta(D^\prime)v^{-1}\\
&=\sum_{v\in W_0}D_v\beta(D^\prime)=\beta(D)\beta(D^\prime),
\end{split}
\end{equation*}
hence $\beta$, restricted to $\mathbb{D}_{\mathcal{R}}(W)^{H_0}$,
is an algebra homomorphism.
\end{proof}
%%%%%%%%%%%%%%%%%%%%%%%%%%%%%%%%%%%%%%%%%%%%%%%%%%%%%%%%%%%%%%%%%%%
In the next subsection we describe Cherednik's commuting family of
difference operators $\beta(h)\in\mathbb{D}_{\mathcal{R}}(L^\prime)^{W_0}$ 
($h\in
Z(H(\underline{\tau}))$) and give the explicit expression of
the Macdonald difference operators inside
$\beta(Z(H(\underline{\tau})))$.

%%%%%%%%%%%%%%%%%%%%%%%%%%%%%%%%%%%%%%%%%%%%%%%%%%%%%%%%%%%%%%%%%%%%%%%%%%%
%%                                                                       %%
%%        The center $Z(H(\underline{\tau})$ and the Macdonald
%%                    difference operators                               %%
%%                                                                       %%
%%%%%%%%%%%%%%%%%%%%%%%%%%%%%%%%%%%%%%%%%%%%%%%%%%%%%%%%%%%%%%%%%%%%%%%%%%

\subsection{Cherednik-Macdonald commuting difference
operators}\label{CMsub}

Let $V_+\subset V$ be the open dominant Weyl chamber with respect to $R^+$,
\[V_+=\{v\in V \,\, | \,\, \langle v,\alpha\rangle>0 \quad
\forall\,\alpha\in R^+ \},
\]
and write $\overline{V}_+$ for its closure in $V$. We denote
\[L_{++}:=L\cap \overline{V}_+,\qquad
L_{++}^\prime:=L^\prime\cap\overline{V}_+
\]
for the cone of dominant elements in $L$ and $L^\prime$,
respectively. Set
\[Y^{\lambda^\prime}:=T_{t(\lambda^\prime)}\in
H(\underline{\tau}),\qquad \lambda^\prime\in L_{++}^\prime,
\]
and for arbitrary $\lambda^\prime\in L^\prime$ define
\[Y^{\lambda^\prime}=Y^{\mu^\prime}\bigl(Y^{\nu^\prime}\bigr)^{-1}\in
H(\underline{\tau})
\]
if $\lambda^\prime=\mu^\prime-\nu^\prime$ with $\mu^\prime,\nu^\prime\in
L_{++}^\prime$. The $Y^{\lambda^\prime}$ ($\lambda^\prime\in
L^\prime$) are well defined and satisfy
\[Y^0=1,\qquad
Y^{\lambda^\prime}Y^{\mu^\prime}=Y^{\lambda^\prime+\mu^\prime}
\]
for $\lambda^\prime,\mu^\prime\in L^\prime$. The
subalgebra $A^\prime(Y)$ of $H(\underline{\tau})$ spanned by the
$Y^{\lambda^\prime}$ ($\lambda^\prime\in L^\prime$) is isomorphic to
$A^\prime$ by the map $Y^{\lambda^\prime}\mapsto \xi^{\lambda^\prime}$
for $\lambda^\prime\in L^\prime$. The multiplication map induces
a linear isomorphism
\[H(\underline{\tau})=H_0(\underline{\tau})\otimes_{\mathbb{C}}A^\prime(Y)
\]
and
\[Z(H(\underline{\tau}))=A_0^\prime(Y),
\]
where $A_0^\prime$ is the subalgebra of
$W_0$-invariant elements in $A^\prime$ and $A_0^\prime(Y)\subset
H(\underline{\tau})$ is the corresponding subspace in $A^\prime(Y)$
via the above mentioned isomorphism $A^\prime\simeq A^\prime(Y)$.
We write $p(Y)\in A^\prime(Y)$ for the element
corresponding to $p(\xi)\in A^\prime$.

Similarly we write $A_0$ for the subalgebra of $W_0$-invariant
elements in $A$. For $\lambda\in L$ and $\lambda^\prime\in L^\prime$ we define
the monomial symmetric functions by
\[m_\lambda(z)=\sum_{\mu\in W_0\lambda}z^{\mu}\in A_0,\qquad
m_{\lambda^\prime}(\xi)=\sum_{\mu^\prime\in
W_0\lambda^\prime}\xi^{\mu^\prime}\in A_0^\prime.
\]
The functions $m_\lambda(z)$ ($\lambda\in L_{++}$)
(respectively $m_{\lambda^\prime}(\xi)$ ($\lambda^\prime\in
L^\prime_{++}$)) form a linear basis of $A_0$ (respectively $A_0^\prime$).
Correspondingly, the elements $m_{\lambda^\prime}(Y)$ ($\lambda^\prime\in
L_{++}^\prime$) form a linear basis of $Z(H(\underline{\tau}))$.

%%%%%%%%%%%%%%%%%%%%%%%%%%%%%%%%%%%%%%%%%%%%%%%%%%%%%%%%%%%%%%%%%%%%%%%%
\begin{defi}\label{Dpdef}
For $p(\xi)\in A_0^\prime$ we write
\[D_p:=\beta(p(Y))\in\mathbb{D}_{\mathcal{R}}(L^\prime)^{W_0}
\]
for the corresponding $W_0$-invariant difference operator.
For $p(\xi)=m_{\lambda^\prime}(\xi)\in A_0^\prime$ with 
$\lambda^\prime\in L^\prime$ we simplify the
notations by writing
\[D_{\lambda^\prime}:=D_{m_{\lambda^\prime}}=
\beta(m_{\lambda^\prime}(Y))\in\mathbb{D}_{\mathcal{R}}(L^\prime)^{W_0}.
\]
\end{defi}
%%%%%%%%%%%%%%%%%%%%%%%%%%%%%%%%%%%%%%%%%%%%%%%%%%%%%%%%%%%%%%%%%%%%%%%%%

Note that $D_{\lambda^\prime}\in\mathbb{D}_{\mathcal{R}}(L^\prime)^{W_0}$
only depends on the $W_0$-orbit $W_0\lambda^\prime$ of $\lambda^\prime\in
L^\prime$. In practice, it will be convenient to take the antidominant
representative of the orbit $W_0\lambda^\prime$. The following fact (cf.
\cite[Thm. 3.2]{Ch0} for case {\bf a}, and \cite{N} for case {\bf c}) 
is now immediate from Lemma \ref{resn}.

%%%%%%%%%%%%%%%%%%%%%%%%%%%%%%%%%%%%%%%%%%%%%%%%%%%%%%%%%%%%%%%%%%%%%%%%%
\begin{cor}\label{commfam}
The $W_0$-invariant difference operators
$D_p\in\mathbb{D}_{\mathcal{R}}(L^\prime)^{W_0}$ \textup{(}$p(\xi)\in
A_0^\prime$\textup{)} pair-wise commute.
\end{cor}
%%%%%%%%%%%%%%%%%%%%%%%%%%%%%%%%%%%%%%%%%%%%%%%%%%%%%%%%%%%%%%%%%%%%%%%%%%

We are now in a position to define the Macdonald difference operator as a 
difference operator $D_p$ for a special
choice of $p(\xi)\in A_0^\prime$.
We call a co-weight $\pi^\prime\in P(R^\vee)$ minuscule
(respectively quasi-minuscule) if it satisfies 
$|\langle \pi^\prime,\alpha\rangle|\leq 1$
for all $\alpha\in R$
(respectively if 
$\pi^\prime\in R^\vee$ and $|\langle\pi^\prime,\alpha\rangle|\leq 1$ for all
$\alpha\in R\setminus \{\pm \pi^{\prime\vee}\}$).
%%%%%%%%%%%%%%%%%%%%%%%%%%%%%%%%%%%%%%%%%%%%%%%%%%%%%%%%%%%%%%%%%%%%%%%%%%%%%
\begin{defi}\label{Macdonalddodef}
The operator $D_{\pi^\prime}\in\mathbb{D}_{\mathcal{R}}(L^\prime)^{W_0}$
with $\pi^\prime\in L^\prime\subseteq P(R^\vee)$ a nonzero minuscule or 
quasi-minuscule
co-weight, is called a Macdonald difference operator.
\end{defi}
%%%%%%%%%%%%%%%%%%%%%%%%%%%%%%%%%%%%%%%%%%%%%%%%%%%%%%%%%%%%%%%%%%%%%%%%%%%%
\begin{rema}
The Macdonald difference operators, 
which can be computed explicitly (as we shall
recall below), serve as the quantum Hamiltonians of relativistic
versions of the quantum trigonometric Calogero-Moser systems
associated to root systems in the sense of Ruijsenaars and
Schneider \cite{RS}. {}From this viewpoint Corollary \ref{commfam} 
reflects the
complete quantum integrability of the corresponding quantum
relativistic system, the higher quantum Hamiltonians being the
difference operators $D_p\in\mathbb{D}_{\mathcal{R}}(L^\prime)^{W_0}$ 
($p(\xi)\in A_0^\prime$).
\end{rema}
%%%%%%%%%%%%%%%%%%%%%%%%%%%%%%%%%%%%%%%%%%%%%%%%%%%%%%%%%%%%%%%%%%%%%%%%%%%

A key property of the difference operators $D_p$ ($p\in A_0^\prime$) is their
triangular action on $A_0$, which we now proceed to recall. Let
\begin{equation}\label{rho}
\rho_{\underline{k}^\prime}=\frac{1}{2}\sum_{\alpha\in
R^+}\underline{k}^\prime(\alpha^\vee)\alpha
\end{equation}
be the deformed half sum of positive roots
associated to the dual multiplicity label $\underline{k}^\prime$.
Let $\widetilde{p}(\xi)\in A^\prime$ for $p(\xi)\in A^\prime$ be the
associated
$\rho_{\underline{k}^\prime}$-twisted trigonometric Laurent
polynomial,
characterized by
\begin{equation}
\label{rhotwist}
\widetilde{p}(\lambda)=p(-\lambda-\rho_{\underline{k}^\prime}),\qquad
\forall\,\lambda\in V_{\mathbb{C}},
\end{equation}
where we use \eqref{planewave} to interpret $p(\xi)$ as function 
on $V_{\mathbb{C}}$.
Concretely, if $p(\xi)$ expands as
$p(\xi)=\sum_{\lambda^\prime\in
L^\prime}K(\lambda^\prime)\xi^{\lambda^\prime}$,
then $\widetilde{p}(\xi)=\sum_{\lambda^\prime\in
L^\prime}K(\lambda^{\prime})
q^{-\langle\rho_{\underline{k}^\prime},\lambda^\prime\rangle}\xi^{-\lambda^\prime}$.

By Theorem \ref{fundamental} and Lemma \ref{resn},
$D_p|_{A_0}$ is an endomorphism of $A_0$ for all $p(\xi)\in A_0^\prime$.
In fact, by \cite[(4.6.13)]{M} they are
triangular in the sense that for $p(\xi)\in A_0^\prime$ and
$\lambda\in L_{++}$,
\begin{equation}\label{triangular}
D_p(m_{\lambda}(z))=\widetilde{p}(\lambda)m_{\lambda}(z)+
\sum_{\mu\in L_{++}: \mu<\lambda}d_{\mu}(p)m_{\mu}(z)
\end{equation}
for some $d_{\mu}(p)\in\mathbb{C}$, with
$\leq$ the dominance order defined by $u\leq v$ for $u,v\in V_{\mathbb{C}}$
if $v-u\in\mathbb{Z}_{\geq 0}\Delta$.

We end this subsection by recalling the explicit form of the
Macdonald difference operators for each of the three cases {\bf a}, 
{\bf b} and {\bf c}.
The minuscule (respectively quasi-minuscule) co-weights form a
$W_0$-invariant subset of $P(R^\vee)$, so we first describe the
nonzero dominant minuscule and quasi-minuscule co-weights.

Let $\pi_i^\prime\in P(R^\vee)$ ($i=1,\ldots,n$) be the
fundamental co-weights with respect to the basis $\{\alpha_j\}_{j=1}^n$
of $R$ (so $\langle\pi_i^\prime,\alpha_j\rangle=\delta_{i,j}$ for 
$i,j=1,\ldots,n$
with $\delta_{i,j}$ the Kronecker delta function). Using the
expansion
\[\varphi=\sum_{j=1}^nm_j\alpha_j
\]
of the highest root $\varphi\in R$ 
in simple roots (in which all coefficients $m_j$ are strictly
positive integers), it follows that
\[\{\pi_j^\prime \,\, | \,\, j\in \{1,\ldots,n\}\,\,\hbox{ and
}\,\, m_j=1\}
\]
is the (possibly empty) set of nonzero dominant minuscule co-weights.
Macdonald difference operators of the first
type are the $D_{\pi^\prime}\in\mathbb{D}_{\mathcal{R}}(L^\prime)^{W_0}$
with $\pi^\prime=w_0\pi_j^\prime$ for some $j\in J_0$, where 
$w_0\in W_0$ is the longest Weyl
group element and $J_0$ is
the (possibly empty) set
\[J_0:=\{j\in\{1,\ldots,n\} \,\, | \,\, \pi_j^\prime\in
L^\prime\,\,\hbox{ and }\,\, m_j=1\}.
\]
%%%%%%%%%%%%%%%%%%%%%%%%%%%%%%%%%%%%%%%%%%%%%%%%%%%%%%%%%%%%%%%%%%%%%%%%%%%
\begin{eg}\label{example}
{\bf i)} For case {\bf a} and case {\bf b} we have $L^\prime=P(R^\vee)$, hence
$J_0$ pa\-ra\-me\-tri\-zes the nonzero dominant minuscule co-weights of
$R$. Using the classification of root systems $R$ one verifies that 
$J_0\not=\emptyset$
unless $R$ is of type $E_8$, $F_4$ or $G_2$.\\
{\bf ii)} For case {\bf c} there is precisely one nonzero
dominant minuscule co-weight, namely the
fundamental co-weight $\pi_j^\prime\in P(R^\vee)$ 
corresponding to the unique long simple root
$\alpha_j$ from the basis $\{\alpha_1,\ldots,\alpha_n\}$ of $R$.
But $\pi_j^\prime\not\in L^\prime=Q(R^\vee)$, hence $J_0=\emptyset$.
\end{eg}
%%%%%%%%%%%%%%%%%%%%%%%%%%%%%%%%%%%%%%%%%%%%%%%%%%%%%%%%%%%%%%%%%%%%%%%%%%%%%

It is easy to show that $\varphi^\vee$ is the only nonzero dominant 
quasi-minuscule co-weight.
Since a quasi-minuscule co-weight lies in $R^\vee$ by definition,
it automatically lies in $L^\prime$.
The Macdonald difference operator of the second type now is
$D_{-\varphi^\vee}\in\mathbb{D}_{\mathcal{R}}(L^\prime)^{W_0}$, which thus exists
for all the three cases {\bf a}, {\bf b} and {\bf c} under
consideration.

For $\lambda^\prime\in L^\prime$ we 
write $W_{0,\lambda^\prime}\subseteq W_0$ for the isotropy subgroup of
$\lambda^\prime$ and $W_0^{\lambda^\prime}$
for a complete set of representatives of $W_0/W_{0,\lambda^\prime}$.
%%%%%%%%%%%%%%%%%%%%%%%%%%%%%%%%%%%%%%%%%%%%%%%%%%%%%%%%%%%%%%
\begin{prop}\label{formexplicit}
Let $\pi^\prime\in L^\prime$ be a nonzero minuscule
or quasi-minuscule co-weight. Then
\[D_{\pi^\prime}=m_{\pi^\prime}(-\rho_{\underline{k}^\prime})+
\sum_{w\in
W_0^{\pi^\prime}}(wf_{\pi^\prime})(z)\bigl(t(w\pi^\prime)-1\bigr),
\]
with $f_{\pi^\prime}(z)\in\mathcal{R}$ given by
\[f_{\pi^\prime}(z)=\prod_{a\in S_1(t(-\pi^\prime))}c_a(z)
\]
and with $S_1(w)=S_1^+\cap w^{-1}S_1^-$ for $w\in W$.
If $\pi^\prime$ is minuscule then
\[D_{\pi^\prime}=\sum_{w\in
W_0^{\pi^\prime}}(wf_{\pi^\prime})(z)t(w\pi^\prime).
\]
\end{prop}
%%%%%%%%%%%%%%%%%%%%%%%%%%%%%%%%%%%%%%%%%%%%%%%%%%%%%%%%%%%%%%%%%
\begin{proof}
By \cite[(4.4.12)]{M} and the concluding paragraph of \cite[\S
4.4]{M} we have
\begin{equation}\label{initial}
D_{\pi^\prime}=g(z)+\sum_{w\in
W_0^{\pi^\prime}}(wf_{\pi^\prime})(z)t(w\pi^\prime)
\end{equation}
for some $g(z)\in\mathcal{R}$, and $g(z)$ is zero if 
$\pi^\prime$ is minuscule.
This yields the second formula for $D_{\pi^\prime}$.

To prove the first formula for $D_{\pi^\prime}$ we use \eqref{triangular} 
for $\lambda=0$,
which gives the expression
\[g(z)=m_{\pi^\prime}(-\rho_{\underline{k}^\prime})-\sum_{w\in
W_0^{\pi^\prime}}(wf_{\pi^\prime})(z).
\]
Combined with \eqref{initial} we obtain the first formula for $D_{\pi^\prime}$.
\end{proof}
%%%%%%%%%%%%%%%%%%%%%%%%%%%%%%%%%%%%%%%%%%%%%%%%%%%%%%%%%%%%%%%%%%
Recall that for case {\bf a}, $S_1=S(R)$ and $R^{\prime\vee}=R$,
while for case {\bf b} and {\bf c}, $S_1=S(R)^\vee$ and
$R^{\prime\vee}=R^\vee$. The following lemma follows now 
by a direct computation.
%%%%%%%%%%%%%%%%%%%%%%%%%%%%%%%%%%%%%%%%%%%%%%%%%%%%%%%%%%%%%%%%%
\begin{lem}\label{crosspre}
{\bf i)} For $\pi^\prime=w_0\pi_j^\prime$ with $j\in J_0$ we have
\[S_1(t(-\pi^\prime))=\{\alpha\in R^{\prime\vee} \,\, | \,\,
\langle \pi^\prime,\alpha\rangle<0\}.
\]
{\bf ii)} For $\pi^\prime=-\varphi^\vee$ we have
\[S_1(t(-\pi^\prime))=\{\alpha\in R^{\prime\vee}\,\, | \,\,
\langle\pi^\prime,\alpha\rangle<0\}\cup \{(\varphi+c)^\vee\},
\]
where for case {\bf a}, we note that $(\varphi+c)^\vee=\varphi+c\in S(R)=S_1$,
since $\varphi$ has squared norm equal to two by convention.
\end{lem}
%%%%%%%%%%%%%%%%%%%%%%%%%%%%%%%%%%%%%%%%%%%%%%%%%%%%%%%%%%%%%%%%%%%
We conclude by writing out the coefficients $f_{\pi^\prime}(z)\in\mathcal{R}$
explicitly. In the list below, $\pi^\prime\in L^\prime$ is a
nonzero anti-dominant minuscule or quasi-minuscule co-weight, 
so $\pi^\prime=w_0\pi_j^\prime$
($j\in J_0$) or $\pi^\prime=-\varphi^\vee$. 
We use that $z^{\varphi+c}=qz^\varphi$
for case {\bf a}, $z^{(\varphi+c)^\vee}=qz^{\varphi^\vee}$ for case
{\bf b} (since $\|\varphi\|^2=2$), and 
$z^{(\varphi+c)^\vee}=q^{\frac{1}{2}}z^{\varphi^\vee}$
for case {\bf c} (since $\|\varphi\|^2=4$).
We furthermore use that
\begin{equation*}
\{\beta\in R \,\, | \,\, \langle\pi^\prime,\beta\rangle=-2\}=
\begin{cases}
\emptyset,\qquad &\hbox{if }\, \pi^\prime=w_0\pi_j^\prime\,\,\, (j\in J_0),\\
\{\varphi\},\qquad &\hbox{if }\, \pi^\prime=-\varphi^\vee,
\end{cases}
\end{equation*}
hence the second factor in the expression of $f_{\pi^\prime}(z)$
for case {\bf a} and case {\bf b} below is one if $\pi^\prime$ is 
minuscule, while it is simply the product
over the singleton $\{\varphi\}$ if $\pi^\prime$ is quasi-minuscule
(cf. Chalykh's \cite[\S 2.2]{C} notational conventions).
%%%%%%%%%%%%%%%%%%%%%%%%%%%%%%%%%%%%%%%%%%%%%%%%%%%%%%%%%%%%%%%%%%%%
\begin{enumerate}
\item[Case {\bf a.}] We have
\[f_{\pi^\prime}(z)=\prod_{\stackrel{\alpha\in
R:}{\langle\pi^\prime,\alpha\rangle<0}}
\frac{(1-\tau_\alpha^2z^\alpha)}{\tau_\alpha(1-z^\alpha)}
\prod_{\stackrel{\beta\in
R:}{\langle\pi^\prime,\beta\rangle=-2}}
\frac{(1-q\tau_\beta^2z^\beta)}{\tau_\beta(1-qz^\beta)}.
\]
\item[Case {\bf b.}] We have
\[f_{\pi^\prime}(z)=\prod_{\stackrel{\alpha\in
R:}{\langle\pi^\prime,\alpha\rangle<0}}
\frac{(1-\tau_{\alpha^\vee}^2z^{\alpha^\vee})}
{\tau_{\alpha^\vee}(1-z^{\alpha^\vee})}
\prod_{\stackrel{\beta\in
R:}{\langle\pi^\prime,\beta\rangle=-2}}
\frac{(1-q\tau_{\beta^\vee}^2z^{\beta^\vee})}
{\tau_{\beta^\vee}(1-qz^{\beta^\vee})}.
\]
\item[Case {\bf c.}] For this case, 
$\pi^\prime=-\varphi^\vee$ is quasi-minuscule
as remarked before. Furthermore, 
$S_1(t(-\pi^\prime))\cap\mathcal{O}_1=\{\varphi^\vee\}$,
$S_1(t(-\pi^\prime))\cap\mathcal{O}_3=\{(\varphi+c)^\vee\}$, 
and the remaining
$2(n-1)$ elements in $S_1(t(-\pi^\prime))$ are 
from the orbit $\mathcal{O}_5$.
To make the expression for $f_{\pi^\prime}$ as explicit 
as possible, we choose without
loss of generality the particular basis
\[\{\alpha_1,\ldots,\alpha_n\}=
\{\epsilon_1-\epsilon_2,\ldots,\epsilon_{n-1}-\epsilon_n,2\epsilon_n\}
\]
for $R=R_C$, so that the corresponding highest root is
$\varphi=2\epsilon_1$ and $\pi^\prime=-\varphi^\vee=-\epsilon_1$. Then
\begin{equation*}
\begin{split}
f_{\pi^\prime}(z)&=q^{-\kappa_1'-(n-1)\kappa_5}\prod_{j=2}^n
\frac{(1-q^{\kappa_5}z^{\epsilon_1+\epsilon_j})(1-q^{\kappa_5}z^{\epsilon_1-\epsilon_j})}
{(1-z^{\epsilon_1+\epsilon_j})(1-z^{\epsilon_1-\epsilon_j})}\\
&\times\frac{(1-q^{\kappa_1}z^{\epsilon_1})(1+q^{\kappa_2}z^{\epsilon_1})
(1-q^{\frac{1}{2}+\kappa_3}z^{\epsilon_1})(1+q^{\frac{1}{2}+\kappa_4}z^{\epsilon_1})}
{(1-z^{2\epsilon_1})(1-qz^{2\epsilon_1})}
\end{split}
\end{equation*}
and the resulting difference operator 
$D_{\pi^\prime}\in\mathbb{D}_{\mathcal{R}}(L^\prime)^{W_0}$
is Koornwinder's \cite{Ko} difference operator, see also \cite{N}.
\end{enumerate}
%%%%%%%%%%%%%%%%%%%%%%%%%%%%%%%%%%%%%%%%%%%%%%%%%%%%%%%%%%%%%%%

%%%%%%%%%%%%%%%%%%%%%%%%%%%%%%%%%%%%%%%%%%%%%%%%%%%%%%%%%%%%%%%%%%%%%%%
%%                                                                   %%
%%           The Harish-Chandra homomorphism                         %%
%%                                                                   %%
%%%%%%%%%%%%%%%%%%%%%%%%%%%%%%%%%%%%%%%%%%%%%%%%%%%%%%%%%%%%%%%%%%%%%%%
\subsection{The Harish-Chandra homomorphism}\label{HChom}
In this subsection we describe the asymptotics of the difference
operators $D_p\in\mathbb{D}_{\mathcal{R}}(L^\prime)^{W_0}$ ($p\in A_0^\prime$) 
deep in
the negative Weyl chamber
\[V_-:=\{ v\in V \,\, | \,\, \langle v,\alpha\rangle<0 \qquad
\forall\,\alpha\in R^+ \}.
\]
We give an algebraic formalization 
in terms of a constant term map $\gamma(\underline{k})$
(see \cite[\S 1.2]{HS} for this notion in the trigonometric 
differential degeneration),
whose definition below is based on the 
fact that $z^{-x}(v)=q^{-\langle x,v\rangle}\rightarrow 0$
for $x\in\mathbb{Z}_{\geq 0}\Delta\setminus\{0\}$ if 
$\langle v,\alpha\rangle\rightarrow-\infty$
for all $\alpha\in R^+$. In this algebraic formulation we
incorporate a $\rho_{\underline{k}^\prime}$-shift and we
identify the algebra $\mathbb{C}[t(L^\prime)]$ of constant
coefficient difference operators with $A^\prime$ by 
$t(\lambda^\prime)\leftrightarrow \xi^{\lambda^\prime}$
for $\lambda^\prime\in L^\prime$. For case {\bf a}, the results
discussed in this section are essentially due to Cherednik, 
see e.g. \cite[\S 3]{Ch1}.

%%%%%%%%%%%%%%%%%%%%%%%%%%%%%%%%%%%%%%%%%%%%%%%%%%%%%%%%%%%%%%%%%%%%%%%%%%
\begin{lem}\label{constanttermmap}
For
\[
D=\sum_{\lambda^\prime\in
L^\prime}\Bigl(\sum_{x\in\mathbb{Z}_{\geq
0}\Delta}C_x(\lambda^\prime)z^{-x}\Bigr)t(\lambda^\prime)
\in\mathbb{D}_{\mathbb{C}[[z^{-\Delta}]]}(L^\prime)
\]
with $\{C_x(\lambda^\prime)\}_{x\in \mathbb{Z}_{\geq 0}\Delta}\subset\mathbb{C}$
the zero set for all but finitely many $\lambda^\prime\in
L^\prime$, we define the constant term $\gamma(\underline{k})(D)\in A^\prime$ 
of $D$
by
\begin{equation}\label{gammak}
\gamma(\underline{k})(D)=\sum_{\lambda^\prime\in
L^\prime}C_0(\lambda^\prime)
q^{\langle\rho_{\underline{k}^\prime},\lambda^\prime\rangle}\xi^{\lambda^\prime}.
\end{equation}
The resulting linear map 
$\gamma(\underline{k}): \mathbb{D}_{\mathbb{C}[[z^{-\Delta}]]}(L^\prime)
\rightarrow A^\prime$ is an algebra homomorphism, 
called the Harish-Chandra homomorphism.
\end{lem}
%%%%%%%%%%%%%%%%%%%%%%%%%%%%%%%%%%%%%%%%%%%%%%%%%%%%%%%%%%%%%%%%%%%%%%%%
\begin{proof}
This follows from a direct computation.
\end{proof}
%%%%%%%%%%%%%%%%%%%%%%%%%%%%%%%%%%%%%%%%%%%%%%%%%%%%%%%%%%%%%%%%%%%%%%%
Recall the
$\mathbb{D}_{\mathbb{C}[[z^{-\Delta}]]}(L^\prime)$-module
$\overline{M}$ from Lemma \ref{M} and the 
$\rho_{\underline{k}^\prime}$-twist \eqref{rhotwist}.

%%%%%%%%%%%%%%%%%%%%%%%%%%%%%%%%%%%%%%%%%%%%%%%%%%%%%%%%%%%
\begin{lem}\label{ltlem}
For $u\in V_{\mathbb{C}}$ and 
$D\in \mathbb{D}_{\mathbb{C}[[z^{-\Delta}]]}(L^\prime)$
we have
\[D(z^u)=\widetilde{\bigl(\gamma(\underline{k})(D)\bigr)}(u)z^u
+\sum_{x\in \mathbb{Z}_{\geq 0}\Delta\setminus \{0\}}K_x(u)z^{u-x}\in
\overline{M}
\]
for certain $K_x(u)\in\mathbb{C}$.
\end{lem}
%%%%%%%%%%%%%%%%%%%%%%%%%%%%%%%%%%%%%%%%%%%%%%%%%%%%%%%%%%
\begin{proof}
Let $D\in\mathbb{D}_{\mathbb{C}[[z^{-\Delta}]]}(L^\prime)$, 
written out explicitly as
\begin{equation}\label{Dgeneric}
D=\sum_{\lambda^\prime\in
L^\prime}\Bigl(\sum_{x\in\mathbb{Z}_{\geq
0}\Delta}C_x(\lambda^\prime)z^{-x}\Bigr)
t(\lambda^\prime)\in\mathbb{D}_{\mathbb{C}[[z^{-\Delta}]]}(L^\prime)
\end{equation}
with $\{C_x(\lambda^\prime)\}_{x\in \mathbb{Z}_{\geq 0}\Delta}\subset\mathbb{C}$
the zero set for all but finitely many $\lambda^\prime\in
L^\prime$. Then
\[D(z^u)=\sum_{x\in\mathbb{Z}_{\geq 0}\Delta}K_{x}(u)z^{u-x}\in\overline{M}
\]
by Lemma \ref{M}, with coefficients
\[K_x(u)=\sum_{\lambda^\prime\in
L^\prime}C_{x}(\lambda^\prime)q^{-\langle\lambda^\prime,u\rangle},\qquad
x\in\mathbb{Z}_{\geq 0}\Delta.
\]
On the other hand,
\[
\widetilde{\bigl(\gamma(\underline{k})(D)\bigr)}(u)=\sum_{\lambda^\prime\in
L^\prime}C_0(\lambda^\prime)q^{-\langle\lambda^\prime,u\rangle}
\]
by \eqref{planewave}, \eqref{rhotwist} and \eqref{gammak}, hence
$K_0(u)=\widetilde{\bigl(\gamma(\underline{k})(D)\bigr)}(u)$.
\end{proof}
%%%%%%%%%%%%%%%%%%%%%%%%%%%%%%%%%%%%%%%%%%%%%%%%%%%%%%%%%%%%%%%%%%%%

We are now in the position to explicitly compute the
constant terms $\gamma(\underline{k})(D_p)$ for the commuting family of
difference operators
$D_p\in\mathbb{D}_{\mathcal{R}}(L^\prime)^{W_0}
\subset\mathbb{D}_{\mathbb{C}[[z^{-\Delta}]]}(L^\prime)$
($p(\xi)\in A_0^\prime$).
%%%%%%%%%%%%%%%%%%%%%%%%%%%%%%%%%%%%%%%%%%%%%%%%%%%%%%%%%%%%%%%%%%%%%%%
\begin{prop}\label{gammainverse}
We have $\bigl(\gamma(\underline{k})(D_p)\bigr)(\xi)=p(\xi)$
for $p(\xi)\in A_0^\prime$.
In particular, the restriction of the algebra homomorphism
\[\gamma(\underline{k})\circ\beta:
\mathbb{D}_{\mathcal{R}}(W)^{H_0}\rightarrow A^\prime
\]
to $A_0^\prime(Y)\subset \mathbb{D}_{\mathcal{R}}(W)^{H_0}$
is the algebra isomorphism $A_0^\prime(Y)\overset{\sim}{\longrightarrow}
A_0^\prime$ which maps $p(Y)$ to $p(\xi)$ \textup{(}$p\in A_0^\prime$\textup{)}.
\end{prop}
%%%%%%%%%%%%%%%%%%%%%%%%%%%%%%%%%%%%%%%%%%%%%%%%%%%%%%%%%%%%%%%%%%%%%%%%%%%%
\begin{proof}
We embed $A$ as subspace in $\overline{M}$ as the complex subspace spanned by
the monomials $z^{\mu}\in\overline{M}$ ($\mu\in L$).
The restriction of the action of 
$D_p\in\hbox{End}_{\mathbb{C}}(\overline{M})$ ($p(\xi)\in A_0^\prime$)
to $A_0$ coincides with the action of $D_p$ on $A_0$ as described in 
Subsections 2.3 and 2.4.

Fix $p\in A_0^\prime$. For $\lambda\in L_{++}$ 
and $\mu\in W_0\lambda$ we have $\mu\leq\lambda$
(i.e. $\lambda-\mu\in\mathbb{Z}_{\geq 0}\Delta$), 
hence \eqref{triangular} and the fact
that $\mathbb{C}[[z^{-\Delta}]]z^u\subset\overline{M}$ is a
$\mathbb{D}_{\mathbb{C}[[z^{-\Delta}]]}(L^\prime)$-submodule 
for all $u\in V_{\mathbb{C}}$
imply
\[D_p(z^\lambda)=\widetilde{p}(\lambda)z^{\lambda}+
\sum_{\mu\in L: \mu<\lambda}K_\mu(\lambda)z^\mu,\qquad \forall\,\lambda\in
L_{++}
\]
for certain coefficients $K_\mu(\lambda)\in\mathbb{C}$. Combined with Lemma
\ref{ltlem} we conclude that
\[\widetilde{\bigl(\gamma(\underline{k})(D_p)\bigr)}(\lambda)=
\widetilde{p}(\lambda),
\qquad \forall\,\lambda\in L_{++}.
\]
This implies that $\bigl(\gamma(\underline{k})(D_p)\bigr)(\xi)=p(\xi)$ 
in $A_0^\prime$.
\end{proof}
%%%%%%%%%%%%%%%%%%%%%%%%%%%%%%%%%%%%%%%%%%%%%%%%%%%%%%%%%%%%%%%%%%%%%%%%%%%%%%
Recall the explicit form of the Macdonald difference operators $D_{\pi^\prime}$
from Proposition \ref{formexplicit}, where $\pi^\prime=w_0\pi_j^\prime$
($j\in J_0$) or $\pi^\prime=-\varphi^\vee$. The following technical result
will be used in the next section in the analysis of the centralizer 
of $D_{\pi^\prime}$
in $\mathbb{D}_{\mathcal{R}}(L^\prime)^{W_0}$.
%%%%%%%%%%%%%%%%%%%%%%%%%%%%%%%%%%%%%%%%%%%%%%%%%%%%%%%%%%%%%%%%%%%%%%%%%%%%
\begin{cor}\label{Macdonaldreturn}
For $w\in W_0^{\pi^\prime}$ we have
\[\gamma(\underline{k})(wf_{\pi^\prime})=
q^{-\langle\rho_{\underline{k}^\prime},w\pi^\prime\rangle}.
\]
In particular,
\[D_{\pi^\prime}=\sum_{\mu^\prime\in
W_0\pi^\prime}q^{-\langle\rho_{\underline{k}^\prime},\mu^\prime\rangle}t(\mu^\prime)+
\sum_{\mu^\prime\in W_0\pi^\prime}g_{\mu^\prime}(z)\bigl(t(\mu^\prime)-1\bigr)\]
where for $w\in W_0^{\pi^\prime}$,
\[g_{w\pi^\prime}(z):=(wf_{\pi^\prime})(z)-\gamma(\underline{k})(wf_{\pi^\prime})=
\sum_{x\in\mathbb{Z}_{\geq
0}\Delta\setminus\{0\}}K_x(w\pi^\prime)z^{-x}
\]
for certain $K_x(w\pi^\prime)\in\mathbb{C}$.
\end{cor}

%%%%%%%%%%%%%%%%%%%%%%%%%%%%%%%%%%%%%%%%%%%%%%%%%%%%%%%%%%%
%%                                                       %%
%%               Centralizer algebras                    %%
%%                                                       %%
%%%%%%%%%%%%%%%%%%%%%%%%%%%%%%%%%%%%%%%%%%%%%%%%%%%%%%%%%%%

\section{Centralizers of Macdonald difference
operators}\label{section2}

In this section we fix
a nonzero anti-dominant minuscule or quasi-minuscule co-weight
$\pi^\prime\in L^\prime$. For the associated Macdonald difference
operator $D_{\pi^\prime}\in\mathbb{D}_{\mathcal{R}}(L^\prime)^{W_0}$
we analyze the centralizer algebra
$\mathbb{D}_{\mathcal{R}}(L^\prime)^{W_0,D_{\pi^\prime}}$
(respectively $\mathbb{D}_{\mathcal{R}}(L^\prime)^{D_{\pi^\prime}}$) consisting of
the difference operators
$D\in\mathbb{D}_{\mathcal{R}}(L^\prime)^{W_0}$ 
(respectively $D\in\mathbb{D}_{\mathcal{R}}(L^\prime)$)
which commute with $D_{\pi^\prime}$.

We already observed that
\[\beta(A_0^\prime(Y))\subseteq
\mathbb{D}_{\mathcal{R}}(L^\prime)^{W_0,D_{\pi^\prime}}\subseteq
\mathbb{D}_{\mathcal{R}}(L^\prime)^{D_{\pi^\prime}}.
\]
We show in this section that the first inclusion is an equality.
As observed in e.g. \cite{ES} and \cite{C}, the second
inclusion is strict for special values of the multiplicity labels,
in which case one speaks of algebraic integrability. We analyze the
second inclusion using a simple symmetrization procedure for difference
operators.

%%%%%%%%%%%%%%%%%%%%%%%%%%%%%%%%%%%%%%%%%%%%%%%%%%%%%%%%%%%%%%%%
%%                                                            %%
%%          Commutatitivity                                   %%
%%                                                            %%
%%%%%%%%%%%%%%%%%%%%%%%%%%%%%%%%%%%%%%%%%%%%%%%%%%%%%%%%%%%%%%%%

\subsection{Commutativity}\label{Comm}
We start by showing that $\mathbb{D}_{\mathcal{R}}(L^\prime)^{D_{\pi^\prime}}$ 
is commutative.
Denote $\mathbb{D}_{\mathbb{C}[[z^{-\Delta}]]}(L^\prime)^{D_{\pi^\prime}}$
for the centralizer algebra of $D_{\pi^\prime}$ in
$\mathbb{D}_{\mathbb{C}[[z^{-\Delta}]]}(L^\prime)$. It
contains $\mathbb{D}_{\mathcal{R}}(L^\prime)^{D_{\pi^\prime}}$
as a subalgebra.

%%%%%%%%%%%%%%%%%%%%%%%%%%%%%%%%%%%%%%%%%%%%%%%%%%%%%%%%%%%%%%%%%%%
\begin{prop}\label{Commutativityprop}
The Harish-Chandra homomorphism $\gamma(\underline{k})$ 
restricts to an injective algebra
homomorphism
\[\gamma(\underline{k}):
\mathbb{D}_{\mathbb{C}[[z^{-\Delta}]]}(L^\prime)^{D_{\pi^\prime}}
\hookrightarrow A^\prime.
\]
In particular, $\mathbb{D}_{\mathbb{C}[[z^{-\Delta}]]}(L^\prime)^{D_{\pi^\prime}}$
and $\mathbb{D}_{\mathcal{R}}(L^\prime)^{D_{\pi^\prime}}$ are commutative.
\end{prop}
%%%%%%%%%%%%%%%%%%%%%%%%%%%%%%%%%%%%%%%%%%%%%%%%%%%%%%%%%%%%%%%%%%%%
\begin{proof}
We express $D\in\mathbb{D}_{\mathbb{C}[[z^{-\Delta}]]}(L^\prime)$
as \eqref{Dgeneric}. We can then 
explicitly compute the $z^{-x}\mathbb{C}[t(L^\prime)]$-term of the commutant
$\lbrack D, 
D_{\pi^\prime}\rbrack\in\mathbb{D}_{\mathbb{C}[[z^{-\Delta}]]}(L^\prime)$
($x\in\mathbb{Z}_{\geq 0}\Delta$) using the expression
from Corollary \ref{Macdonaldreturn} for the Macdonald 
difference operator $D_{\pi^\prime}$.
It follows that
$D\in\mathbb{D}_{\mathbb{C}[[z^{-\Delta}]]}(L^\prime)^{D_{\pi^\prime}}$ if and
only if
\begin{equation}\label{recuroperator}
\begin{split}
&\Bigl(\sum_{\mu^\prime\in W_0\pi^\prime}
q^{-\langle\rho_{\underline{k}^\prime},\mu^\prime\rangle}
\bigl(1-q^{\langle \mu^\prime,x\rangle}\bigr)t(\mu^\prime)\Bigr)
\Bigl(\sum_{\lambda^\prime\in L^{\prime}}
C_x(\lambda^\prime)t(\lambda^\prime)\Bigr)\\
&\qquad\qquad=\sum_{\stackrel{\lambda^\prime\in L^{\prime}}{\mu^\prime\in
W_0\pi^\prime}}\sum_{0<y\leq x}
C_{x-y}(\lambda^\prime)K_y(\mu^\prime)\bigl(q^{\langle\mu^\prime,x-y\rangle}-
q^{\langle\lambda^\prime,y\rangle}\bigr)t(\mu^\prime+\lambda^\prime)\\
&\qquad\qquad\quad+\sum_{\stackrel{\lambda^\prime\in L^{\prime}}{\mu^\prime\in
W_0\pi^\prime}}
\sum_{0<y\leq x}
C_{x-y}(\lambda^\prime)K_y(\mu^\prime)
\bigl(q^{\langle\lambda^\prime,y\rangle}-1\bigr)t(\lambda^\prime)
\end{split}
\end{equation}
in $\mathbb{C}[t(L^\prime)]$ for all $x\in\mathbb{Z}_{\geq 0}\Delta$.
Observe that the first factor
\[\sum_{\mu^\prime\in W_0\pi^\prime}q^{-\langle\rho_{k^\prime},\mu^\prime\rangle}
\bigl(1-q^{\langle
\mu^\prime,x\rangle}\bigr)t(\mu^\prime)
\]
in the left hand side of \eqref{recuroperator} is nonzero for all 
$x\in\mathbb{Z}_{\geq
0}\Delta\setminus\{0\}$ since $W_0$ acts irreducibly on $V$.

Let now $D\in\mathbb{D}_{\mathbb{C}[[z^{-\Delta}]]}(L^\prime)^{D_{\pi^\prime}}$
with $\gamma(\underline{k})(D)=0$. By the definition of the
Harish-Chandra homomorphism (see Lemma \ref{constanttermmap}), we
obtain
\begin{equation}\label{kern}
C_0(\lambda^\prime)=0,\qquad \forall\,\lambda^\prime\in
L^\prime.
\end{equation}
By induction on the height of $x\in\mathbb{Z}_{\geq 0}\Delta$
(with respect to the simple roots $\Delta$), it follows from
\eqref{recuroperator} and \eqref{kern} that $C_x(\lambda^\prime)=0$
for all $\lambda^\prime\in L^\prime$ and all $x\in\mathbb{Z}_{\geq
0}\Delta$. Hence $D=0$, as desired.
\end{proof}
%%%%%%%%%%%%%%%%%%%%%%%%%%%%%%%%%%%%%%%%%%%%%%%%%%%%%%%%%%%%%%%%%%%%%%%%%%%%

%%%%%%%%%%%%%%%%%%%%%%%%%%%%%%%%%%%%%%%%%%%%%%%%%%%%%%%%%%%%%%%%
%%                                                            %%
%%        Symmetrization of difference operators              %%
%%                                                            %%
%%%%%%%%%%%%%%%%%%%%%%%%%%%%%%%%%%%%%%%%%%%%%%%%%%%%%%%%%%%%%%%%

\subsection{Symmetrization of difference operators}\label{Symmsub}

Define a surjective linear projection
$\pi:\mathbb{D}_{\mathcal{R}}(L^\prime)
\rightarrow\mathbb{D}_{\mathcal{R}}(L^\prime)^{W_0}$
by
\[\pi(D):=\frac{1}{\#W_0}\sum_{w\in W_0}wDw^{-1},\qquad
D\in\mathbb{D}_{\mathcal{R}}(L^\prime).
\]
Observe that
$\pi(DD^\prime)=\pi(D)D^\prime$ and $\pi(D^\prime D)=D^\prime\pi(D)$
for $D\in\mathbb{D}_{\mathcal{R}}(L^\prime)$ and
$D^\prime\in\mathbb{D}_{\mathcal{R}}(L^\prime)^{W_0}$.
Let $e$ be the trivial idempotent of $\mathbb{C}[W_0]$,
\[e=\frac{1}{\#W_0}\sum_{w\in
W_0}w\in\mathbb{C}[W_0]\subset\mathbb{D}_{\mathcal{R}}(W).
\]
%%%%%%%%%%%%%%%%%%%%%%%%%%%%%%%%%%%%%%%%%%%%%%%%%%%%%%%%%%%%%%%%%%%%%%%%%%%%
\begin{lem}
{\bf i)} We have $eDe=\pi(D)e$ for $D\in\mathbb{D}_{\mathcal{R}}(L^\prime)$.\\
{\bf ii)} The projection map $\pi$ restricts to a surjection
\[\pi: \mathbb{D}_{\mathcal{R}}(L^\prime)^{D_{\pi^\prime}}\rightarrow
\mathbb{D}_{\mathcal{R}}(L^\prime)^{W_0,D_{\pi^\prime}}.
\]
\end{lem}
%%%%%%%%%%%%%%%%%%%%%%%%%%%%%%%%%%%%%%%%%%%%%%%%%%%%%%%%%%%%%%%%%%%%%%%%%%%%
\begin{proof}
{\bf i)} Since $we=e$ for $w\in W_0$ we have
\[
eDe=\frac{1}{\#W_0}\sum_{w\in W_0}wDe=\frac{1}{\#W_0}\sum_{w\in
W_0}wDw^{-1}e=\pi(D)e.
\]
{\bf ii)} Note that
$D_{\pi^\prime}\in
\beta(A_0^\prime(Y))\subseteq\mathbb{D}_{\mathcal{R}}(L^\prime)^{W_0}$,
hence
\[\lbrack \pi(D),D_{\pi^\prime}\rbrack= \pi\bigl(\lbrack
D,D_{\pi^\prime}\rbrack\bigr),\qquad
D\in\mathbb{D}_{\mathcal{R}}(L^\prime),
\]
which immediately leads to the desired result.
\end{proof}
%%%%%%%%%%%%%%%%%%%%%%%%%%%%%%%%%%%%%%%%%%%%%%%%%%%%%%%%%%%%%%%%%%%%%%%%%%%
\begin{prop}\label{centralizergen}
Let $D\in\mathbb{D}_{\mathcal{R}}(L^\prime)^{D_{\pi^\prime}}$.
If $A_0\subset\overline{M}$ is invariant under the action of
the difference operator $D$ on $\overline{M}$,
then $D\in\mathbb{D}_{\mathcal{R}}(L^\prime)^{W_0, D_{\pi^\prime}}$.
\end{prop}
%%%%%%%%%%%%%%%%%%%%%%%%%%%%%%%%%%%%%%%%%%%%%%%%%%%%%%%%%%%%%%%%%%%%%%%%%%%%
\begin{proof}
Let
$D\in\mathbb{D}_{\mathcal{R}}(L^\prime)^{D_{\pi^\prime}}$
satisfying $D(A_0)\subseteq A_0$, so that
\[
D(m_\lambda(z))=\pi(D)(m_\lambda(z)),\qquad \forall\,\lambda\in L_{++}.
\]
Equating the $z^\lambda$-coefficient using Lemma \ref{ltlem} we
obtain
\[
\widetilde{\bigl(\gamma(\underline{k})(D)\bigr)}(\lambda)=
\widetilde{\bigl(\gamma(\underline{k})(\pi(D))\bigr)}(\lambda),
\qquad\forall\,\lambda\in
L_{++},
\]
compare with the proof of Proposition \ref{gammainverse}.
We conclude that
$\gamma(\underline{k})(D)=\gamma(\underline{k})(\pi(D))$, hence
\[D=\pi(D)\in\mathbb{D}_{\mathcal{R}}(L^\prime)^{W_0,D_{\pi^\prime}}
\]
since the Harish-Chandra homomorphism $\gamma(\underline{k})$ is injective on
$\mathbb{D}_{\mathcal{R}}(L^\prime)^{D_{\pi^\prime}}$.
\end{proof}
%%%%%%%%%%%%%%%%%%%%%%%%%%%%%%%%%%%%%%%%%%%%%%%%%%%%%%%%%%%%%%%%%%%%%%%%%%%%%
\begin{cor}\label{ok}
If $A_0$ is a
$\mathbb{D}_{\mathcal{R}}(L^\prime)^{D_{\pi^\prime}}$-submodule
of $\overline{M}$ then
\[\mathbb{D}_{\mathcal{R}}(L^\prime)^{D_{\pi^\prime}}=
\mathbb{D}_{\mathcal{R}}(L^\prime)^{W_0,D_{\pi^\prime}}.
\]
\end{cor}
%%%%%%%%%%%%%%%%%%%%%%%%%%%%%%%%%%%%%%%%%%%%%%%%%%%%%%%%%%%%%%%%%%%%%%%%%%%%%
We will see in Subsection \ref{C34} that the converse of Corollary
\ref{ok} is also true. Furthermore, in Subsection \ref{polynomialsection}
we use Harish-Chandra series and the theory of Macdonald polynomials
to give conditions on the multiplicity label $\underline{k}$ for
which $A_0$ is a $\mathbb{D}_{\mathcal{R}}(L^\prime)^{D_{\pi^\prime}}$-submodule
of $\overline{M}$.

%%%%%%%%%%%%%%%%%%%%%%%%%%%%%%%%%%%%%%%%
%%                                    %%
%%        Rank reduction              %%
%%                                    %%
%%%%%%%%%%%%%%%%%%%%%%%%%%%%%%%%%%%%%%%%

\subsection{Rank reduction}\label{Rrsub}

In this subsection we compute the asymptotics of the
Macdonald difference operators $D_{\pi^\prime}$ along co-dimension one facets
of the negative Weyl chamber. It allows us to reduce the analysis of the
centralizer algebra $\mathbb{D}_{\mathcal{R}}(L^\prime)^{W_0,D_{\pi^\prime}}$
to the case of rank one root systems $R$. This technique is
reminiscent of the trigonometric differential case, see e.g. \cite{HS} 
and references therein.

For a subset $F\subseteq\Delta$ we define the facet $V_{-,F}$ of
the negative Weyl chamber $V_-$ by
\[V_{-,F}=\{v\in V \, | \, \langle v,\alpha\rangle=0\quad
(\alpha\in F),\qquad \langle v,\beta\rangle<0\quad
(\beta\in\Delta\setminus F) \}.
\]
We write
$R^{\prime\vee}_F=R^{\prime\vee}\cap \mathbb{Z}F$ for the
corresponding standard parabolic root subsystem of $R^{\prime\vee}$,
$R^{\prime\vee,\pm}_F$
for the set of positive respectively negative roots in $R^{\prime\vee}$
with respect to the basis $F$ of $R^{\prime\vee}_F$,
and
$W_{0,F}$ for the standard parabolic subgroup of $W_0$
generated by the simple reflections $s_\alpha$ ($\alpha\in F$).
%Similarly,
We furthermore set
$R_F=R\cap\bigoplus_{\alpha\in F}\mathbb{Q}\alpha$,
$R^{\pm}_F=R^{\pm}\cap\bigoplus_{\alpha\in F}\mathbb{Q}\alpha$ and
\[\rho_{\underline{k}^\prime,F}=
\frac{1}{2}\sum_{\alpha\in R^+\setminus R_F^+}\underline{k}^\prime(\alpha^\vee)
\alpha\in V_{\mathbb{C}}.
\]
Note that $\rho_{\underline{k}^\prime}=\rho_{\underline{k}^\prime}(\emptyset)$.

We define $\mathcal{R}_F\subseteq\mathcal{R}$ to be the
subalgebra generated by
\[\frac{1}{1-rz^\alpha},\qquad r\in\mathbb{C},\,\, \alpha\in R^{\prime\vee}_F.
\]
The embedding $\mathcal{R}\hookrightarrow
\mathbb{C}[[z^{-\Delta}]]$ (see Subsection 2.2) maps 
$\mathcal{R}_F$ into the subalgebra
$\mathbb{C}[[z^{-F}]]:=\mathbb{C}[[z^{-\alpha}\, | \, \alpha\in F]]$ of
$\mathbb{C}[[z^{-\Delta}]]$.
For a lattice $X\subset V$ we write
$\mathbb{D}_{\mathcal{R}_F}(X)=\mathcal{R}_F\# t(X)$
(respectively $\mathbb{D}_{\mathbb{C}[[z^{-F}]]}(X)=
\mathbb{C}[[z^{-F}]]\# t(X)$),
which is a subalgebra of $\mathbb{D}_{\mathcal{R}}(X)$ (respectively
$\mathbb{D}_{\mathbb{C}[[z^{-\Delta}]]}(X)$).
Note that $\mathbb{D}_{\mathcal{R}_F}(X)\subseteq\mathbb{D}_{\mathcal{R}}(X)$
is a $W_{0,F}$-invariant subalgebra if $X$ is a $W_{0,F}$-invariant lattice 
in $V$.
In the special case $F=\emptyset$ 
we have $\mathcal{R}_\emptyset=\mathbb{C}[[z^{-\emptyset}]]=\mathbb{C}$, hence
$\mathbb{D}_{\mathcal{R}_{\emptyset}}(X)=\mathbb{C}[t(X)]$
is the algebra of constant coefficient difference operators with
step-sizes from $X$.

%%%%%%%%%%%%%%%%%%%%%%%%%%%%%%%%%%%%%%%%%%%%%%%%%%%%%%%%%%%%%%%%%%%%%%%%%%%%
\begin{defi}\label{gammaFdef}
Let $F\subseteq \Delta$. The constant term map along the facet $V_{-,F}$
is the algebra homomorphism $\gamma_F(\underline{k}):
\mathbb{D}_{\mathbb{C}[[z^{-\Delta}]]}(L^\prime)\rightarrow
\mathbb{D}_{\mathbb{C}[[z^{-F}]]}(L^\prime)$ defined by
\begin{equation}\label{gammaFk}
\gamma_F(\underline{k})(D)=\sum_{\lambda^\prime\in L^\prime}
\Bigl(\sum_{x\in\mathbb{Z}_{\geq 0}F}C_x(\lambda^\prime)z^{-x}\Bigr)
q^{\langle\rho_{\underline{k}^\prime,F},\lambda^\prime\rangle}t(\lambda^\prime)
\end{equation}
for
\[
D=\sum_{\lambda^\prime\in
L^\prime}\Bigl(\sum_{x\in\mathbb{Z}_{\geq
0}\Delta}C_x(\lambda^\prime)z^{-x}\Bigr)t(\lambda^\prime)
\in\mathbb{D}_{\mathbb{C}[[z^{-\Delta}]]}(L^\prime),
\]
where $\{C_x(\lambda^\prime)\}_{x\in \mathbb{Z}_{\geq 0}\Delta}\subset\mathbb{C}$
is the zero set for all but finitely many $\lambda^\prime\in
L^\prime$.
\end{defi}
%%%%%%%%%%%%%%%%%%%%%%%%%%%%%%%%%%%%%%%%%%%%%%%%%%%%%%%%%%%%%%%%%%%%%%%%%%%
\begin{lem}\label{gammaFklemma}
The map $\gamma_F(\underline{k})$ 
restricts to a $W_{0,F}$-equivariant algebra homomorphism
\[\gamma_F(\underline{k}): 
\mathbb{D}_{\mathcal{R}}(L^\prime)\rightarrow \mathbb{D}_{\mathcal{R}_F}(L^\prime).
\]
\end{lem}
%%%%%%%%%%%%%%%%%%%%%%%%%%%%%%%%%%%%%%%%%%%%%%%%%%%%%%%%%%%%%%%%%%%%%%%%%%%%
\begin{proof}
The map $\gamma_F(\underline{k})$ restricts to an algebra homomorphism
$\gamma_F(\underline{k}): \mathbb{D}_{\mathcal{R}}(L^\prime)\rightarrow 
\mathbb{D}_{\mathcal{R}_F}(L^\prime)$ since
for $r\in\mathbb{C}^\times$,
\begin{equation}\label{restoRF}
\gamma_F(\underline{k})\Bigl(\frac{1}{1-rz^\alpha}\Bigr)=
\begin{cases}
\frac{1}{1-rz^\alpha},\qquad &\alpha\in R^{\prime\vee}_F,\\
1,\qquad &\alpha\in R^{\prime\vee,-}\setminus R^{\prime\vee,-}_F,\\
0,\qquad &\alpha\in R^{\prime\vee,+}\setminus R^{\prime\vee,+}_F.
\end{cases}
\end{equation}
Since $R^{\prime\vee}_F$ and $R^{\prime\vee,\pm}\setminus R^{\prime\vee,\pm}_F$
are $W_{0,F}$-invariant subsets of $R^{\prime\vee}$,
\eqref{restoRF} implies that the restriction of 
$\gamma_F(\underline{k})$ to $\mathcal{R}\subset 
\mathbb{D}_{\mathcal{R}}(L^\prime)$
is $W_{0,F}$-equivariant. Furthermore, 
$\rho_{\underline{k}^\prime,F}\in V_{\mathbb{C}}$
is $W_{0,F}$-stable since $R^+\setminus R_F^+$
is $W_{0,F}$-invariant. It follows that 
$\gamma_F(\underline{k}): 
\mathbb{D}_{\mathcal{R}}(L^\prime)\rightarrow \mathbb{D}_{\mathcal{R}_F}(L^\prime)$
is $W_{0,F}$-equivariant.
\end{proof}
%%%%%%%%%%%%%%%%%%%%%%%%%%%%%%%%%%%%%%%%%%%%%%%%%%%%%%%%%%%%%%%%%%%%%%%%%%%%

Observe that $\gamma(\underline{k})
=\delta_F(\underline{k})\circ\gamma_F(\underline{k})$ with
$\delta_F(\underline{k}): \mathbb{D}_{\mathbb{C}[[z^{-F}]]}(L^\prime)\rightarrow
A^\prime$ the algebra homomorphism defined by
\begin{equation}\label{deltaFk}
\delta_F(\underline{k})(D)=\sum_{\lambda^\prime\in L^\prime}
C_0(\lambda^\prime)
q^{\langle\rho_{\underline{k}^\prime}-\rho_{\underline{k}^\prime,F},\lambda^\prime\rangle}
\xi^{\lambda^\prime}
\end{equation}
for
\[
D=\sum_{\lambda^\prime\in
L^\prime}\Bigl(\sum_{x\in\mathbb{Z}_{\geq
0}F}C_x(\lambda^\prime)z^{-x}\Bigr)t(\lambda^\prime)
\in\mathbb{D}_{\mathbb{C}[[z^{-F}]]}(L^\prime),
\]
where $\{C_x(\lambda^\prime)\}_{x\in \mathbb{Z}_{\geq 0}F}\subset\mathbb{C}$
is the zero set for all but finitely many $\lambda^\prime\in
L^\prime$.

%%%%%%%%%%%%%%%%%%%%%%%%%%%%%%%%%%%%%%%%%%%%%%%%%%%%%%%%%%%%%%%%%%%%%%%%%%
Our aim is to compute the rank one reductions 
$\gamma_{\{\alpha\}}(\underline{k})(D_{\pi^\prime})$
($\alpha\in\Delta$) of the Macdonald difference operator $D_{\pi^\prime}$.
For $F\subseteq \Delta$, set $S_F=\{a\in S \,\, | \,\, Da\in\mathbb{Z}F \}$.
We need several preparatory lemmas.
%%%%%%%%%%%%%%%%%%%%%%%%%%%%%%%%%%%%%%%%%%%%%%%%%%%%%%%%%%%%%%%%%%%%%%%%%%%%%
\begin{lem}\label{cross}
Let $\pi^\prime\in L^\prime$ be a nonzero 
antidominant minuscule or quasi-minuscule co-weight.
Fix $i\in\{1,\ldots,n\}$ and $w\in W_0^{\pi^\prime}$.

{\bf a)} If $\pi^\prime=w_0\pi_j^\prime$ \textup{(}$j\in J_0$\textup{)}, then
\begin{equation*}
S_{\{a_i\}}\cap wS_1(t(-\pi^\prime))=
\begin{cases}
\{a_i\} \,\, &\hbox{ if }\,\,\langle a_i,w\pi^\prime\rangle<0,\\
\{-a_i\} \,\, &\hbox{ if }\,\, \langle a_i,w\pi^\prime\rangle>0,\\
\emptyset\,\, &\hbox{ if }\,\,\langle a_i,w\pi^\prime\rangle=0.
\end{cases}
\end{equation*}

{\bf b)} If $\pi^\prime=-\varphi^\vee$, then
\begin{equation*}
S_{\{a_i\}}\cap wS_1(t(-\pi^\prime))=
\begin{cases}
\{a_i, a_i+2c/\|\varphi\|^2 \}\,\, &\hbox{ if }\,\,
w\pi^\prime=-a_i,\\
\{-a_i,-a_i+2c/\|\varphi\|^2 \}\,\, &\hbox{ if }\,\,
w\pi^\prime=a_i,\\
\{a_i\} \,\, &\hbox{ if }\,\, \langle a_i,w\pi^\prime\rangle<0\,
\hbox{ and }\, w\pi^\prime\not=-a_i,\\
\{-a_i\} \,\, &\hbox{ if }\,\, 
\langle a_i,w\pi^\prime\rangle>0\,\hbox{ and }\,
w\pi^\prime\not=a_i,\\
\emptyset\,\, &\hbox{ if }\,\, \langle a_i,w\pi^\prime\rangle=0.
\end{cases}
\end{equation*}
\end{lem}
%%%%%%%%%%%%%%%%%%%%%%%%%%%%%%%%%%%%%%%%%%%%%%%%%%%%%%%%%%%%%%%%%%%%%%%%
\begin{proof}
This follows directly from Lemma \ref{crosspre}.
\end{proof}
%%%%%%%%%%%%%%%%%%%%%%%%%%%%%%%%%%%%%%%%%%%%%%%%%%%%%%%%%%%%%%%%%%%%%%%%%%%
Recall from Proposition \ref{formexplicit}
the rational function $f_{\pi^\prime}(z)\in\mathcal{R}$
occurring as coefficient in the Macdonald difference 
operator $D_{\pi^\prime}$.
%%%%%%%%%%%%%%%%%%%%%%%%%%%%%%%%%%%%%%%%%%%%%%%%%%%%%%%%%%%%%%%%%%%%%%%%%%%%%
\begin{lem}\label{cross2}
Let $i\in\{1,\ldots,n\}$, 
$w\in W_0^{\pi^\prime}$  and let $\pi^\prime\in L^\prime$ be a
nonzero anti-dominant minuscule or quasi-minuscule co-weight. We
have
\[\gamma_{\{a_i\}}(\underline{k})\bigl((wf_{\pi^\prime})(z)\bigr)=
q^{-\langle\rho_{\underline{k}^\prime,\{a_i\}},
w\pi^\prime\rangle}\prod_{a\in S_{\{a_i\}}\cap
wS_1(t(-\pi^\prime))}c_a(z).
\]
\end{lem}
%%%%%%%%%%%%%%%%%%%%%%%%%%%%%%%%%%%%%%%%%%%%%%%%%%%%%%%%%%%%%%%%%%%%%%%%%%
\begin{proof}
We define a function $\epsilon: S\rightarrow \{\pm 1\}$ by
$\epsilon(a)=1$ if $Da\in\mathbb{Z}_{\geq 0}\Delta$ and $\epsilon(a)=-1$
if $a\in \mathbb{Z}_{\leq 0}\Delta$. For $a\in S_1$ we then have
\begin{equation}\label{casymp}
\gamma_{\{a_i\}}(\underline{k})(c_a(z))=
\begin{cases}
c_a(z),\qquad &\hbox{ if } a\in S_{\{a_i\}},\\
\tau_a^{\epsilon(a)},\qquad &\hbox{ if } a\in S\setminus S_{\{a_i\}},
\end{cases}
\end{equation}
hence
\begin{equation}\label{fasymp}
\gamma_{\{a_i\}}(\underline{k})\bigl((wf_{\pi^\prime})(z)\bigr)=
r_i(w)\prod_{a\in S_{\{a_i\}}\cap
wS_1(t(-\pi^\prime))}c_a(z)
\end{equation}
for certain $r_i(w)\in\mathbb{C}^\times$.
It remains to show that
$r_i(w)=q^{-\langle\rho_{\underline{k}^\prime,\{a_i\}},w\pi^\prime\rangle}$.
By Corollary \ref{Macdonaldreturn}, \eqref{casymp}, \eqref{fasymp} 
and the fact
that $\gamma(\underline{k})=
\delta_{\{a_i\}}(\underline{k})\circ\gamma_{\{a_i\}}(\underline{k})$,
we have
\[
r_i(w)=q^{-\langle\rho_{\underline{k}^\prime},w\pi^\prime\rangle}\prod_{a\in
S_{\{a_i\}}\cap wS_1(t(-\pi^\prime))}\tau_a^{-\epsilon(a)}.
\]
Using that
$\rho_{\underline{k}^\prime}=
\rho_{\underline{k}^\prime,\{a_i\}}+
\frac{1}{2}\underline{k}^\prime(\alpha_i^\vee)\alpha_i$,
it suffices to show that
\begin{equation}\label{identitybasic}
\prod_{a\in S_{\{a_i\}}\cap
wS_1(t(-\pi^\prime))}\tau_a^{\epsilon(a)}=
q^{-\frac{1}{2}\underline{k}^\prime(\alpha_i^\vee)
\langle\alpha_i,w\pi^\prime\rangle}.
\end{equation}
This follows from straightforward 
computations using Lemma \ref{cross}. As an
example we discuss the proof for case {\bf c} in detail. 
For case {\bf c}, $w\pi^\prime=-w\varphi^\vee$ is
quasi-minuscule. Suppose first that $w\pi^\prime\not=\pm
a_i$. Since $a_i=\alpha_i^\vee$, we have 
$(w\pi^\prime)^\vee\not=\pm \alpha_i$, hence
$\langle \alpha_i,w\pi^\prime\rangle\in\{-1,0,1\}$. 
If $\langle \alpha_i,w\pi^\prime\rangle=0$,
then both sides of \eqref{identitybasic} 
are equal to one. Suppose now that 
$\langle\alpha_i,w\pi^\prime\rangle\in\{\pm
1\}$. Since $\pi^\prime$ is anti-dominant,
the left hand side of \eqref{identitybasic} 
equals $\tau_{a_i}^{\epsilon(w^{-1}(a_i))}$
by Lemma \ref{cross}. On the other hand, 
the assumptions imply that
$\alpha_i\in R=R_C$ is a short root, hence $\|\alpha_i\|^2=2$
and $a_i=\alpha_i^\vee\in\mathcal{O}_5$. Consequently, 
$\underline{k}^\prime(a_i)=\underline{k}(a_i)=\kappa_5$
and $\tau_{a_i}=q^{\frac{1}{2}\underline{k}^\prime(\alpha_i^\vee)}$.
The right hand side of \eqref{identitybasic} thus also equals 
$\tau_{a_i}^{\epsilon(w^{-1}(a_i))}$.
If $w\pi^\prime\in\{\pm a_i\}$, then 
by Lemma \ref{cross} the left hand side of
\eqref{identitybasic} equals
\begin{equation}\label{expression}
\bigl(\tau_{a_i}\tau_{a_i+c/2}\bigr)^{\epsilon(w^{-1}(a_i))}=
q^{\epsilon(w^{-1}(a_i))\kappa_1^\prime},
\end{equation}
where we have used the fact that 
$a_i\in \mathcal{O}_1$ and $a_i+c/2\in\mathcal{O}_3$.
For the right hand side of \eqref{identitybasic} we use 
that $\alpha_i^\vee=a_i\in\mathcal{O}_1$
and $\langle \alpha_i,w\pi^\prime\rangle=-2\epsilon(w^{-1}(a_i))$ to
arrive at the same expression \eqref{expression}.
\end{proof}
%%%%%%%%%%%%%%%%%%%%%%%%%%%%%%%%%%%%%%%%%%%%%%%%%%%%%%%%%%%%%%%%%%%%%%%%%%%%%
For $i\in \{1,\ldots,n\}$ we write
\[L^\prime_i=L^\prime+\mathbb{Z}\frac{\alpha_i^\vee}{2}
\]
unless we are in case {\bf c} and $a_i\in W_0\varphi^\vee$, in
which case we set $L^\prime_i=L^\prime$.
Set
\[X_i=\{\mu^\prime\in L_i^\prime \,\, | \,\,
\langle\mu^\prime,\alpha_i\rangle=0\}.
\]
The lattice $L_i^\prime$ has the following elementary properties.
%%%%%%%%%%%%%%%%%%%%%%%%%%%%%%%%%%%%%%%%%%%%%%%%%%%%%%%%%%%%%%%%%%%%%%%%%%
\begin{lem}\label{latticeextension}
The lattice $L_i^\prime$ is $s_i$-invariant. 
It decomposes as the orthogonal direct sum
\begin{equation}\label{orthodecomp1}
L_i^\prime=\mathbb{Z}\frac{\alpha_i^\vee}{2}\oplus X_i
\end{equation}
unless we are in case {\bf c} and $a_i\in W_0\varphi^\vee$, in
which case it decomposes as the orthogonal direct sum
\begin{equation}\label{orthodecomp2}
L_i^\prime=\mathbb{Z}a_i\oplus X_i.
\end{equation}
\end{lem}
%%%%%%%%%%%%%%%%%%%%%%%%%%%%%%%%%%%%%%%%%%%%%%%%%%%%%%%%%%%%%%%%%%%%%%%%%%
\begin{proof}
If we are in case {\bf c} and 
$a_i\in W_0\varphi^\vee$, then $a_i=\alpha_i^\vee\in Q(R^\vee)=L^\prime$
with $\alpha_i$ the unique long simple root from the basis 
$\{\alpha_1,\ldots,\alpha_n\}$
of $R=R_C$, hence $\langle L^\prime,\alpha_i\rangle=2\mathbb{Z}$
(cf. \cite[(2.1.6)]{M}). Since $a_i=\alpha_i^\vee=\alpha_i/2$, 
we conclude that
$\langle L^\prime,a_i\rangle\in\mathbb{Z}$. Combined with the
observation that $a_i$ has squared length one, we obtain the
orthogonal decomposition \eqref{orthodecomp2} for $L_i^\prime=L^\prime$.

Suppose now that we are in case {\bf a} or case {\bf b}, or that we are in
case {\bf c} and $a_i\not\in W_0\varphi^\vee$.
Since $L^\prime\subseteq P(R^\vee)$ and 
$\langle \alpha_i^\vee/2,\alpha_i\rangle=1$
we have $\langle L_i^\prime, \alpha_i\rangle=\mathbb{Z}$. The
orthogonal decomposition \eqref{orthodecomp1} follows now immediately.
\end{proof}
%%%%%%%%%%%%%%%%%%%%%%%%%%%%%%%%%%%%%%%%%%%%%%%%%%%%%%%%%%%%%%%%%%%%%%%%%%%%
We now consider the algebra $\mathbb{D}_{\mathcal{R}_{\{a_i\}}}(L_i^\prime)$ of
difference operators with coefficients from $\mathcal{R}_{\{a_i\}}$ and 
with step-sizes
from the enlarged lattice $L_i^\prime$. Since $L_i^\prime$ is
$s_i$-invariant, $\mathbb{D}_{\mathcal{R}_{\{a_i\}}}(L_i^\prime)$
is a $W_{0,\{a_i\}}$-module algebra. Lemma \ref{latticeextension} directly
implies
the following result.
%%%%%%%%%%%%%%%%%%%%%%%%%%%%%%%%%%%%%%%%%%%%%%%%%%%%%%%%%%%%%%%%%%%%%%%%%%%
\begin{cor}\label{algebrared}
The center $Z\bigl(\mathbb{D}_{\mathcal{R}_{\{a_i\}}}(L_i^\prime)\bigr)$
of $\mathbb{D}_{\mathcal{R}_{\{a_i\}}}(L_i^\prime)$
is the algebra of constant coefficient difference operators
with step-sizes from $X_i$,
\[Z\bigl(\mathbb{D}_{\mathcal{R}_{\{a_i\}}}(L_i^\prime)\bigr)=
\mathbb{C}[t(X_i)].
\]
Furthermore,
\begin{equation}\label{rankonesplitting}
\mathbb{D}_{\mathcal{R}_{\{a_i\}}}(L_i^\prime)\simeq
\mathbb{D}_{\mathcal{R}_{\{a_i\}}}
\Bigl(\mathbb{Z}\frac{\alpha_i^\vee}{2}\Bigr)\otimes
\mathbb{C}[t(X_i)]
\end{equation}
as algebras unless we are in case {\bf c} and $a_i\in W_0\varphi^\vee$, in
which case
\begin{equation}\label{rankonesplitting2}
\mathbb{D}_{\mathcal{R}_{\{a_i\}}}(L_i^\prime)\simeq
\mathbb{D}_{\mathcal{R}_{\{a_i\}}}\bigl(\mathbb{Z}a_i\bigr)\otimes
\mathbb{C}[t(X_i)]
\end{equation}
as algebras, with the isomorphisms 
in \eqref{rankonesplitting} and \eqref{rankonesplitting2}
realized by the multiplication map.
\end{cor}
%%%%%%%%%%%%%%%%%%%%%%%%%%%%%%%%%%%%%%%%%%%%%%%%%%%%%%%%%%%%%%%%%%%%%%%%%%%

We now decompose the operator
$\gamma_{\{a_i\}}(\underline{k})(D_{\pi^\prime})
\in\mathbb{D}_{\mathcal{R}_{\{a_i\}}}(L_i^\prime)$
according to the decomposition \eqref{rankonesplitting} 
respectively \eqref{rankonesplitting2}.
The result can be described in terms of an explicit rank one
difference operator, which we now define first (its definition
depends on the case involved).
%%%%%%%%%%%%%%%%%%%%%%%%%%%%%%%%%%%%%%%%%%%%%%%%%%%%%%%%%%%%%%%%%%%%%%%%%%%%%%
\begin{defi}\label{rankoneL}
For $i\in\{1,\ldots,n\}$ we define 
$\mathcal{L}_i\in
\mathbb{D}_{\mathcal{R}_{\{a_i\}}}
\bigl(\mathbb{Z}\alpha_i^\vee/2\bigr)^{W_{0,\{a_i\}}}$
by
\begin{equation}
\label{mathcalL}
\mathcal{L}_i=c_{a_i}(z)t(-\alpha_i^\vee/2)+c_{-a_i}(z)t(\alpha_i^\vee/2)
\end{equation}
unless we are in case {\bf c} and $a_i\in W_0\varphi^\vee$, in which 
case $\mathcal{L}_i\in
\mathbb{D}_{\mathcal{R}_{\{a_i\}}}\bigl(\mathbb{Z}a_i\bigr)^{W_{0,\{a_i\}}}$ 
is defined by
\begin{equation}\label{mathcalLAW}
\mathcal{L}_i=c_{a_i}(z)c_{a_i+c/2}(z)\bigl(t(-a_i)-1\bigr)+
c_{-a_i}(z)c_{-a_i+c/2}(z)\bigl(t(a_i)-1\bigr).
\end{equation}
\end{defi}
%%%%%%%%%%%%%%%%%%%%%%%%%%%%%%%%%%%%%%%%%%%%%%%%%%%%%%%%%%%%%%%%%%%%%%%%%%%%
Polynomial eigenfunctions of the rank one difference operator 
\eqref{mathcalL} (respectively
\eqref{mathcalLAW}) are the continuous $q$-ultraspherical polynomials
(respectively the Askey-Wilson polynomials), see \cite[Chapter 6]{M}
and references therein.

We define constant coefficient difference operators 
$y_i, z_i\in\mathbb{C}[X_i]$ by
\begin{equation}\label{yz}
y_i=\sum_{\stackrel{w\in W_0^{\pi^\prime}:}{\langle
\alpha_i,w\pi^\prime\rangle=-1}}t\Bigl(\frac{\alpha_i^\vee}{2}+
w\pi^\prime\Bigr),
\qquad z_i=\sum_{\stackrel{w\in
W_0^{\pi^\prime}:}{\langle\alpha_i,w\pi^\prime\rangle=0}}t(w\pi^\prime).
\end{equation}
If we are in case {\bf c} and if $a_i\in W_0\varphi^\vee$, then $y_i$
should be read as zero. This convention is in accordance with the
following slight variation of \cite[(2.1.6)]{M}.
%%%%%%%%%%%%%%%%%%%%%%%%%%%%%%%%%%%%%%%%%%%%%%%%%%%%%%%%%%%%%%%%%%
\begin{lem}\label{ynonzero}
The set
\begin{equation}\label{set}
\{w\in W_0^{\pi^\prime} \,\,\, | \,\,\,
\langle\alpha_i,w\pi^\prime\rangle=-1\}
\end{equation}
is nonempty
unless $R$ is of type $C$, 
$\pi^\prime=-\varphi^\vee$ and $a_i\in W_0\varphi^\vee$.
\end{lem}
%%%%%%%%%%%%%%%%%%%%%%%%%%%%%%%%%%%%%%%%%%%%%%%%%%%%%%%%%%%%%%%%%%%%
\begin{proof}
Suppose that $\pi^\prime$ is minuscule, 
or that $\pi^\prime=-\varphi^\vee$ and $a_i\not\in
W_0\varphi^\vee$. Then $\langle
\alpha_i,w\pi^\prime\rangle\in\{-1,0,1\}$ for all $w\in W_{0}^{\pi^\prime}$.
Since $W_0$ acts irreducibly on $V$, there exists a
$w\in W_0^{\pi^\prime}$ such that $\langle \alpha_i,w\pi^\prime\rangle\not=0$,
hence \eqref{set} is nonempty.

Suppose now that $\pi^\prime=-\varphi^\vee$ and $a_i\in W_0\varphi^\vee$.
Then $\alpha_i$ is a long root in $R$. It follows from the root system 
classification that
$\langle \alpha_i,\beta^\vee\rangle=-1$ for some long root $\beta\in
R$, unless $R$ is of type $C$. Since long roots in $R$ are
$W_0$-conjugate to $\varphi$,
we conclude that \eqref{set} is nonempty unless $R$ is of type $C$.
\end{proof}
%%%%%%%%%%%%%%%%%%%%%%%%%%%%%%%%%%%%%%%%%%%%%%%%%%%%%%%%%%%%%%%%%

%%%%%%%%%%%%%%%%%%%%%%%%%%%%%%%%%%%%%%%%%%%%%%%%%%%%%%%%%%%%%%%%%%%%%%%%%%%
\begin{prop}\label{casei}
Fix $i\in\{1,\ldots,n\}$.\\
{\bf i)} Suppose that 
$\pi^\prime=w_0\pi_j^\prime$ \textup{(}$j\in J_0$\textup{)}, or that
$\pi^\prime=-\varphi^\vee$ and $a_i\not\in W_0\varphi^\vee$. Then
\[\gamma_{\{a_i\}}(\underline{k})(D_{\pi^\prime})=
y_i\mathcal{L}_i+z_i.
\]
{\bf ii)} For cases {\bf a} and 
{\bf b} with $\pi^\prime=-\varphi^\vee$ and $a_i\in
W_0\varphi^\vee$ we have
\[\gamma_{\{a_i\}}(\underline{k})\bigl(D_{\pi^\prime}\bigr)=
\mathcal{L}_i^2+y_i\mathcal{L}_i+z_i-2.
\]
{\bf iii)} For case {\bf c} with 
$\pi^\prime=-\varphi^\vee$ and $a_i\in W_0\varphi^\vee$ we
have
\[\gamma_{\{a_i\}}(\underline{k})\bigl(D_{\pi^\prime}\bigr)=\mathcal{L}_i+z_i+
q^{\kappa_1^\prime}+q^{-\kappa_1^\prime}.
\]
\end{prop}
%%%%%%%%%%%%%%%%%%%%%%%%%%%%%%%%%%%%%%%%%%%%%%%%%%%%%%%%%%%%%%%%%%%%%%%%%%%%
\begin{proof}
Note that
Proposition \ref{formexplicit} and Lemma \ref{cross2} gives the
initial expression
\begin{equation}\label{initialexp}
\begin{split}
\gamma_{\{a_i\}}(\underline{k})(D_{\pi^\prime})&=
m_{\pi^\prime}(-\rho_{\underline{k}^\prime})+\sum_{w\in W_0^{\pi^\prime}}g_{w,i}(z)
\Bigl(t(w\pi^\prime)-q^{-\langle\rho_{\underline{k}^\prime,\{a_i\}},w\pi^\prime\rangle}
\Bigr),\\
g_{w,i}(z)&=\prod_{a\in S_{\{a_i\}}\cap wS_1(t(-\pi^\prime))}c_a(z).
\end{split}
\end{equation}
The further computations depend on the three different cases.\\
{\bf i)} Under the present assumptions, 
$a_i=\alpha_i^\vee\in\mathcal{O}_5$ for case {\bf c},
$\langle a_i,w\pi^\prime\rangle<0$ for $w\in
W_0^{\pi^\prime}$ implies $\langle
\alpha_i,w\pi^\prime\rangle=-1$, and
\[q^{-\langle\rho_{\underline{k}^\prime,\{a_i\}},\mu^\prime\rangle}=
q^{-\langle\rho_{\underline{k}^\prime},\mu^\prime\rangle}
\tau_{a_i}^{\langle\alpha_i,\mu^\prime\rangle},
\qquad \mu^\prime\in L_i^\prime
\]
(compare with the proof of Lemma \ref{cross2}).
Combined with Lemma \ref{cross} and the $s_i$-invariance of
$\rho_{\underline{k}^\prime,\{a_i\}}$,
the expression \eqref{initialexp} becomes
\begin{equation*}
\begin{split}
\gamma_{\{a_i\}}(\underline{k})(D_{\pi^\prime})&=
m_{\pi^\prime}(-\rho_{\underline{k}^\prime})+\sum_{\stackrel{w\in
W_0^{\pi^\prime}:}{\langle\alpha_i,w\pi^\prime\rangle=0}}
\Bigl(t(w\pi^\prime)-
q^{-\langle\rho_{\underline{k}^\prime},w\pi^\prime\rangle}\Bigr)\\
&+\sum_{\stackrel{w\in
W_0^{\pi^\prime}:}{\langle\alpha_i,w\pi^\prime\rangle=-1}}
\Bigl(c_{a_i}(z)t(w\pi^\prime)+c_{-a_i}(z)t(s_iw\pi^\prime)\Bigr)\\
&-\sum_{\stackrel{w\in
W_0^{\pi^\prime}:}{\langle\alpha_i,w\pi^\prime\rangle=-1}}
\Bigl(c_{a_i}(z)+c_{-a_i}(z)\Bigr)\tau_{a_i}^{-1}
q^{-\langle\rho_{\underline{k}^\prime},w\pi^\prime\rangle}.
\end{split}
\end{equation*}
The second line is $y_i\mathcal{L}_i$. To show that the remaining
terms sum up to $z_i$, it suffices to note that
\[\sum_{\stackrel{w\in
W_0^{\pi^\prime}:}{\langle\alpha_i,w\pi^\prime\rangle=-1}}
\Bigl(c_{a_i}(z)+c_{-a_i}(z)\Bigr)\tau_{a_i}^{-1}
q^{-\langle\rho_{\underline{k}^\prime},w\pi^\prime\rangle}=
\sum_{\stackrel{w\in W_0^{\pi^\prime}:}{\langle
\alpha_i,w\pi^\prime\rangle\not=0}}
q^{-\langle\rho_{\underline{k}^\prime},w\pi^\prime\rangle}
\]
by \eqref{ceq}.\\
{\bf ii)} \& {\bf iii)} Under the present 
assumptions, $\langle w\pi^\prime,a_i\rangle<0$
and $w\pi^\prime\not\in\{\pm a_i\}$ imply
$\langle\alpha_i,w\pi^\prime\rangle=-1$.
As in the proof of part {\bf i)}, Lemma \ref{cross}
and \eqref{initialexp} give
\begin{equation}\label{initialexp2}
\begin{split}
\gamma_{\{a_i\}}(\underline{k})\bigl(D_{\pi^\prime}\bigr)&=
c_{a_i}(z)c_{a_i+\frac{2c}{\|\varphi\|^2}}(z)\bigl(t(-a_i)-1\bigr)+
c_{-a_i}(z)c_{-a_i+\frac{2c}{\|\varphi\|^2}}(z)\bigl(t(a_i)-1\bigr)\\
&+y_i\mathcal{L}_i+z_i+
q^{\langle\rho_{\underline{k}^\prime},a_i\rangle}+
q^{-\langle\rho_{\underline{k}^\prime},a_i\rangle}.
\end{split}
\end{equation}
For case {\bf c}, $\|\varphi\|^2=4$ and 
$q^{\langle\rho_{\underline{k}^\prime},a_i\rangle}=q^{\kappa_1^\prime}$,
while Lemma \ref{ynonzero} implies that $y_i=0$. This proves {\bf iii)}.

For case {\bf a} and case {\bf b}, the root $a_i$ has squared length two,
hence $a_i=\alpha_i=\alpha_i^\vee$. Since the labeling of $S$ only 
depends on the
gradient root system $R^{\prime\vee}$ of $S$, it follows that
\[t\bigl(\pm \alpha_i^\vee/2\bigr)(c_{a_i}(z))=c_{a_i\mp
c}(z),\qquad t\bigl(\pm\alpha_i^\vee/2\bigr)(c_{-a_i}(z))=c_{-a_i\pm
c}(z).
\]
Combined with \eqref{ceq}, we obtain
\begin{equation*}
\mathcal{L}_i^2=c_{a_i}(z)c_{a_i+c}(z)\bigl(t(-a_i)-1\bigr)+
c_{-a_i}(z)c_{-a_i+c}(z)\bigl(t(a_i)-1\bigr)
+\bigl(\tau_{a_i}+\tau_{a_i}^{-1}\bigr)^2.
\end{equation*}
Since $\|\varphi\|^2=2$ and
$q^{\langle\rho_{\underline{k}^\prime},a_i\rangle}=\tau_{a_i}^2$,
part {\bf ii)} of the proposition now follows from
\eqref{initialexp2}.
\end{proof}
%%%%%%%%%%%%%%%%%%%%%%%%%%%%%%%%%%%%%%%%%%%%%%%%%%%%%%%%%%%%%%%%%%%%%%%%%%%%%

%%%%%%%%%%%%%%%%%%%%%%%%%%%%%%%%%%%%%%%%%%%%%%%%%%%%%%
%%                                                  %%
%%  The centralizer theorem                         %%
%%                                                  %%
%%%%%%%%%%%%%%%%%%%%%%%%%%%%%%%%%%%%%%%%%%%%%%%%%%%%%%

\subsection{Centralizers}\label{C34}

Let $\pi^\prime\in L^\prime$ be a nonzero anti-dominant minuscule
or quasi-minuscule co-weight. For $i\in\{1,\ldots,n\}$ we denote
$\mathbb{D}_i^{\pi^\prime}(\underline{k})$ for the centralizer of
$\gamma_{\{a_i\}}(\underline{k})(D_{\pi^\prime})$
in $\mathbb{D}_{\mathcal{R}_{\{a_i\}}}(L_i^\prime)^{W_{0,\{a_i\}}}$.
Write $\mathbb{C}[\mathcal{L}_i]\subseteq
\mathbb{D}_{\mathcal{R}_{\{a_i\}}}(L_i^\prime)^{W_{0,\{a_i\}}}$
for the unital subalgebra generated by the difference 
operator $\mathcal{L}_i$.

%%%%%%%%%%%%%%%%%%%%%%%%%%%%%%%%%%%%%%%%%%%%%%%%%%%%%%%%%%%%%%%%%%%%%%%%%%%%
\begin{prop}\label{commutantprop}
Let $i\in\{1,\ldots,n\}$ and let $\pi^\prime\in L^\prime$ be a
nonzero anti-dominant minuscule or quasi-minuscule co-weight.
The restriction of the multiplication map 
$D\otimes D^\prime\mapsto DD^\prime$
\textup{(}$D,D^\prime\in \mathbb{D}_{\mathcal{R}_{\{a_i\}}}(L_i^\prime)$\textup{)}
to $\mathbb{C}[\mathcal{L}_i]\otimes\mathbb{C}[X_i]$ defines an
algebra isomorphism
\[\mu_i: \mathbb{C}[\mathcal{L}_i]
\otimes_{\mathbb{C}}\mathbb{C}[X_i]\overset{\sim}{\longrightarrow}
\mathbb{D}_i^{\pi^\prime}(\underline{k}).
\]
\end{prop}
%%%%%%%%%%%%%%%%%%%%%%%%%%%%%%%%%%%%%%%%%%%%%%%%%%%%%%%%%%%%%%%%%%%%%%%%%%%%
\begin{proof}
By Corollary \ref{algebrared}, 
the multiplication map restricts to an injective algebra homomorphism
\[\mu_i: \mathbb{C}[\mathcal{L}_i]\otimes_{\mathbb{C}}\mathbb{C}[X_i]
\hookrightarrow
\mathbb{D}_{\mathcal{R}_{\{a_i\}}}(L_i^\prime).
\]
A constant coefficient difference operator with step-sizes 
from $X_i$ commutes with
the difference operator 
$\gamma_{\{a_i\}}(\underline{k})(D_{\pi^\prime})$ by Corollary
\ref{algebrared}, and it is $s_i$-invariant since $s_i$ 
fixes $X_i$ point-wise.
Furthermore, by Corollary \ref{algebrared} and Proposition
\ref{casei} we have $\mathbb{C}[\mathcal{L}_i]\subseteq
\mathbb{D}_i^{\pi^\prime}(\underline{k})$.
We conclude that the image of $\mu_i$ is 
contained in $\mathbb{D}_i^{\pi^\prime}(\underline{k})$.

Before proving that $\mu_i$ is a linear isomorphism, 
we first introduce some convenient notations and terminology.
We set $\upsilon_i=\alpha_i^\vee/2$ unless we are in case {\bf
c} with $a_i\in W_0\varphi^\vee$, in which case we set
$\upsilon_i=a_i$. For a nonzero difference operator 
$D\in \mathbb{D}_i^{\pi^\prime}(\underline{k})$,
consider its unique expansion
\[D=\sum_{\stackrel{m\in\mathbb{Z}}{x\in
X_i}}g_m^x(z)t\bigl(x+m\upsilon_i\bigr)
\]
with $g_m^x(z)\in\mathcal{R}_{\{a_i\}}$ nonzero for at most
finitely many pairs $(m,x)\in\mathbb{Z}\times X_i$.
Since $s_i$ fixes $X_i$
point-wise, the $s_i$-invariance of $D$ implies
$s_i(g_m^x(z))=g_{-m}^x(z)$. We call
\[D^{(m)}:=\sum_{x\in
X_i}g_m^x(z)t(x+m\upsilon_i)
\]
the $m$th order term of the difference operator $D$.
We write $M(D)$ for the largest integer $m$ for which
$D^{(m)}\not=0$. Since $s_i(D^{(m)})=D^{(-m)}$ we have
$M(D)\geq 0$.

Fix now a nonzero difference operator
$D\in\mathbb{D}_i^{\pi^\prime}(\underline{k})$ and set $M=M(D)$.
Consider the decomposition $D=\sum_{x\in X_i}D_xt(x)$
with $D_x\in \mathbb{D}_{\mathcal{R}_{\{a_i\}}}(\mathbb{Z}\upsilon_i)$
the $s_i$-invariant difference operator
\[D_x:=\sum_{m\in\mathbb{Z}}g_m^x(z)t\bigl(m\upsilon_i\bigr),\qquad
x\in X_i.
\]
We have to show that $D_x\in \mathbb{C}[\mathcal{L}_i]$ for all $x\in X_i$.

By Corollary \ref{algebrared}, Lemma \ref{ynonzero} and
Proposition \ref{casei}, the fact that $D$ centralizes 
$\gamma_{\{a_i\}}(\underline{k})(D_{\pi^\prime})$
implies
\begin{equation}\label{identitycomm1}
\lbrack D,\mathcal{L}_i\rbrack=0
\end{equation}
unless we are in case {\bf a} or 
case {\bf b} with $\pi^\prime=-\varphi^\vee$ and
$a_i\in W_0\varphi^\vee$, in which case it implies
\begin{equation}\label{identitycomm}
\lbrack D,\mathcal{L}_i^2\rbrack
+\lbrack D,\mathcal{L}_i\rbrack y_i=0.
\end{equation}
We set $d_i(z)=c_{-a_i}(z)$ unless we are in case {\bf c} with $a_i\in
W_0\varphi^\vee$, in which case we set
$d_i(z)=c_{-a_i}(z)c_{-a_i+c/2}(z)$. With this notation, the highest
order term of $\mathcal{L}_i$ is $d_i(z)t(\upsilon_i)$.
Considering the highest order term of the identity
\eqref{identitycomm1} and \eqref{identitycomm} respectively, we obtain
\[\lbrack D^{(M)}, d_i(z)t(\upsilon_i)\rbrack=0
\]
unless we are in case {\bf a} or case {\bf b} with $\pi^\prime=-\varphi^\vee$
and $a_i\in W_0\varphi^\vee$, in which case we obtain
\begin{equation}\label{relevantcommutant}
\lbrack D^{(M)}, \bigl(d_i(z)t(\upsilon_i)\bigr)^2\rbrack=0.
\end{equation}
In particular, the highest order term $D^{(M)}$ satisfies
\eqref{relevantcommutant} for all the cases under consideration.

It follows from \eqref{relevantcommutant} that the coefficients
$g_M^x(z)\in\mathcal{R}_{\{a_i\}}$
($x\in X_i$) are solutions of the difference equation
\begin{equation}\label{firstorderdiff}
t(2\upsilon_i)\bigl(f(z)\bigr)=
\left(\frac{t(M\upsilon_i)(e_i(z))}{e_i(z)}\right)f(z),
\end{equation}
where $e_i(z)=
d_i(z)\bigl(t(\upsilon_i)(d_i(z))\bigr)\in\mathcal{R}_{\{a_i\}}$.
The space of functions $f(z)\in\mathcal{R}_{\{a_i\}}$
satisfying \eqref{firstorderdiff} is an one-dimensional 
complex vector space
spanned by
\[
f_M(z)=\prod_{j=0}^{M-1}
\Bigl(t(j\upsilon_i)\bigl(d_i(z)\bigr)\Bigr)\in\mathcal{R}_{\{a_i\}},
\]
hence $g_M^x(z)=K_M^xf_M(z)$ for some $K_M^x\in\mathbb{C}$ ($x\in X_i$).
In particular, $D_x\in\mathbb{C}$
for all $x\in X_i$ if $M=0$, which proves 
that $D_x\in\mathbb{C}[\mathcal{L}_i]$ ($x\in X_i$)
for $M=0$.

Let $M>0$ and suppose that 
$D_x^\prime\in\mathbb{C}[\mathcal{L}_i]$ for all $x\in
X_i^\prime$ if $0\not=D^\prime\in\mathbb{D}_i^{\pi^\prime}(\underline{k})$ 
and $M(D^\prime)<M$.
Let $0\not=D\in\mathbb{D}_i^{\pi^\prime}(\underline{k})$ with
$M(D)=M$. Since
\[\mathcal{L}_i^{M}=\sum_{m=-M}^{M}\bigl(\mathcal{L}_i^M\bigr)^{(m)}
\]
with $M$th order term given by
$\bigl(\mathcal{L}_i^M\bigr)^{(M)}=f_M(z)t(M\upsilon_i)$, 
it follows that
\[D^\prime:=D-\mathcal{L}_i^M\sum_{x\in
X_i}K_M^xt(x)\in\mathbb{D}_i^{\pi^\prime}(\underline{k})
\]
is either zero, or it is nonzero and $M(D^\prime)<M$. By the
induction hypothesis it follows that
$D_x\in\mathbb{C}[\mathcal{L}_i]$ for all $x\in X_i$.
\end{proof}
%%%%%%%%%%%%%%%%%%%%%%%%%%%%%%%%%%%%%%%%%%%%%%%%%%%%%%%%%%%%%%%%%%%%%%%%
\begin{thm}\label{centralizerthm}
Let $\pi^\prime\in L^\prime$ be a
nonzero anti-dominant minuscule or quasi-minuscule co-weight.\\
{\bf i)}
The Harish-Chandra homomorphism 
$\gamma(\underline{k})$ restricts to an algebra isomorphism
\[\gamma(\underline{k}):
\mathbb{D}_{\mathcal{R}}(L^\prime)^{W_0,D_{\pi^\prime}}
\overset{\sim}{\longrightarrow} A_0^\prime.
\]
{\bf ii)} The map $\beta$ restricts to an algebra isomorphism
\[\beta: A_0^\prime(Y)\overset{\sim}{\longrightarrow}
\mathbb{D}_{\mathcal{R}}(L^\prime)^{W_0,D_{\pi^\prime}}.
\]
\end{thm}
%%%%%%%%%%%%%%%%%%%%%%%%%%%%%%%%%%%%%%%%%%%%%%%%%%%%%%%%%%%%%%%%%%%%%%%%%%
\begin{proof}
{\bf i)} Let $i\in\{1,\ldots,n\}$ 
and write $\mathcal{A}_i^\prime$ for the group algebra
$\mathbb{C}[t(L_i^\prime)]$, with canonical basis 
denoted by $\xi^{\lambda^\prime}$
($\lambda^\prime\in L_i^\prime$). 
The group algebra $\mathcal{A}_i^\prime$ is a $W_{0,\{a_i\}}$-module algebra
with action defined by $s_i(\xi^{\lambda^\prime}):=\xi^{s_i\lambda^\prime}$
($\lambda^\prime\in L_i^\prime$). It contains $A^\prime$ as 
$W_{0,\{a_i\}}$-module
algebra.

The map $\delta_{\{a_i\}}(\underline{k})$ 
(see \eqref{deltaFk}) extends to an algebra homomorphism
\[\delta_{\{a_i\}}(\underline{k}):
\mathbb{D}_{\mathcal{R}_{\{a_i\}}}(L_i^\prime)\rightarrow
\mathcal{A}_i^\prime,
\]
defined by the same formula 
\eqref{deltaFk} (with the finite sum over $\lambda^\prime$ now
running over $L_i^\prime$). By a direct computation we have
\[\delta_{\{a_i\}}(\underline{k})(\mathcal{L}_i)=
\xi^{\alpha_i^\vee/2}+\xi^{-\alpha_i^\vee/2}
\]
unless we are in case {\bf c} and $a_i\in W_0\varphi^\vee$, in
which case we have
\[\delta_{\{a_i\}}(\underline{k})(\mathcal{L}_i)=
\xi^{a_i}+\xi^{-a_i}-q^{\kappa_1^\prime}
-q^{-\kappa_1^\prime}.
\]
Combined with Proposition \ref{commutantprop}, we conclude that
$\delta_{\{a_i\}}(\underline{k})$ maps the centralizer subalgebra
$\mathbb{D}_i^{\pi^\prime}(\underline{k})$ onto the subalgebra
of $W_{0,\{a_i\}}$-invariant elements in $\mathcal{A}_i^\prime$.

Fix $D\in\mathbb{D}_{\mathcal{R}}(L^\prime)^{W_0,D_{\pi^\prime}}$.
Then
$\gamma_{\{a_i\}}(\underline{k})(D)\in
\mathbb{D}_i^{\pi^\prime}(\underline{k})\cap
\mathbb{D}_{\mathcal{R}_{\{a_i\}}}(L^\prime)$ by 
Lemma \ref{gammaFklemma}, hence
the constant term
\[
\gamma(\underline{k})(D)=
\delta_{\{a_i\}}(\underline{k})\bigl(\gamma_{\{a_i\}}(\underline{k})(D)\bigr)
\in A^\prime
\]
of $D$ is $W_{0,\{a_i\}}$-invariant 
for all $i\in\{1,\ldots,n\}$ by the previous paragraph. It
follows that $\gamma(\underline{k})(D)\in A_0^\prime$. 
Proposition \ref{Commutativityprop}
now completes the proof of part {\bf i)}.\\
Part {\bf ii)} of the proposition follows 
from part {\bf i)} and Proposition \ref{gammainverse}.
\end{proof}
%%%%%%%%%%%%%%%%%%%%%%%%%%%%%%%%%%%%%%%%%%%%%%%%%%%%%%%%%%%%%%%%%%%%%%%%%%
We now have the following stronger version of Corollary \ref{ok}.
%%%%%%%%%%%%%%%%%%%%%%%%%%%%%%%%%%%%%%%%%%%%%%%%%%%%%%%%%%%%%%%%%%%%%%%%%%%%
\begin{cor}\label{ok2}
We have
\begin{equation}\label{okformula}
\mathbb{D}_{\mathcal{R}}(L^\prime)^{D_{\pi^\prime}}=
\mathbb{D}_{\mathcal{R}}(L^\prime)^{W_0,D_{\pi^\prime}}
\end{equation}
if and only if $A_0$ is a 
$\mathbb{D}_{\mathcal{R}}(L^\prime)^{D_{\pi^\prime}}$-submodule of $\overline{M}$.
\end{cor}
%%%%%%%%%%%%%%%%%%%%%%%%%%%%%%%%%%%%%%%%%%%%%%%%%%%%%%%%%%%%%%%%%%%%%%%%%%
\begin{proof}
By \cite[\S 6.4]{M}, the difference 
operator $D_p|_{A_0}$ is an endomorphism of $A_0$ for all $p\in A_0^\prime$.
The previous theorem thus implies that $A_0$ is a
$\mathbb{D}_{\mathcal{R}}(L^\prime)^{W_0,D_{\pi^\prime}}$-submodule
of $\overline{M}$. The result follows now directly from Corollary \ref{ok}.
\end{proof}
%%%%%%%%%%%%%%%%%%%%%%%%%%%%%%%%%%%%%%%%%%%%%%%%%%%%%%%%%%%%%%%%%%%%%%%

%%%%%%%%%%%%%%%%%%%%%%%%%%%%%%%%%%%%%%%%%%%%%%%%%%%%%%%%%%%
%%                                                       %%
%%      Harish-Chandra series                            %%
%%                                                       %%
%%%%%%%%%%%%%%%%%%%%%%%%%%%%%%%%%%%%%%%%%%%%%%%%%%%%%%%%%%%

\section{Harish-Chandra series}\label{HCsection}

%%%%%%%%%%%%%%%%%%%%%%%%%%%%%%%%%%%%%%%%%%%%%%%%%%%%%%%%%%%
%%                                                       %%
%% Harish-Chandra series with formal spectral parameter  %%
%%                                                       %%
%%%%%%%%%%%%%%%%%%%%%%%%%%%%%%%%%%%%%%%%%%%%%%%%%%%%%%%%%%%

\subsection{Harish-Chandra series with formal spectral
parameter}\label{HCformalsub}

In the trigonometric differential case, Harish-Chandra series are power
series solutions to the differential analogues of the commuting
difference operators $D_p$ ($p\in A_0^\prime$), see \cite[\S 4.2]{HS}
and references therein. In this subsection we construct the natural difference
analogue of the Harish-Chandra series with formal spectral parameter.
For $R$ of type $A$, the difference analogues of
the Harish-Chandra series were considered in \cite[\S 6]{EK} in
the context of weighted traces of quantum group intertwiners, 
see also \cite[\S 9]{EV}.

Let $\mathcal{B}^\prime\subseteq\mathcal{Q}^\prime$ be an $A^\prime$-submodule.
We write $\mathcal{B}^\prime[[z^{-\Delta}]]$ for the
$A^\prime$-module of formal power series
\begin{equation}\label{F}
F(z,\xi)=\sum_{x\in\mathbb{Z}_{\geq 0}\Delta}f_x(\xi)z^{-x},\qquad
f_x(\xi)\in\mathcal{B}^\prime.
\end{equation}
The following lemma is easily checked.

%%%%%%%%%%%%%%%%%%%%%%%%%%%%%%%%%%%%%%%%%%%%%%%%%%%%%%%%%%%%%%%%
\begin{lem}\label{actionlemma}
The canonical $\mathbb{C}[[z^{-\Delta}]]$-action on
$\mathcal{B}^\prime[[z^{-\Delta}]]$ uniquely extends to an
action of $\mathbb{D}_{\mathbb{C}[[z^{-\Delta}]]}(L^\prime)$
on the $A^\prime$-module $\mathcal{B}^\prime[[z^{-\Delta}]]$ by
\[t(\lambda^\prime)\bigl(F(z,\xi)\bigr)=\sum_{x\in\mathbb{Z}_{\geq
0}\Delta}f_x(\xi)\xi^{-\lambda^\prime}q^{\langle\lambda^\prime,x\rangle}z^{-x},\qquad
\lambda^\prime\in L^\prime,
\]
with $F(z,\xi)\in\mathcal{Q}^\prime[[z^{-\Delta}]]$ given by \eqref{F}.
\end{lem}
%%%%%%%%%%%%%%%%%%%%%%%%%%%%%%%%%%%%%%%%%%%%%%%%%%%%%%%%%%%%%%%%%%%

Let $\pi^\prime\in L^\prime$ be a nonzero anti-dominant minuscule or
quasi-minuscule co-weight.
Recall from Proposition \ref{Commutativityprop} that the
centralizer algebra
$\mathbb{D}_{\mathcal{R}}(L^\prime)^{D_{\pi^\prime}}$
is a commutative algebra of difference operators containing the subalgebra
$\beta(A_0^\prime(Y))=\mathbb{D}_{\mathcal{R}}(L^\prime)^{W_0,D_{\pi^\prime}}$
of difference operators $D_p$ ($p\in A_0^\prime$). Recall the 
$\rho_{\underline{k}^\prime}$-twist
\eqref{rhotwist}.

%%%%%%%%%%%%%%%%%%%%%%%%%%%%%%%%%%%%%%%%%%%%%%%%%%%%%%%%%%%%%%%%%%%%
\begin{thm}\label{universalHC}
There exists a unique
\begin{equation}\label{seriesform}
\Phi(z,\xi)=\sum_{x\in\mathbb{Z}_{\geq 0}\Delta}
\Gamma_x(\xi)z^{-x}\in\mathcal{Q}^\prime[[z^{-\Delta}]],\qquad
\Gamma_x(\xi)\in\mathcal{Q}^\prime
\end{equation}
normalized by $\Gamma_0(\xi)=1$ and satisfying the difference equations
\[D\bigl(\Phi(z,\xi)\bigr)=
\widetilde{\bigl(\gamma(\underline{k})(D)\bigr)}(\xi)\Phi(z,\xi),\qquad
\forall\, D\in \mathbb{D}_{\mathcal{R}}(L^\prime)^{D_{\pi^\prime}}
\]
with respect to the action from Lemma \ref{actionlemma}. In
particular,
$\Phi(z,\xi)$ satisfies the difference equations
\[D_p\bigl(\Phi(z,\xi)\bigr)=\widetilde{p}(\xi)\Phi(z,\xi),
\qquad \forall\,p(\xi)\in
A_0^\prime.
\]
\end{thm}
%%%%%%%%%%%%%%%%%%%%%%%%%%%%%%%%%%%%%%%%%%%%%%%%%%%%%%%%%%%%%%%%%%%%
\begin{proof}
Let $D\in \mathbb{D}_{\mathcal{R}}(L^\prime)^{D_{\pi^\prime}}\subseteq
\mathbb{D}_{\mathbb{C}[[z^{-\Delta}]]}(L^\prime)^{D_{\pi^\prime}}$,
written out explicitly as
\[D=\sum_{\lambda^\prime\in L^\prime}
\Bigl(\sum_{y\in\mathbb{Z}_{\geq 0}\Delta}
d_{\lambda^\prime}(y)z^{-y}\Bigr)t(\lambda^\prime)
\]
with $d_{\lambda^\prime}(y)\in\mathbb{C}$ 
and with the first sum over finitely many
$\lambda^\prime\in L^\prime$. We use the 
shorthand notation $r_D(\xi):=\bigl(\gamma(\underline{k})(D)\bigr)(\xi)\in
A^\prime$ for the constant term of $D$. By the definition of
$\gamma(\underline{k})$ (see Lemma \ref{constanttermmap}) we then have
\[\widetilde{r}_D(\xi)=\sum_{\lambda^\prime\in
L^\prime}d_{\lambda^\prime}(0)\xi^{-\lambda^\prime}.
\]

A direct computation now shows that 
$D\bigl(\Phi(z,\xi)\bigr)=\widetilde{r}_D(\xi)\Phi(z,\xi)$,
with $\Phi(z,\xi)\in\mathcal{Q}^\prime[[z^{\Delta}]]$ a series
of the form \eqref{seriesform}, if and only if
\begin{equation}\label{recuruniversal}
\bigl(\widetilde{r}_D(\xi)-
t(x)\bigl(\widetilde{r}_D(\xi)\bigr)\bigr)\Gamma_x(\xi)
=\sum_{\stackrel{\lambda^\prime\in L^\prime}{0\leq
y<x}}d_{\lambda^\prime}(x-y)q^{\langle
\lambda^\prime,y\rangle}\xi^{-\lambda^\prime}\Gamma_y(\xi)
\end{equation}
for all $x\in\mathbb{Z}_{\geq 0}\Delta$.

We now explore the recurrence relations \eqref{recuruniversal}
first for the Macdonald difference operator $D_{\pi^\prime}$, in
which case $r_{D_{\pi^\prime}}(\xi)=m_{\pi^\prime}(\xi)$ by 
Proposition \ref{gammainverse}.
We have $t(x)\bigl(\widetilde{m}_{\pi^\prime}(\xi)\bigr)\not=
\widetilde{m}_{\pi^\prime}(\xi)$
for all $x\in\mathbb{Z}_{\geq 0}\Delta\setminus\{0\}$ 
since $W_0$ acts irreducibly on $V$.
Hence for $D=D_{\pi^\prime}$ the
recurrence relations \eqref{recuruniversal} has a unique solution
$\Gamma_x(\xi)\in\mathcal{Q}^\prime$
($x\in\mathbb{Z}_{\geq 0}\Delta$) normalized by $\Gamma_0(\xi)=1$,
and the
resulting formal power series
\[\Phi(z,\xi)=\sum_{x\in\mathbb{Z}_{\geq 0}\Delta}\Gamma_x(\xi)z^{-x}
\in\mathcal{Q}^\prime[[z^{-\Delta}]]
\]
thus satisfies 
$D_{\pi^\prime}\bigl(\Phi(z,\xi)\bigr)=
\widetilde{m}_{\pi^\prime}(\xi)\Phi(z,\xi)$.

Fix now 
$D\in\mathbb{D}_{\mathcal{R}}(L^\prime)^{D_{\pi^\prime}}$ arbitrary and consider
\[\Phi^\prime(z,\xi):=D\bigl(\Phi(z,\xi)\bigr)-\widetilde{r}_D(\xi)\Phi(z,\xi)
\in\mathcal{Q}^\prime[[z^{-\Delta}]]
\]
with $\Phi(z,\xi)$ as defined in the previous paragraph.
Since $\lbrack D, D_{\pi^\prime}\rbrack=0$ we have
\[D_{\pi^\prime}\bigl(\Phi^\prime(z,\xi)\bigr)=
\widetilde{m}_{\pi^\prime}(\xi)
\Phi^\prime(z,\xi),
\]
hence the coefficients $\Gamma_x^\prime(\xi)\in\mathcal{Q}^\prime$ in
the expansion $\Phi^\prime(z,\xi)=
\sum_{x\in\mathbb{Z}_{\geq 0}\Delta}\Gamma_x^\prime(\xi)z^{-x}$
satisfies the recurrence relations \eqref{recuruniversal} for
$D=D_{\pi^\prime}$. By a direct computation (analogous e.g. to 
the computations
in the proof of Lemma \ref{ltlem}) 
we furthermore have $\Gamma_0^\prime(\xi)=0$,
hence we conclude that $\Gamma_x^\prime(\xi)=0$
for all $x\in\mathbb{Z}_{\geq 0}\Delta$. Consequently
$D\bigl(\Phi(z,\xi)\bigr)=\widetilde{r}_D(\xi)\Phi(z,\xi)$.
Noting finally that $r_{D_p}(\xi)=p(\xi)$ for $p\in A_0^\prime$ 
by Proposition \ref{gammainverse},
we obtain the desired results.
\end{proof}
%%%%%%%%%%%%%%%%%%%%%%%%%%%%%%%%%%%%%%%%%%%%%%%%%%%%%%%%%%%%%%%%%%%%%%%%%%%%%
\begin{rema}\label{oneenough}
It follows from (the proof of) 
Theorem \ref{universalHC} that the Harish-Chandra
series $\Phi(z,\xi)\in\mathcal{Q}^\prime[[z^{-\Delta}]]$
is uniquely characterized up to $\mathcal{Q}^\prime$-multiples
by the single difference equation
\[D_{{\pi^\prime}}\bigl(\Phi(z,\xi)\bigr)=
\widetilde{m}_{\pi^\prime}(\xi)\Phi(z,\xi)
\]
involving the Macdonald difference operator $D_{\pi^\prime}$.
\end{rema}
%%%%%%%%%%%%%%%%%%%%%%%%%%%%%%%%%%%%%%%%%%%%%%%%%%%%%%%%%%%%%%%%%%%%%%%%%%%%%%

%%%%%%%%%%%%%%%%%%%%%%%%%%%%%%%%%%%%%%%%%%%%%%%%%%%%%%%%%%%%%%%%%%
%%                                                              %%
%%   Harish-Chandra series with specialized spectral parameter  %%
%%                                                              %%
%%%%%%%%%%%%%%%%%%%%%%%%%%%%%%%%%%%%%%%%%%%%%%%%%%%%%%%%%%%%%%%%%%

\subsection{Harish-Chandra series with specialized spectral
parameter}\label{HCspecialsub}

We view $\mathcal{Q}^\prime$ as rational trigonometric
functions on $V_{\mathbb{C}}$ using $\xi^{\lambda^\prime}(v)=
q^{\langle \lambda^\prime,v\rangle}$
for $v\in V_{\mathbb{C}}$ and 
$\lambda^\prime\in L^\prime$, cf. \eqref{planewave}.
We write $q=e^{-2\pi\sigma}$ with $\sigma\in\mathbb{R}_{>0}$, and
\[Z=\{\lambda\in V \, | \, 
\langle \lambda^\prime,\lambda\rangle
\in\mathbb{Z}\quad\forall\,\lambda^\prime\in L^\prime\}
\]
for the lattice in $V$ dual to $L^\prime$. Note that $Z=Q(R)$
for case {\bf a} and {\bf b}, and $Z=P(R)$ for case {\bf c}.
%%%%%%%%%%%%%%%%%%%%%%%%%%%%%%%%%%%%%%%%%%%%%%%%%%%%%%%%%%%%%%%%%%%%%%%%
\begin{lem}\label{sing}
The singularities of the coefficients $\Gamma_x(\xi)\in\mathcal{Q}^\prime$
\textup{(}$x\in\mathbb{Z}_{\geq 0}\Delta$\textup{)} of the
Harish-Chandra series $\Phi(z,\xi)$ are contained in the subset
\[\mathcal{D}_{\underline{k}}=\bigcup_{w\in W_0}\{\lambda\in 
V_{\mathbb{C}} \, | \,
\lambda-w\cdot\lambda\in \mathbb{Z}_{\geq 0}\Delta\setminus
\{0\}+\sqrt{-1}Z/\sigma\}\subset V_{\mathbb{C}},
\]
where the dot-action of $W_0$ on $V_{\mathbb{C}}$ is defined by
$w\cdot\lambda=w(\lambda+\rho_{\underline{k}^\prime})-\rho_{\underline{k}^\prime}$
for $w\in W_0$ and $\lambda\in V_{\mathbb{C}}$.
\end{lem}
%%%%%%%%%%%%%%%%%%%%%%%%%%%%%%%%%%%%%%%%%%%%%%%%%%%%%%%%%%%%%%%%%%%%%%%%%%
\begin{proof}
For $\lambda,\mu\in V_{\mathbb{C}}$ we have
\begin{equation}\label{ischar}
\widetilde{p}(\lambda)=\widetilde{p}(\mu)\quad \forall\, p(\xi)\in
A_0^\prime \,\,\,\Leftrightarrow \,\,\, \lambda\in
W_0\cdot\mu+\sqrt{-1}Z/\sigma.
\end{equation}
We now prove the lemma by induction 
to the height $\sum_{\alpha\in\Delta}l_\alpha\in\mathbb{Z}_{\geq 0}$
of an element $x=\sum_{\alpha\in\Delta}l_\alpha\alpha\in\mathbb{Z}_{\geq
0}\Delta$. If $\lambda\in V_{\mathbb{C}}\setminus\mathcal{D}_{\underline{k}}$
and $x\in\mathbb{Z}_{\geq 0}\Delta\setminus\{0\}$, then
\eqref{ischar} implies the existence of a $p\in A_0^\prime$ such
that $\widetilde{p}(\lambda)\not=\widetilde{p}(\lambda-x)$.
By \eqref{recuruniversal} applied to $D=D_p$ the regularity of
$\Gamma_x(\xi)$ at $\lambda$ is implied by the regularity of $\Gamma_y(\xi)$
at $\lambda$ for elements $y\in\mathbb{Z}_{\geq 0}\Delta$ with height 
strictly smaller than $x$.
\end{proof}
%%%%%%%%%%%%%%%%%%%%%%%%%%%%%%%%%%%%%%%%%%%%%%%%%%%%%%%%%%%%%%%%%%%%%%%%%%%
We denote $\mathcal{B}^\prime_{\underline{k}}$
for the $A^\prime$-submodule of $\mathcal{Q}^\prime$
consisting of the rational functions $f(\xi)\in\mathcal{Q}^\prime$
with singularities contained in $\mathcal{D}_{\underline{k}}$.
The previous lemma shows that $\Phi(z,\xi)\in
\mathcal{B}^\prime_{\underline{k}}[[z^{-\Delta}]]$.
The following lemma will allow us to derive difference 
equations for the Harish-Chandra
series when the spectral parameter is 
specialized to an element in the open and dense subset
$V_{\mathbb{C}}\setminus \mathcal{D}_{\underline{k}}$ of $V_{\mathbb{C}}$.
Recall the $\mathbb{D}_{\mathbb{C}[[z^{-\Delta}]]}(L^\prime)$-module 
$\overline{M}$
from Subsection \ref{DRO}.

%%%%%%%%%%%%%%%%%%%%%%%%%%%%%%%%%%%%%%%%%%%%%%%%%%%%%%%%%%%%%%%%%%%%%%
\begin{lem}\label{iotalemma}
For $\lambda\in V_{\mathbb{C}}\setminus\mathcal{D}_{\underline{k}}$ the
assignment
\[\sum_{x\in \mathbb{Z}_{\geq 0}\Delta}f_x(\xi)z^{-x}\mapsto
\sum_{x\in\mathbb{Z}_{\geq
0}\Delta}f_x(\lambda)z^{\lambda-x},\qquad
(f_x(\xi)\in\mathcal{B}^\prime_{\underline{k}}),
\]
defines a morphism
$\iota_\lambda:\mathcal{B}^\prime_{\underline{k}}[[z^{-\Delta}]]\rightarrow
\overline{M}$ of $\mathbb{D}_{\mathbb{C}[[z^{-\Delta}]]}(L^\prime)$-modules.
\end{lem}
%%%%%%%%%%%%%%%%%%%%%%%%%%%%%%%%%%%%%%%%%%%%%%%%%%%%%%%%%%%%%%%%%%%%%%%
\begin{proof}
Direct verification.
\end{proof}
%%%%%%%%%%%%%%%%%%%%%%%%%%%%%%%%%%%%%%%%%%%%%%%%%%%%%%%%%%%%%%%%%%%%%%
\begin{thm}\label{HCspecial}
The Harish-Chandra series with spectral
parameter $\lambda\in V_{\mathbb{C}}\setminus\mathcal{D}_{\underline{k}}$,
defined by
\[\Phi_\lambda(z):=\iota_\lambda\bigl(\Phi(z,\xi)\bigr)=
\sum_{x\in\mathbb{Z}_{\geq
0}\Delta}\Gamma_x(\lambda)z^{\lambda-x}\in\overline{M},
\]
satisfies
\begin{equation}\label{deq}
D\bigl(\Phi_\lambda(z)\bigr)=\widetilde{\bigl(\gamma(\underline{k})(D)\bigr)}
(\lambda)\Phi_\lambda(z),\qquad
\forall\, D\in\mathbb{D}_{\mathcal{R}}(L^\prime)^{D_{\pi^\prime}},
\end{equation}
hence in particular
\[
D_p\bigl(\Phi_\lambda(z)\bigr)=\widetilde{p}(\lambda)\Phi_\lambda(z),\qquad
\forall\, p(\xi)\in A_0^\prime.
\]
The latter system of difference equations, together with the
normalization $\Gamma_0(\lambda)=1$, uniquely characterizes
$\Phi_\lambda(z)$ within the 
$\mathbb{D}_{\mathbb{C}[[z^{-\Delta}]]}(L^\prime)$-submodule
$\mathbb{C}[[z^{-\Delta}]]z^\lambda$ of $\overline{M}$.
\end{thm}
%%%%%%%%%%%%%%%%%%%%%%%%%%%%%%%%%%%%%%%%%%%%%%%%%%%%%%%%%%%%%%%%%%%%%%%%%
\begin{proof}
The first part follows immediately from the results from the previous section.
For uniqueness, note that the coefficients $\Gamma_x(\lambda)$ of
a solution $\sum_{x\in\mathbb{Z}_{\geq 0}\Delta}
\Gamma_x(\lambda)z^{\lambda-x}\in\overline{M}$
of the difference equations \eqref{deq} satisfy the homogeneous recurrence
relations \eqref{recuruniversal} specialized to 
$\lambda\in V_{\mathbb{C}}\setminus
\mathcal{D}_{\underline{k}}$. With a 
similar argument as in the proof of Lemma
\ref{sing} it follows that the coefficients 
$\Gamma_x(\lambda)$ ($x\in\mathbb{Z}_{\geq
0}\Delta$) are determined by $\Gamma_0(\lambda)$.
\end{proof}
%%%%%%%%%%%%%%%%%%%%%%%%%%%%%%%%%%%%%%%%%%%%%%%%%%%%%%%%%%%%%%%%%%%%%%%%%%%%%

%%%%%%%%%%%%%%%%%%%%%%%%%%%%%%%%%%%%%%%%%%%%%%%%%%%%%%%%%%%%%%%%%%
\begin{thm}\label{basis}
Let $\lambda\in V_{\mathbb{C}}$ such that
\begin{equation}\label{genericcondition}
\lambda-w\cdot\lambda\not\in\mathbb{Z}\Delta+\sqrt{-1}Z/\sigma,\qquad
\forall w\in W_0\setminus \{e\}.
\end{equation}
Then $\{\Phi_\mu\,\, | \,\, \mu\in W_0\cdot\lambda+\sqrt{-1}Z/\sigma\}$
is a basis of the common eigenspace
\[\overline{M}_\lambda=\{F(z)\in\overline{M} \,\,\, | \,\,\,
D_p\bigl(F(z)\bigr)=\widetilde{p}(\lambda)F(z)\quad\forall\,p(\xi)\in
A_0^\prime\}.
\]
\end{thm}
%%%%%%%%%%%%%%%%%%%%%%%%%%%%%%%%%%%%%%%%%%%%%%%%%%%%%%%%%%%%%%%%%%%%
\begin{proof}
By \eqref{genericcondition}, $W_0\cdot\lambda+\sqrt{-1}Z/\sigma$
consists of elements from
$V_{\mathbb{C}}\setminus\mathcal{D}_{\underline{k}}$ 
which are pair-wise incomparable
with respect to the dominance order $\leq$.
In particular, the Harish-Chandra series $\Phi_\mu(z)\in\overline{M}$
($\mu\in W_0\cdot\lambda+\sqrt{-1}Z/\sigma)$ are well defined and
linearly independent. Furthermore, $\Phi_\mu(z)\in\overline{M}_\lambda$
($\mu\in W_0\cdot\lambda+\sqrt{-1}Z/\sigma)$ by \eqref{ischar}.

For $F(z)=\sum_\nu K_\nu z^\nu\in \overline{M}$ 
we write $\hbox{Supp}(F(z))=\{\nu\in V_{\mathbb{C}} \, | \,
K_\nu\not=0\}$. We claim that if $F(z)\in\overline{M}_\lambda$ and
if $\nu\in \hbox{Supp}(F(z))$ is a maximal element with respect to 
the dominance order $\leq$,
then $\nu\in W_0\cdot\lambda+\sqrt{-1}Z/\sigma$. Before proving the
claim, we first show that it implies that $F(z)$
is a finite linear combination of the Harish-Chandra 
series $\Phi_\mu(z)$ ($\mu\in W_0\cdot\lambda+\sqrt{-1}Z/\sigma$).

If $0\not=F(z)\in\overline{M}_\lambda$ then we can
choose a maximal element $\nu_1\in\hbox{Supp}(F(z))$ 
with respect to $\leq$ by Zorn's Lemma.
Then $\nu_1\in W_0\cdot\lambda+\sqrt{-1}Z/\sigma$
by the claim, hence $\nu_1\not\in \mathcal{D}_{\underline{k}}$ 
and $\Phi_{\nu_1}(z)\in\overline{M}_\lambda$.
Set $F_1(z)=F(z)-K_{\nu_1}\Phi_{\nu_1}(z)\in\overline{M}_\lambda$.
If $F_1(z)\not=0$ then we choose a maximal element $\nu_2$ in 
$\hbox{Supp}(F_1(z))$, which necessarily satisfies
$\nu_1\not=\nu_2\in W_0\cdot\lambda+\sqrt{-1}Z/\sigma$. 
In particular, $\nu_1$ and $\nu_2$ are incomparable with respect
to $\leq$. We proceed to define
\[F_2(z):=F_1(z)-K_{\nu_2}\Phi_{\nu_2}(z)=
F(z)-K_{\nu_1}\Phi_{\nu_1}(z)-K_{\nu_2}\Phi_{\nu_2}(z)\in\overline{M}_\lambda.
\]
Repeating the above procedure we 
construct $F_m(z)\in\overline{M}_\lambda$ and $\nu_m\in V_{\mathbb{C}}$
from $F_{m-1}(z)\in\overline{M}_\lambda$ inductively,
as long as $F_{m-1}(z)\not=0$. Since the resulting 
spectral parameters $\mu_1,\ldots,\mu_m$ are
pair-wise incomparable with respect to $\leq$,
we have $F_m(z)=0$ for some $m\in\mathbb{Z}_{\geq 0}$
in view of the definition of $\overline{M}$.
Consequently, $F(z)$ is a linear combination 
of the Harish-Chandra series $\Phi_{\nu_j}(z)$ ($j=1,\ldots,m$).

It remains to prove the claim. Let 
$\mu\in\hbox{Supp}(F(z))$ be a maximal element with respect to $\leq$ and set
\[G(z):=F(z)-K_{\mu}z^{\mu}\in\overline{M},
\]
so that $\hbox{Supp}(G(z))=\hbox{Supp}(F(z))\setminus \{\mu\}$.
Fix $p\in A_0^\prime$. Since $D_p(z^\nu)\in z^\nu\mathbb{C}[[z^{-\Delta}]]$ 
for $\nu\in V_{\mathbb{C}}$,
we have $\mu\not\in\hbox{Supp}(D_p(G(z)))$. 
Combined with Lemma \ref{ltlem} and Proposition \ref{gammainverse},
we conclude that the coefficient of $z^\mu$ in
\[
D_p(F(z))=D_p(G(z))+K_\mu D_p(z^\mu)
\]
is $K_\mu\widetilde{p}(\mu)$.
On the other hand, the coefficient of 
$z^{\mu}$ in $D_p(F(z))=\widetilde{p}(\lambda)F(z)$
is $K_\mu \widetilde{p}(\lambda)$.
Hence $\widetilde{p}(\mu)=\widetilde{p}(\lambda)$ for all $p\in
A_0^\prime$, which implies that $\mu\in
W_0\cdot\lambda+\sqrt{-1}Z/\sigma$ by \eqref{ischar}.
\end{proof}
%%%%%%%%%%%%%%%%%%%%%%%%%%%%%%%%%%%%%%%%%%%%%%%%%%%%%%%%%%%%%%%%%%%%%%%%%%
\begin{rema}
The Harish-Chandra series $\Phi_\lambda(z)$ for generic real spectral
values $\lambda\in V$ is contained in the 
$\mathbb{D}_{\mathbb{C}[[z^{-\Delta}]]}(L^\prime)$-submodule
$\overline{M}^{re}\subset \overline{M}$ consisting of the formal series
$F(z)=\sum_{u\in\mathcal{C}}K_uz^u$ in $\overline{M}$
with $\mathcal{C}\subset V$. Suppose that $\lambda\in V$ satisfies
\[\lambda-w\cdot\lambda\not\in\mathbb{Z}\Delta,\qquad\forall\,w\in
W_0\setminus\{e\},
\]
and suppose that the multiplicity label $\underline{k}^\prime$ is real-valued,
so that $\rho_{\underline{k}^\prime}\in V$ and the dot-action of $W_0$ preserves
$V$. Then Theorem \ref{basis} implies that
the common eigenspace $\overline{M}^{re}_\lambda:=
\overline{M}^{re}\cap\overline{M}_\lambda$
is $\#W_0$-dimensional with basis $\{\Phi_{w\cdot\lambda}\}_{w\in
W_0}$, cf. \cite[Cor. 4.2.6]{HS} for the analogous statement in 
the trigonometric differential set-up.
For type A, this result essentially is \cite[Thm. 5, part 2]{EK}.
\end{rema}

%%%%%%%%%%%%%%%%%%%%%%%%%%%%%%%%%%%%%%%%%%%%%%%%%%%%%%%%%%%%%%%%%%%%%%
%%                                                                  %%
%%  Relation to Macdonald polynomials                               %%
%%                                                                  %%
%%%%%%%%%%%%%%%%%%%%%%%%%%%%%%%%%%%%%%%%%%%%%%%%%%%%%%%%%%%%%%%%%%%%%%

\subsection{Relation to Macdonald polynomials}\label{polynomialsection}

The properties \eqref{triangular} and 
\eqref{ischar} lead to the following well-known
definition of the Macdonald polynomials \cite{M}, \cite{M2}, \cite{Ko}
(known as Koornwinder \cite{Ko} polynomials for case {\bf c}).
%%%%%%%%%%%%%%%%%%%%%%%%%%%%%%%%%%%%%%%%%%%%%%%%%%%%%%%%%%%%%%%%%%%%%%%%
\begin{defi}\label{Macdonalddef}
Suppose that
\begin{equation}\label{Macdcond}
\lambda\not\in W_0\cdot\mu+\sqrt{-1}Z/\sigma\qquad
\forall\,\lambda,\mu\in L_{++}: \lambda\not=\mu.
\end{equation}
The monic Macdonald polynomial 
$P_\lambda(z)\in A_0$ of degree $\lambda\in L_{++}$ is
the unique $W_0$-invariant Laurent polynomial satisfying
\[P_\lambda(z)=m_\lambda(z)+\sum_{\mu\in L_{++}:
\mu<\lambda}k_{\lambda,\mu}m_{\mu}(z)
\]
for certain coefficients $k_{\lambda,\mu}\in\mathbb{C}$ and
satisfying the system of difference equations
\[D_p\bigl(P_\lambda(z)\bigr)=\widetilde{p}(\lambda)P_\lambda(z),\qquad
\forall\,p(\xi)\in A_0^\prime.
\]
\end{defi}
%%%%%%%%%%%%%%%%%%%%%%%%%%%%%%%%%%%%%%%%%%%%%%%%%%%%%%%%%%%%%%%%%%%%%
Note that if the multiplicity label $\underline{k}^\prime$ is real 
valued then the conditions
\eqref{Macdcond} reduce to
\begin{equation}\label{Macdcondreal}
\lambda\not\in W_0\cdot\mu\qquad\quad 
\forall\,\lambda,\mu\in L_{++}: \lambda\not=\mu.
\end{equation}
For instance, 
if $k_a^\prime\geq 0$ for all $a\in S^\prime$, then $\rho_{\underline{k}^\prime}\in
\overline{V}_+$, hence \eqref{Macdcondreal} is satisfied.

%%%%%%%%%%%%%%%%%%%%%%%%%%%%%%%%%%%%%%%%%%%%%%%%%%%%%%%%%%%%%%%%%%%%%%
\begin{prop}\label{same}
Suppose that \eqref{Macdcond} is satisfied.
For $\lambda\in L_{++}$ not contained in $\mathcal{D}_{\underline{k}}$
we have
\[\Phi_\lambda(z)=P_\lambda(z).
\]
\end{prop}
%%%%%%%%%%%%%%%%%%%%%%%%%%%%%%%%%%%%%%%%%%%%%%%%%%%%%%%%%%%%%%%%%%
\begin{proof}
By the assumptions, the Harish-Chandra series $\Phi_\lambda(z)$
and the Macdonald polynomial $P_\lambda(z)$ are well defined.
Furthermore, both $\Phi_\lambda(z)$ and $P_\lambda(z)$ are elements from
$\overline{M}_\lambda$ with the coefficient of $z^\lambda$ equal to
one. Hence $\Phi_\lambda(z)=P_\lambda(z)$ by Theorem
\ref{HCspecial}.
\end{proof}
%%%%%%%%%%%%%%%%%%%%%%%%%%%%%%%%%%%%%%%%%%%%%%%%%%%%%%%%%%%%%%%%%%%%%%%
We now return to the centralizers of Macdonald difference operators.
%%%%%%%%%%%%%%%%%%%%%%%%%%%%%%%%%%%%%%%%%%%%%%%%%%%%%%%%%%%%%%%%%%%%%%
\begin{cor}\label{algintcor}
Let $\pi^\prime\in L^\prime$ be a nonzero anti-dominant minuscule
or quasi-minuscule co-weight.
Suppose that \eqref{Macdcond} is 
satisfied and that $L_{++}\subset V_{\mathbb{C}}\setminus
\mathcal{D}_{\underline{k}}$. Then
\[\mathbb{D}_{\mathcal{R}}(L^\prime)^{D_{\pi^\prime}}=
\mathbb{D}_{\mathcal{R}}(L^\prime)^{W_0,D_{\pi^\prime}},
\]
hence the centralizer of $D_{\pi^\prime}$ 
in $\mathbb{D}_{\mathcal{R}}(L^\prime)$
consists only of the Cherednik-Macdonald 
difference operators $D_p$ \textup{(}$p(\xi)\in
A_0^\prime$\textup{)}.
\end{cor}
%%%%%%%%%%%%%%%%%%%%%%%%%%%%%%%%%%%%%%%%%%%%%%%%%%%%%%%%%%%%%%%%%%%%%%%
\begin{proof}
The second part of the statement follows 
from  Theorem \ref{centralizerthm}{\bf ii)}.
For the first statement it suffices to show 
that $A_0\subset \overline{M}$ is a
$\mathbb{D}_{\mathcal{R}}(L^\prime)^{D_{\pi^\prime}}$-submodule
in view of Corollary \ref{ok2}.
Since the Macdonald polynomials $P_\lambda(z)$ ($\lambda\in
L_{++}$) form a linear basis of $A_0$, it suffices to note that
\begin{equation*}
\begin{split}
D\bigl(P_\lambda(z)\bigr)&=D\bigl(\Phi_\lambda(z)\bigr)\\
&=\widetilde{\bigl(\gamma(\underline{k})(D)\bigr)}(\lambda)\Phi_\lambda(z)\\
&=\widetilde{\bigl(\gamma(\underline{k})(D)\bigr)}(\lambda)P_\lambda(z)
\end{split}
\end{equation*}
for $D\in\mathbb{D}_{\mathcal{R}}(L^\prime)^{D_{\pi^\prime}}$ and $\lambda\in
L_{++}$, where the first and the last equality follow from
Proposition \ref{same} and the second equality follows from Theorem
\ref{HCspecial}.
\end{proof}
%%%%%%%%%%%%%%%%%%%%%%%%%%%%%%%%%%%%%%%%%%%%%%%%%%%%%%%%%%%%%%%%%%%%%%%%%
As an example, consider a multiplicity label 
$\underline{k}^\prime: S^\prime\rightarrow
\mathbb{R}^\times$ such that
\begin{equation}\label{condexpl}
\rho_{\underline{k}^\prime}-w(\rho_{\underline{k}^\prime})=
\sum_{\stackrel{\alpha\in
R^+:}{w^{-1}\alpha\in R^-}}
\underline{k}^\prime(\alpha^\vee)\alpha\not\in L,\qquad
\forall\, w\in W_0\setminus\{e\}
\end{equation}
(the alternative expression for 
$\rho_{\underline{k}^\prime}-w(\rho_{\underline{k}^\prime})$
follows from \cite[(1.5.3)]{M}). Then $L_{++}$ is contained 
in $V_{\mathbb{C}}\setminus
\mathcal{D}_{\underline{k}}$ and \eqref{Macdcondreal} is
satisfied, hence Corollary \ref{algintcor} implies that
$\mathbb{D}_{\mathcal{R}}(L^\prime)^{D_{\pi^\prime}}=
\mathbb{D}_{\mathcal{R}}(L^\prime)^{W_0,D_{\pi^\prime}}$.
For example, for case {\bf a} we have $L=P(R)$ so condition
\eqref{condexpl} implies $\underline{k}^\prime(\alpha^\vee)\not\in\mathbb{Z}$
for all $\alpha\in R$.

Examples for which the conditions of Corollary \ref{algintcor} are
violated, are discussed in the next subsection.

%%%%%%%%%%%%%%%%%%%%%%%%%%%%%%%%%%%%%%%%%%%%%%%%%%%%%%%%%%%%%%%%%%%%%%
%%                                                                  %%
%%  Relation to Baker-Akhiezer functions                            %%
%%                                                                  %%
%%%%%%%%%%%%%%%%%%%%%%%%%%%%%%%%%%%%%%%%%%%%%%%%%%%%%%%%%%%%%%%%%%%%%%

\subsection{Relation to Baker-Akhiezer functions}\label{BAsub}

For special discrete values of the multiplicity label $\underline{k}$,
the commuting Cherednik-Macdonald difference operators $D_{p}$ 
($p(\xi)\in A_0^\prime$)
are algebraically integrable in the sense of e.g. \cite{C} and \cite{ES}.
For such multiplicity labels, Chalykh \cite{C} defines and studies for 
case {\bf a} and case {\bf c}
eigenfunctions of the difference operators $D_{p}$ ($p(\xi)\in A_0^\prime$)
called Baker-Akhiezer functions (see \cite{ES} for $R$ of type $A$).
In this subsection we relate the Harish-Chandra series
to the Baker-Akhiezer functions for case {\bf a}. 
Case {\bf c} can be treated in a similar fashion.

We assume throughout the remainder of the subsection that we
are in case {\bf a}, so that $\Delta=\{\alpha_1,\ldots,\alpha_n\}$
is the basis of $R$, $\mathbb{Z}\Delta=Q(R)$ respectively 
$\mathbb{Z}_{\geq 0}\Delta=Q_+(R)$
is the root lattice respectively 
its cone of positive integral linear combinations
of positive roots, and $(L,L^\prime)=(P(R),P(R^\vee))$.
We furthermore assume throughout the remainder of the subsection that the
multiplicity label $\underline{k}$ satisfies
\begin{equation}\label{algintcondition}
k_\alpha\in\mathbb{Z}_{\leq 0},\qquad \forall\,\alpha\in R.
\end{equation}
Note that \eqref{algintcondition} implies that
$\rho_{\underline{k}^\prime}\in P(R)$ and
\[\rho_{\underline{k}}^\prime:=\frac{1}{2}\sum_{\alpha\in
R^+}k_\alpha\alpha^\vee\in P(R^\vee).
\]

With this choice of multiplicity label, the Macdonald
polynomials $P_\lambda(z)$ are not defined for low degree
$\lambda\in L_{++}$ (see \cite[\S 5.4]{C}), and the
Harish-Chandra series $\Phi_\lambda(z)$ ($\lambda\in L_{++}$)
are not well defined for large degree $\lambda\in L_{++}$ (specifically, for
$\lambda\in L_{++}$ such that $\lambda+\rho_{\underline{k}^\prime}\in V_+$).
In particular, Proposition \ref{same} and
Corollary \ref{algintcor} are no longer valid. The theory in this
set-up requires a completely different approach, which was developed
by Chalykh in \cite{C}. In this approach a key role is played
by the normalized Baker-Akhiezer function $\psi_\lambda(z)$, whose
definition we now shortly recall from \cite{C}.

Our present notations are matched with the ones from \cite{C} as
follows; the parameters $(q,\tau_\alpha)$ correspond to $(q^2,t_\alpha)$ 
in \cite{C},
and the minuscule or quasi-minuscule co-weight $\pi^\prime$ corresponds
to $-\pi$  in \cite{C} (with these correspondences, our Macdonald difference
operator $D_{\pi^\prime}$ turns into the Macdonald difference operator $D^{\pi}$
from \cite[\S 2.2]{C}). Set
\[\mathcal{N}=\{\sum_{\alpha\in R^+}l_\alpha\alpha\,\,\, | \,\,\,
l_\alpha\in\mathbb{Z}\,\hbox{ and } 0\leq l_\alpha\leq -k_\alpha\,\,\,\forall\, 
\alpha\in R\}\subset Q_+(R).
\]
Chalykh's \cite[Thm.\,4.7]{C} Baker-Akhiezer function
$\psi^\vee(\lambda,z)$ associated to the co-root lattice $R^\vee$
is now defined as follows.
%%%%%%%%%%%%%%%%%%%%%%%%%%%%%%%%%%%%%%%%%%%%%%%%%%%%%%%%%%%%%%%%%%%%%%%%%%
\begin{defi}[\cite{C}]\label{Bafunction}
The Baker-Akhiezer functions $\psi_\lambda(z)\in M$ \textup{(}$\lambda\in
V_{\mathbb{C}}$\textup{)} are the unique functions of the form
\begin{equation}\label{supportBA}
\psi_\lambda(z)=\sum_{x\in\mathcal{N}}
K^{BA}_x(\lambda)z^{\lambda-\rho_{\underline{k}^\prime}-x}\qquad
(K_x^{BA}(\xi)\in A^\prime),
\end{equation}
satisfying the equalities
\begin{equation}\label{equalitypropertiesBA}
\psi_\lambda(v+r\alpha^\vee/2)=\psi_\lambda(v-r\alpha^\vee/2)\quad
\hbox{ for } q^{\langle v,\alpha\rangle}=1
\end{equation}
if $\alpha\in R$ and $r=1,\ldots,-k_\alpha$, and normalized by
\begin{equation}\label{normalizationBA}
K_0^{BA}(\xi)=
\xi^{\rho_{\underline{k}}^\prime}
\prod_{\alpha\in
R^+}\prod_{j=1}^{-k_\alpha}\bigl(q^{j/2}-q^{-j/2}\xi^{\alpha^\vee}\bigr).
\end{equation}
\end{defi}
%%%%%%%%%%%%%%%%%%%%%%%%%%%%%%%%%%%%%%%%%%%%%%%%%%%%%%%%%%%%%%%%%%%%%%%%%%%%%
A key property of the Baker-Akhiezer functions $\psi_\lambda(z)\in M\subset 
\overline{M}$
is the fact that
\begin{equation}\label{eigenvalueeqnBA}
D_{\pi^\prime}\bigl(\psi_\lambda(z)\bigr)=m_{\pi^\prime}(-\lambda)\psi_\lambda(z)
\end{equation}
for a nonzero antidominant minuscule or 
quasi-minuscule co-weight $\pi^\prime\in L^\prime$,
cf. \cite[Thm. 3.7(iv) \& Thm. 4.7]{C}.
Comparing with the properties of the Harish-Chandra series
$\Phi_\lambda(z)=\sum_{x\in Q_+(R)}\Gamma_x(\lambda)z^{\lambda-x}\in\overline{M}$
($\lambda\in V_{\mathbb{C}}\setminus\mathcal{D}_{\underline{k}}$)
we obtain the following result.

%%%%%%%%%%%%%%%%%%%%%%%%%%%%%%%%%%%%%%%%%%%%%%%%%%%%%%%%%%%%%%%%%%%%%%%
\begin{prop}\label{Baconn}
For case {\bf a} and for multiplicity labels $\underline{k}$
satisfying \eqref{algintcondition} we have
\[\Phi_\lambda(z)=K_0^{BA}(\lambda+\rho_{\underline{k}^\prime})^{-1}
\psi_{\lambda+\rho_{\underline{k}^\prime}}(z)
\]
if $\lambda\in
V_{\mathbb{C}}\setminus \mathcal{D}_{\underline{k}}$. In
particular, the rational functions $\Gamma_x(\xi)\in\mathcal{Q}^\prime$ 
of the Harish-Chandra
series $\Phi(z,\xi)=\sum_{x\in Q_+(R)}\Gamma_x(\xi)z^{-x}$ 
with formal spectral parameter
satisfy $\Gamma_x(\xi)=0$ if $x\not\in\mathcal{N}$
and
\[
\Gamma_x(\xi)=
t(-\rho_{\underline{k}^\prime})\bigl(K_x^{BA}(\xi)/K_0^{BA}(\xi)\bigr),\qquad
\forall\, x\in\mathcal{N}.
\]
\end{prop}
%%%%%%%%%%%%%%%%%%%%%%%%%%%%%%%%%%%%%%%%%%%%%%%%%%%%%%%%%%%%%%%%%%%%%%%
\begin{proof}
First note that \eqref{eigenvalueeqnBA} implies that 
$\psi_{\lambda+\rho_{\underline{k}^\prime}}(z)\in\overline{M}_\lambda$,
by a similar argument as in the proof of Theorem \ref{universalHC} 
(see also Remark \ref{oneenough}).

Fix now $\lambda\in V_{\mathbb{C}}
\setminus\mathcal{D}_{\underline{k}}$. Then 
$K_0^{BA}(\lambda+\rho_{\underline{k}^\prime})\not=0$
and both $\Phi_\lambda(z)$ and 
$\psi_{\lambda+\rho_{\underline{k}^\prime}}(z)$ are in the common
eigenspace $\overline{M}_\lambda$, hence 
$\psi_{\lambda+\rho_{\underline{k}^\prime}}(z)=
K_0^{BA}(\lambda+\rho_{\underline{k}^\prime})\Phi_\lambda(z)$
by Theorem \ref{HCspecial}.
\end{proof}
%%%%%%%%%%%%%%%%%%%%%%%%%%%%%%%%%%%%%%%%%%%%%%%%%%%%%%%%%%%%%%%%%%%%%%%%

%%%%%%%%%%%%%%%%%%%%%%%%%%%%%%%%%%%%%%%%%%%%%%%%%%%%%%%%%%%%%%%%%%%%%%%%%
\begin{rema}
{\bf i)} In the trigonometric differential degeneration, the analogue of
Proposition \ref{Baconn} was established in
\cite{C00} and \cite[Section VI.C]{C0}.\\
{\bf ii)}
Proposition \ref{Baconn} suggests that various
properties of the normalized Baker-Akhiezer functions (such as
duality \cite[Thm.\,4.7]{C} and bispectrality \cite[Cor.\,4.8]{C})
can be transferred to Harish-Chandra series for arbitrary
multiplicity labels $\underline{k}$, cf. \cite[\S 6]{C0} 
for the differential set-up.
We return to these issues in future work.
\end{rema}
%%%%%%%%%%%%%%%%%%%%%%%%%%%%%%%%%%%%%%%%%%%%%%%%%%%%%%%%%%%%%%%%%%%%%%%%%%%%

%%%%%%%%%%%%%%%%%%%%%%%%%%%%%%%%%
%%     Appendix of notation    %%
%%%%%%%%%%%%%%%%%%%%%%%%%%%%%%%%%
\section{Appendix:  Commonly used notation}\label{sectionCN}

We provide here a list of notation used throughout the paper. For
each symbol, we give the subsection where it was first introduced and we
provide a brief description. The reader is referred to the
appropriate subsection for more information and explicit definitions.

\medskip
\begin{tabbing}
\noindent
\=Defined in Subsection \ref{Rsd}:\=\\
\>$(V,\langle\cdot,\cdot\rangle)$\>finite dimensional Euclidean
space\\
\>$\widehat{V}$\>space of affine linear real functions on $V$\\
\>$c\in \widehat{V}$\>constant function one\\
\>$D$\>the gradient map from $\widehat V$ to $V$\\
\>$f^\vee$\>co-root $2f/\|f\|^2$\\
\>$R$\>finite reduced irreducible root system contained in $V$\\
\>$s_\alpha$\>orthogonal reflection in the hyperplane $\alpha^\perp\subset V$\\
\>$W_{0}$\>Weyl group for $R$ generated by the $s_\alpha,\alpha\in R$\\
\>$Q=Q(R)$\>root lattice of $R$\\
\>$P=P(R)$\>weight lattice of $R$\\
\>$Q^\vee=Q(R^\vee)$\>co-root lattice of $R$\\
\>$P^\vee=P(R^\vee)$\>co-weight lattice of $R$\\
\>$W_{Q^\vee}$\>affine Weyl group of $R$\\
\>$W_{P^\vee}$\>extended affine Weyl group of $R$\\
\>$t(\lambda), \lambda\in P^\vee$\>translation sending $v$ to $v+\lambda$, for
$v\in V$\\
\>$S(R)$\> affine root system $\{\alpha+rc\,\, | \,\, \alpha\in
R,\,\, r\in\mathbb{Z} \}$\\
\>$S(R)^\vee$\>dual affine root system $\{f^\vee \,\, | \,\, f\in
S(R)\}$\\
\>$S_{nr}$\>roots of the nonreduced affine root system of type
$C^\vee C_n$\\
\>Cases {\bf a,b,c}\>  Cherednik-Macdonald theory cases \\
\>$(R,R^{\prime})$\>pair of root systems defined for each case \\
\>$(L,L')$\>pair of lattices defined for each case\\
\>$(S,S')$\>pair of irreducible affine root systems defined for each
case\\
\>$(W,W')$\>extended Weyl groups associated to $(R,R^\prime)$\\
\>$S_s$\>set of short roots of $S$\\
\>$S_l$\>set of long roots of $S$\\
\>$\mathcal{O}_i$, $i=1,\cdots,5$\>$W$-orbits of $S$ for case {\bf c}\\
\>$S_1$\>reduced affine root subsystem of indivisible affine roots in
$S$\\
\>$\{\alpha_1,\dots, \alpha_n\}$\>basis for the root system $R$\\
\>$\varphi$\>highest root of $R$ with respect to above basis\\
\>$\{a_0,a_1,\dots, a_n\}$\>basis for $S$ (given in this section)\\
\>$\Delta$\>$\{a_1,\dots, a_n\}$\\
\>$R^{\prime\vee}$\>$D(S_1)$\\
\>$s_i$\>reflection associated to the simple root $a_i$\\
\>$S_1^+$\>positive affine roots of $S_1$ with respect to $\{a_0,\dots,
a_n\}$\\
\>$S_1^-$\>negative affine roots of $S_1$ with respect to $\{a_0,\dots,
a_n\}$\\
\>$l(w), w\in W$\>$\#\bigl(S_1^+\cap w^{-1}S_1^-\bigr)$\\
\>$\Omega$\>$\{w\in W\, |\, l(w)=0\}$\\
\>$\mathbb{C}[W]$\>complex group algebra of $W$\\
\>$\mathbb{C}[\Omega]$\>complex group algebra of $\Omega$\\
\>\>\\
\=Defined in Subsection \ref{DRO}:\=\\
\>$A$\>group algebra $\mathbb{C}[L]$ with basis
$\{z^{\lambda}|\ \lambda\in L\}$\\
\>$A^\prime$\>group algebra $\mathbb{C}[L^\prime]$ with basis
$\{\xi^{\lambda^{\prime}}|\ \lambda^{\prime}\in L^\prime\}$\\
\>$V_{\mathbb{C}}$\>complexification of $V$\\
\>$\mathcal{Q}$\>quotient field of $A$\\
\>$\mathcal{Q}^\prime$\>quotient field of $A^\prime$\\
\>$\mathcal{R}$\>subalgebra of $\mathcal{Q}$ defined by Definition
\ref{Rdef}\\
\>$\mathbb{D}_{\mathcal{R}}(L^\prime)$\>smash-product algebra $\mathcal{R}\#
t(L^\prime)$\\
\>$\mathbb{D}_{\mathcal{R}}(W)$\>smash-product algebra $\mathcal{R}\#
W$\\
\>$X$\>lattice in $V$\\
\>$\mathbb{D}_{\mathcal{R}}(X)$\>smash-product algebra $\mathcal{R}\# t(X)$\\
\>$\mathbb{C}[[z^{-\Delta}]]$\>algebra of formal power series \eqref{f}\\
\>$\mathbb{D}_{\mathbb{C}[[z^{-\Delta}]]}(X)$\>smash-product algebra
$\mathbb{C}[[z^{-\Delta}]]\# t(X)$ for lattices $X\subset V$\\
\>$M$\>algebra of analytic functions on $V_{\mathbb{C}}$ spanned
by $z^u, u\in V_{\mathbb{C}}$\\
\>$\overline{M}$\>vector space spanned by elements of
the form \eqref{leadingexponent}\\
\>\>\\
\=Defined in Subsection \ref{Chers}:\=\\
\>$\mathbb{C}(S)^W$\>space of multiplicity labels associated to
$S$\\
\>$\underline{k}$\>multiplicity label associated to $S$\\
\>$k_a$\>$\underline{k}(a)$\\
\>$k_i$\>$\underline{k}(a_i)$\\
\>$\underline{k}^\prime$\>multiplicity label associated to $S^\prime$
and dual to $\underline{k}$\\
\>$\kappa_j$(resp. $\kappa_j^{\prime}$)\>value of $\underline{k}$
(resp. $\underline{k}^\prime$) at the $W$-orbit
$\mathcal{O}_j$ for case {\bf c}\\
\>$\mathbb{C}(S_1)^W$\>space of multiplicity labels associated to
$S_1$\\
\>$\underline{\tau}$, $\underline{\tau}^\prime$\>invertible multiplicity
labels of $S_1$ associated to $\underline{k}$\\
\>$H(\underline{\tau})$\>extended affine Hecke algebra (Definition
\ref{HA})\\
\>$T_i, i=0,\dots, n$\> generators for $H(\underline{\tau})$\\
\>$T_w$\>$\omega T_{i_1}T_{i_2}\cdots T_{i_{l(w)}}$ for reduced
expression $w=\omega s_{i_1}s_{i_2}\cdots s_{i_{l(w)}}$\\
\>$c_a(z)$\>element in $\mathcal{R}$ given by \eqref{c}\\
\>$\pi_{\underline{k}}$\>embedding of $H(\underline{\tau})$ inside
$\mathbb{D}_{\mathcal{R}}(W)$\\
\>$T_i(\underline{k})$\>difference reflection operator
$\pi_{\underline{k}}(T_i)$\\
\>$\beta$\>linear map given by $\beta(\sum_{w\in W_0}D_ww)=\sum_{w\in
W}D_w$\\
\>$H_0=H_0(\underline{\tau})$\>finite Hecke algebra generated by
$T_j, j=1,\dots, n$\\
\>$Z(H(\underline{\tau}))$\>center of $H(\underline{\tau})$\\
\>\>\\
\=Defined in Subsection \ref{CMsub}:\=\\
\>$V_+$\>open dominant Weyl chamber in $V$ with respect to $R^+$\\
\>$\overline{V}_+$\>closure of $V_+$ in $V$\\
\>$L_{++}$\>$L\cap \overline{V}_+$\\
\>$L^\prime_{++}$\>$L^\prime\cap \overline{V}_+$\\
\>$Y^{\lambda^\prime}$, $\lambda^\prime\in L_{++}^\prime$\>
$T_{t(\lambda^\prime)}$\\
\>$Y^{\lambda^{\prime}}$, $\lambda^\prime\in L^\prime$\>
$Y^{\mu^\prime}\bigl(Y^{\nu^\prime}\bigr)^{-1}$  for
$\lambda^{\prime}=\mu^{\prime}-\nu^{\prime}$ with
$\mu^{\prime},\nu^{\prime}\in L_{++}^{\prime}$\\
\>$A^\prime(Y)$\>subalgebra of $H(\underline{\tau})$ spanned by the
$Y^{\lambda^{\prime}}$, $\lambda^{\prime}\in L^{\prime}$\\
\>$A_0^\prime$\>algebra of $W_0$-invariant elements in
$A^{\prime}$\\
\>$A_0^\prime(Y)$\>algebra corresponding
to $A_0^\prime$ via canonical isomorphism $A^\prime\cong A^\prime(Y)$\\
\>$m_{\lambda'}(\xi)$, $\lambda^\prime\in L^\prime$\>monomial symmetric 
function $\sum_{\mu^\prime\in
W_0\lambda^\prime}\xi^{\mu^\prime}$ in $A_0^\prime$\\
\>$m_{\lambda'}(Y)$\>element in $A_0^\prime(Y)$ corresponding to
$m_{\lambda'}(\xi)$\\
\>$D_p, D_{\lambda^\prime}$\>$W_0$-invariant difference operators
defined in Definition \ref{Dpdef}\\
\>$\rho_{\underline{k}^\prime}$\>deformed half sum of positive roots
(see \eqref{rho})\\
\>$\widetilde{p}(\xi)$, 
$p(\xi)\in A^\prime$\>$\rho_{\underline{k}^\prime}$-twisted
Laurent polynomial (see \eqref{rhotwist})\\
\>minuscule co-weight\>$\pi^\prime\in P(R^\vee)$ such that
$|\langle\pi^\prime,\alpha\rangle|\leq 1$ for all $\alpha\in R$\\
\>quasi-minuscule\>$\pi^{\prime}\in R^\vee$ such that 
$|\langle\pi^\prime,\alpha\rangle|\leq 1$ for all
$\alpha\in R\setminus \{\pm \pi^{\prime\vee}\}$\\
\>$D_{\pi^\prime}$\>Macdonald difference operator (Definition
\ref{Macdonalddodef})\\
\>$\pi^{\prime}_i$\>fundamental co-weights\\
\>$w_0$\>longest Weyl group element in $W_0$\\
\>$J_0$\>$\{j\in\{1,\ldots,n\} \,\, | \,\, \pi_j^\prime\in
L^\prime\,\,\hbox{ and }\,\, m_j=1\}$\\
\>$W_{0,\lambda^\prime}$\>isotropy subgroup of  $\lambda^\prime$
inside $W_0$\\
\>$W_0^{\lambda^\prime}$\>complete set of representatives of
$W_0/W_{0,\lambda^\prime}$\\
\>$S_1(w)$\>$S_1^+\cap w^{-1}S_1^-$\\
\>$m_{\lambda^\prime}(-\rho_{\underline{k}^\prime})$\>$\sum_{\mu^\prime\in
W_0\lambda^\prime}q^{\langle
-\rho_{\underline{k}^\prime},\mu^\prime\rangle}$\\
\>$f_{\pi^\prime}(z)$\>coefficients of Macdonald difference
operators (Proposition \ref{formexplicit})\\
\>\>\\
\=Defined in Subsection \ref{HChom}:\=\\
\>$V_-$\>$\{ v\in V \,\, | \,\, \langle v,\alpha\rangle<0 \qquad
\forall\,\alpha\in R^+ \}$\\
\>$\gamma(\underline{k})$\>Harish-Chandra homomorphism \eqref{gammak}\\
\>\>\\
\=Defined in Subsection \ref{Symmsub}:\=\\
\>$\pi(D)$, $D\in\mathbb{D}_{\mathcal{R}}(L^\prime)$\>
$({\#W_0})^{-1}\sum_{w\in W_0}wDw^{-1}$
\\
\>$e$\>trivial idempotent in $\mathbb{C}[W_0]$\\
\>\>\\
\=Defined in Subsection \ref{Rrsub}:\=\\
\>$V_{-,F}$\>$\{v\in V \, | \, \langle v,\alpha\rangle=0\quad
(\alpha\in F),\, \langle v,\beta\rangle<0\quad
(\beta\in\Delta\setminus F) \}$\\
\>$R^{\prime\vee}_F$\>$R^{\prime\vee}\cap \mathbb{Z}F$\\
\>$R^{\prime\vee,+}_F$ (resp. $R^{\prime\vee,-}_F$) \>
set of positive (resp. negative) roots in $R^{\prime\vee}_F$\\
\>$W_{0,F}$\>parabolic subgroup of $W_0$ generated by
$s_{\alpha},\alpha\in F$\\
\>$R_F$\>$R\cap\bigoplus_{\alpha\in F}\mathbb{Q}\alpha$\\
\>$R^{\pm}_F$\>$R^{\pm}\cap\bigoplus_{\alpha\in F}\mathbb{Q}\alpha$\\
\>$\rho_{\underline{k}^\prime,F}$\>
$\frac{1}{2}\sum_{\alpha\in R^+\setminus R_F^+}\underline{k}^\prime(\alpha^\vee)
\alpha$\\
\>$S_F$\>$\{a\in S \,\, | \,\, Da\in \mathbb{Z}F\}$\\
\>${\mathcal R}_F$\>algebra generated by $(1-rz^{\alpha})^{-1}$, $r\in
\mathbb{C}, \alpha\in R^{\prime\vee}_F$\\
\>$\mathbb{C}[[z^{-F}]]$\>$\mathbb{C}[[z^{-\alpha}|\alpha\in F]]$\\
\>$\mathbb{D}_{\mathcal{R}_F}(X)$\>$\mathcal{R}_F\# t(X)$\\
\>$\mathbb{D}_{\mathbb{C}[[z^{-F}]]}(X)$\>$\mathbb{C}[[z^{-F}]]\#
t(X)$\\
\>$\gamma_F(\underline{k})$\>constant term map along $V_{-,F}$
(Definition \ref{gammaFdef})\\
\>$L_i^\prime$\>$L^\prime+\mathbb{Z}\frac{\alpha_i^\vee}{2}$ in cases
{\bf a}, {\bf b}; $L^\prime$ in case {\bf c}\\
\>$X_i$\>$\{\mu^\prime\in L_i^\prime \,\, | \,\,
\langle\mu^\prime,\alpha_i\rangle=0\}$\\
\>$\mathcal{L}_i$\>rank one difference operator defined in Definition
\ref{rankoneL}\\
\>$y_i, z_i$\>elements of $\mathbb{C}[X_i]$ defined by \eqref{yz}\\
\>$\mathbb{D}_i^{\pi^\prime}(\underline{k})$\>centralizer of
$\gamma_{\{a_i\}}(\underline{k})(D_{\pi^\prime})$
in $\mathbb{D}_{\mathcal{R}_{\{a_i\}}}(L_i^\prime)^{W_{0,\{a_i\}}}$\\
\>\>\\
\=Defined in Subsection \ref{HCformalsub}:\=\\
\> $\Phi(z,\xi)$\>difference analogue of the Harish-Chandra series
(see \eqref{seriesform})\\
\>\>\\
\=Defined in Subsection \ref{HCspecialsub}:\=\\
\>$\sigma$\>positive real number such that $q=e^{-2\pi\sigma}$\\
\>$Z$\>$\{\lambda\in V \, | \, \langle
\lambda^\prime,\lambda\rangle\in\mathbb{Z}\quad\forall\,
\lambda^\prime\in L^\prime\}$\\
\>$w\cdot
\lambda$\>$w(\lambda+\rho_{\underline{k}^{\prime}})-\rho_{\underline{k}^{\prime}}$\\
\>$\mathcal{D}_{\underline{k}}$\>$\bigcup_{w\in W_0}\{\lambda\in V_{\mathbb{C}} \, 
| \, \lambda-w\cdot\lambda\in \mathbb{Z}_{\geq 0}\Delta\setminus
\{0\}+\sqrt{-1}Z/\sigma\}$\\
\>$\mathcal{B}^\prime_{\underline{k}}$\>elements of $\mathcal{Q}^\prime$
with singularities  in $\mathcal{D}_{\underline{k}}$\\
\> $\iota_\lambda$, $\lambda\in V_{\mathbb{C}}\setminus 
\mathcal{D}_{\underline{k}}$
\> specialization map (Lemma \ref{iotalemma})\\
\>
$\Phi_\lambda(z)$\>$\iota_\lambda\bigl(\Phi(z,\xi)\bigr)\in \overline{M}$
(see Theorem \ref{HCspecial})\\
\>$\overline{M}_\lambda$\>$\{F(z)\in\overline{M} \,\,\, | \,\,\,
D_p\bigl(F(z)\bigr)=\widetilde{p}(\lambda)F(z)\quad\forall\,p(\xi)\in
A_0^\prime\}$\\
\>\>\\
\=Defined in Subsection \ref{polynomialsection}:\=\\
\>$P_{\lambda}(z)$\>Macdonald polynomial of degree $\lambda\in
L_{++}$ (Definition \ref{Macdonalddef})\\
\>\>\\
\=Defined in Subsection \ref{BAsub}:\=\\
\>$\rho_{\underline{k}}^\prime$\>$\frac{1}{2}\sum_{\alpha\in
R^+}k_\alpha\alpha^\vee$\\
\>$\psi_{\lambda}(z)$\>Baker-Akhiezer function (Definition
\ref{Bafunction})\\

\end{tabbing}

%%%%%%%%%%%%%%%%%%%%%%%%%%%%%%%%%%%%%%%
%%                                   %%
%%      References                   %%
%%                                   %%
%%%%%%%%%%%%%%%%%%%%%%%%%%%%%%%%%%%%%%%

\end{document}